\documentclass[11pt]{article}
\usepackage{amsmath, amssymb, color}
\usepackage{graphicx,graphics}
\usepackage{tikz}
\usepackage{pgflibraryarrows}
\usepackage{pgflibrarysnakes}
\usepackage{mathrsfs}
\usepackage[all,poly,knot]{xy}
\allowdisplaybreaks[3]

\date{}

\newtheorem{Theorem}{Theorem}[section]
\newtheorem{Lemma}[Theorem]{Lemma}
\newtheorem{Cor}[Theorem]{Corollary}
\newtheorem{Prop}[Theorem]{Proposition}
\newcommand{\QED}{\ \hfill \rule{0.5em}{0.5em} }
\numberwithin{equation}{section}

\begin{document}

\title{Graph Inverse Semigroups and Leavitt Path Algebras}

\author{ John Meakin and Zhengpan Wang\thanks{Partially supported by Chongqing Natural Science
Foundation (cstc2019jcyj-msxmX0435).}  }

\maketitle

\begin{abstract}

We study two classes of inverse semigroups built from directed graphs, namely graph inverse semigroups and a new class of  semigroups that we refer to as Leavitt inverse semigroups.
These semigroups are closely related to graph $C^*$-algebras
and Leavitt path algebras.
We  provide a topological
characterization of the universal groups of the local submonoids of
these inverse semigroups. We  study the relationship between the
graph inverse semigroups of two graphs when there is a directed
immersion between the graphs. We describe  the structure of graphs
that admit a directed cover or directed immersion into a circle and
we provide structural information about graph inverse semigroups of
finite graphs that admit a directed cover onto a bouquet of circles.
We also find necessary and sufficient conditions for a homomorphic
image of a graph inverse semigroup to be another graph inverse
semigroup. We find a
presentation for the Leavitt inverse semigroup of a graph  in terms of
generators and relations.  We describe the structure of the Leavitt
inverse semigroup and the Leavitt path algebra of a graph that
admits a directed immersion into a circle. We show that two graphs
that have isomorphic Leavitt inverse semigroups have isomorphic
Leavitt path algebras  and we classify graphs that have isomorphic
Leavitt inverse semigroups. As a consequence, we show that Leavitt path algebras are $0$-retracts of certain matrix algebras.
\renewcommand{\thefootnote}{}
\footnotetext{{\bf Mathematics Subject Classification (2010)}  20M18, 16S99}
\footnotetext{{\bf Keywords:} graph inverse semigroups; Leavitt inverse semigroups; Leavitt path algebras; directed immersions; directed covers; universal groups; congruence pairs; 0-retracts}
\end{abstract}

\section {Introduction}

The notion of a Leavitt path algebra is an outgrowth of the seminal
paper by W.G. Leavitt \cite{L} providing a construction of what is
now referred to as the {\em Leavitt algebra} $L_F(1,n)$
corresponding to a positive integer $n$ and a field $F$. The
algebras $L_F(1,n)$ are the universal examples of algebras that do
not have the invariant basis number property, namely if $R = L_F(1,n)$ then
the free left $R$-modules $R$ and $R^n$ are isomorphic.

If $F$ is a field and $\Gamma$ is a directed graph, then we may form
the Leavitt path algebra $L_F(\Gamma)$, whose elements correspond
roughly to directed paths in the graph. The precise definition is
given in Section 6 below. Leavitt path algebras for $F = \mathbb C$ are closely related to Cuntz-Krieger graph $C^*$-algebras in the sense of Kumjian, Pask and Raeburn  \cite{KPR}. The Leavitt algebra $L_F(1,n)$ is the
Leavitt path algebra constructed from the graph $\Gamma = B_n$, the
bouquet of $n$ circles (that is, the graph with one vertex and $n$
directed edges). The general study of Leavitt path algebras was initiated independently by Abrams and Aranda Pino \cite{AP}
and by Ara, Moreno and Pardo \cite{AMP} around 2004. It has deep connections to ring theory and the theory
of graph $C^*$-algebras. We refer the reader to the survey article by Abrams \cite{A} or the book by  Abrams,
Ara and
Siles Molina \cite{AAM} for much information
about Leavitt path algebras.

A certain amount of structural information about Leavitt path
algebras may be gleaned from the theory of inverse semigroups. We
recall that an {\em inverse semigroup} is a semigroup $S$ such that
for every $a \in S$ there exists a unique element $a^{-1} \in S$
such that $a = aa^{-1}a$ and $a^{-1} = a^{-1}aa^{-1}$. The book by
Lawson \cite{Law1} is a standard reference for the theory of inverse
semigroups and their connections to other fields: any undefined
notation and concepts about inverse semigroups that are used in this
paper may be found in \cite{Law1}.  In particular, we shall make use
(often without comment) of the elementary fact that idempotents of
an inverse semigroup commute. We shall also make use of the natural
partial order on an inverse semigroup $S$ (defined by $a \leq b$ if
$a = eb$ for some idempotent $e$ of $S$).

Most of the inverse semigroups that arise in this paper (except groups!) have a zero, which we denote by $0$ (or $0_S$ if we need to specify the inverse semigroup $S$ under consideration) and all homomorphisms under consideration will map $0$ to $0$. Thus, {\em unless stated otherwise, an inverse semigroup $S$ has a zero, and a homomorphism $f: S \rightarrow T$ between inverse semigroups $S$ and $T$ will be assumed to map $0_S$ onto $0_T$}.

The most obvious way to associate an inverse semigroup to a Leavitt
path algebra is to study the connection between Leavitt path
algebras and graph inverse semigroups. The graph inverse semigroup
$I(\Gamma)$ associated with a directed graph $\Gamma$ is defined in  Section 2
below.  Leavitt path algebras may be viewed as algebras constructed
from the contracted semigroup algebra of a graph inverse semigroup
by imposing some additional algebra relations known as the
Cuntz-Krieger relations. It is  known (see \cite{MeMi}, Theorem 20)
that if two graph inverse semigroups $I(\Gamma)$ and $I(\Delta)$ are
isomorphic, then the corresponding graphs $\Gamma$ and $\Delta$ are
isomorphic, and hence the Leavitt path algebras $L_F(\Gamma)$ and
$L_F(\Delta)$ are isomorphic, but the converse is far from true.

In the present paper we study several structural properties of graph
inverse semigroups. We also introduce another inverse semigroup
$LI(\Gamma)$ naturally associated with a directed graph $\Gamma$.
The inverse semigroup $LI(\Gamma)$ is the multiplicative subsemigroup of the
Leavitt path algebra $L_F(\Gamma)$ generated by the vertices and
edges (and ``inverse edges") of the graph $\Gamma$ (these elements
generate $L_F(\Gamma)$ as an $F$-algebra). It is a quotient of the
graph inverse semigroup $I(\Gamma)$ and again has the property that
$L_F(\Gamma) \cong L_F(\Delta)$ if $LI(\Gamma) \cong LI(\Delta)$.
While the converse is certainly false in general, these inverse
semigroups
 provide significantly more information about Leavitt path
algebras than do graph inverse semigroups.

The definition and basic notation for graph inverse semigroups is
introduced in Section 2 of this paper. In Section 3  we study the
relationship between the graph inverse semigroups
$I(\tilde{\Gamma})$ and $I(\Gamma)$ when there is a directed cover
or directed immersion $f : \tilde{\Gamma} \rightarrow \Gamma$. In
this case the map $f$ induces homomorphisms between corresponding
local submonoids of the graph inverse semigroups (Theorem
\ref{circuitembed}). We provide a description of graphs that admit a
directed cover or directed immersion into a circle (Theorem
\ref{out1}) and we prove a structural property of finite directed
covers of a bouquet of circles (Theorem \ref{finitecover}). Section
4 is concerned with  groups naturally  associated with graph inverse
semigroups. We examine the universal group of a graph inverse
semigroup and provide a topological description of the universal
groups of its local submonoids (Theorems \ref{local} and
\ref{freefact}). In Section 5 we determine necessary and sufficient
conditions for a quotient of a graph inverse semigroup $I(\Gamma)$
to be another graph inverse semigroup (Theorem \ref{preserve}) and
as a consequence we show that the quotient graph inverse semigroup
is a retract of $I(\Gamma)$ (Corollary \ref{congraphinv}).

Section 6 is concerned with Leavitt path algebras and Leavitt
inverse semigroups. We  define the Leavitt inverse semigroup
$LI(\Gamma)$ associated with a directed graph $\Gamma$ and find a
presentation for $LI(\Gamma)$ as an inverse semigroup in terms of
generators and relations (Theorem \ref{leavittinv}). We describe the
structure of the Leavitt inverse semigroup and the Leavitt path
algebra of a graph that admits a directed immersion into a circle
(Theorems \ref{classout1} and \ref{classout2}). We show that two
graphs that have isomorphic Leavitt inverse semigroups have
isomorphic Leavitt path algebras (Theorem \ref{algiso}). In the
final section (Section 7)  we classify graphs that have isomorphic
Leavitt inverse semigroups (Theorem \ref{leavinvclass}). As  a consequence, we obtain  structural results for Leavitt path algebras of a restricted class of graphs and we   show that Leavitt path algebras are $0$-retracts of  matrix algebras of a restricted type (Theorem \ref{algstructure}).


\section{Graph inverse semigroups}

All graphs under consideration in this paper will be {\em directed}
graphs with either finitely many or countably infinitely many
vertices and edges. We denote the set of vertices of a graph
$\Gamma$ by ${\Gamma}^0$ and the set of edges of $\Gamma$ by
${\Gamma}^1$. If $e \in \Gamma^1$ then $e$ is a directed edge from a
vertex that we will denote by $s(e)$ to a vertex that we will denote
by $r(e)$. In fact, $s$ and $r$ can be considered as mappings of
$\Gamma^1$ into $\Gamma^0$, respectively called the {\em source
mapping} and the {\em range mapping} for $\Gamma$. Thus for each vertex
$v \in \Gamma^0, s^{-1}(v) = \{e \in \Gamma^1 : s(e) = v\}$ and  the
{\em out-degree} of a vertex $v$ is $|s^{-1}(v)|$, the number of
directed edges with source $v$. (This is referred to as the {\em
index} of $v$ by some  authors). We allow for the possibility that
$s(e) = r(e) = v \in \Gamma^0$, in which case $e$ is a {\em loop} at
$v$. A {\em directed path} in $\Gamma$ is a finite sequence $p =
e_1e_2\ldots e_n$ of edges $e_i \in \Gamma^1$ with $r(e_i) = s(e_{i+1})$
for $i = 1,\cdots ,n-1$. We define $s(p) = s(e_1)$ and $r(p) = r(e_n)$
and refer to $p$ as a directed path from $s(p)$ to $r(p)$. We also
consider a vertex $v$ as being an {\em empty (directed) path} (i.e. a path with
no edges) based at $v$ and with $s(v) = r(v) = v$.

It is convenient to extend the notation so as to allow
 paths in which edges are read in either the positive or
negative direction. To do this, we associate with each edge $e$ an
``inverse edge" $e^*$ (sometimes called a ``ghost edge" by some
authors)
 with $s(e^*) = r(e)$ and $r(e^*) = s(e)$. Also define $(e^*)^*
= e$. We denote by $(\Gamma^1)^*$ the set $\{e^* : e \in \Gamma^1\}$
and assume that $\Gamma^1 \cap (\Gamma^1)^* = \emptyset$ and that the map $e \rightarrow e^*$ is a bijection from $\Gamma^1$ to $(\Gamma^1)^*$. With this
convention, we can define a {\em path} in $\Gamma$ as a sequence $p
= e_1e_2\ldots e_n$ with $e_i \in \Gamma^1 \cup (\Gamma^1)^*$ and
$r(e_i) = s(e_{i+1})$ for $i = 1,\cdots ,n-1$ and for each path $p =
e_1e_2\ldots e_n$ we define the inverse path to be $p^* =
e_n^*\ldots e_2^*e_1^*$. As usual, $s(p) = s(e_1)$ and $r(p) = r(p_n)$.
 The graph $\Gamma$ is said to be {\em
connected} if for all $v,w \in \Gamma^0$ there is at least one path
$p$ with $s(p) = v$ and $r(p) = w$ while $\Gamma$ is said to be {\em
strongly connected} if for all $v,w \in \Gamma^0$ there is at least
one {\em directed} path $p$ with $s(p) = v$ and $r(p) = w$. A path
$p$ is a {\em circuit at $v$} if $s(p) = r(p) = v$. Thus, for
example, a path of the form $ee^*$ where $e$ is an edge with $s(e) =
v$ is a circuit at $v$. A path $p = e_1e_2\ldots e_n$ is called {\em
reduced} if $e_i \neq e_{i+1}^*$ for each $i$. A {\em reduced circuit} is a circuit
$p = e_1e_2\ldots e_n$ that is a reduced path and such that $e_1 \neq e_n^*$.   A {\em directed
circuit} is a directed path that is a circuit. A {\em cycle} is a
directed circuit $e_1e_2\ldots e_n$ such that $s(e_i) \neq s(e_j)$ if
$i, j \in \{1,2,\cdots ,n\}$ and $i \neq j$. Two cycles $C_1$ and $C_2$
are said to be {\em conjugate} if $C_1 = e_1e_2\ldots e_n$ and $C_2 =
e_ie_{i+1}\ldots e_ne_1\ldots e_{i-1}$ for some $i$.
The graph $\Gamma$ is {\em acyclic} if it has no non-trivial cycles.  $\Gamma$
is called a {\em tree} if it is connected and has no non-trivial
reduced circuits. Equivalently (see for example Hatcher's book
\cite{hatch}), $\Gamma$ is a tree if it is connected and its
fundamental group $\pi_1(\Gamma)$ is trivial. Thus trees are
connected acyclic graphs but connected acyclic graphs are not
necessarily trees.

\medskip

Graph inverse semigroups were first introduced by Ash and Hall
\cite{AH} (for a restricted class of directed graphs) in connection
with their study of the partially ordered set of
$\mathcal{J}$-classes of a semigroup. Graph inverse semigroups generalize the polycyclic
monoids introduced by Nivat and Perrot \cite{NivPer} and arise very
naturally in the extensive theories of graph $C^*$-algebras and Leavitt path algebras.
Graph inverse semigroups
have been studied in their own right by several authors, for example
Costa and Steinberg \cite{CS}, Jones and Lawson \cite{JL}, Lawson
\cite{Law2}, Krieger \cite{K},  Mesyan and Mitchell \cite{MeMi}  and
Wang \cite{Wang}.

\medskip

Define the graph inverse semigroup $I({\Gamma})$ of a directed graph
$\Gamma$ as the semigroup generated by ${\Gamma}^{0} \cup {\Gamma}^1
\cup (\Gamma^1)^*$ together with a zero $0$ subject to the relations

\medskip

(1) $s(e)e = er(e) = e$ for all $e \in \Gamma^0 \cup {\Gamma}^{1}
\cup (\Gamma^1)^*$;

(2) $u v = 0$ if $u,v \in {\Gamma}^{0}$ and $u \neq v$;

(3) $e^{*}  f = 0$ if $e,f \in {\Gamma}^1$ and $e \neq f$;

(4) $e^{*}  e = r(e)$ if $e \in {\Gamma}^1$.

\medskip

We emphasize that condition (1) of the definition above implies that
$v^2 = v$ for all $v \in \Gamma^0$; that is, the vertices of
$\Gamma$ are idempotents in $I(\Gamma)$. Condition (1) also implies
that $e^*s(e) = r(e)e^* = e^*$ for all $e \in \Gamma^1$.

 It
is  not difficult to see that $I(\Gamma)$ is in fact an {\em inverse
semigroup}. A  straightforward argument shows that every non-zero
element of a
 graph inverse semigroup $I(\Gamma)$
may be uniquely written in the form $pq^*$ where $p$ and $q$ are
directed (possibly empty) paths with $r(p) = r(q)$. We refer to this
as the {\em canonical form} of a non-zero element of $I(\Gamma)$.
The inverse of an element $pq^*$ is of course $qp^*$. If $pq^*$ and
$rs^*$ are non-zero elements of $I(\Gamma)$, then the product
$pq^*rs^*$ is non-zero if and only if either $q$ is a prefix of $r$
(i.e. $r = qt$ for some directed (possibly empty) path $t$, in which case
$pq^*rs^* = pts^*$), or else $r$ is a prefix of $q$ (i.e. $q = rt$
for some directed (possibly empty) path $t$, in which case
$pq^*rs^* = p(st)^*$). The non-zero idempotents are of the form
$pp^*$ for some (possibly empty) directed path $p$, and $pp^* \geq
qq^*$ in the natural partial order on $I(\Gamma)$ if and only if $q
= pt$ for some directed (possibly empty) path $t$. Thus the vertices of $\Gamma$ are
the maximal idempotents in the natural partial order on $I(\Gamma)$.

If $\Gamma$ is the graph with one vertex and one edge, then
$I({\Gamma})$ is the bicyclic monoid with a (removable) zero. If
$\Gamma = B_n$ is the bouquet of $ n > 1$ circles (i.e. the graph
with one vertex and $n$ directed edges), then $I({\Gamma})$ is
isomorphic to the polyclic monoid $P_n$. This is the inverse monoid
generated (as an inverse monoid) by a set $A$ with $|A| = n$ subject
to the defining relations  $ a^{-1}a = 1 $ and $ a^{-1}b = 0$ for
all $ a,b \in A$ with $ b \neq a$. (Here we regard $A$ as a set of
labels for the directed edges of $B_n$). The monoid $P_n$ was
introduced by Nivat and Perrot \cite{NivPer} as the syntactic monoid
of the ``correct parenthesis" language with $n$ sets of parentheses.
It was rediscovered in the operator algebra literature where it is
referred to as the Cuntz monoid, used in the construction of the
Cuntz algebra ${\mathcal O}_{n}$ (see Paterson's book \cite{Pat} for
details). The algebra constructed from the graph $B_n$ in the
original paper by Leavitt \cite{L} is what is now referred to as the
Leavitt path algebra of this graph (see \cite{AAM} for details).

\section {Directed Covers and Directed Immersions}

A {\em morphism} from the (directed) graph $\tilde{\Gamma}$ to the
(directed)  graph $\Gamma$ is a function $f : \tilde{\Gamma}
\rightarrow \Gamma$ that takes vertices to vertices and edges to
edges, and preserves incidence and orientation of edges; that is, $f(s(\tilde{e})) = s(f(\tilde{e}))$ and $f(r(\tilde{e})) = r(f(\tilde{e}))$ for all $\tilde{e} \in \tilde{\Gamma}^1$. (Here we
abuse notation slightly by using the same symbol $f$ to denote the
corresponding function that takes vertices to vertices and the
function that takes edges to edges.)  We extend the notation by
defining $f(\tilde{e}^*) = f(\tilde{e})^*$ for all $\tilde{e} \in
\tilde{\Gamma}^1$ and $f(\tilde{e}_1\tilde{e}_2\ldots \tilde{e}_n) =
f(\tilde{e}_1)f(\tilde{e}_2)\ldots f(\tilde{e}_n)$ for each path
$\tilde{p} = \tilde{e}_1\tilde{e}_2\ldots \tilde{e}_n$. In fact we will
often use the notation  $f(\tilde{p}) = p$ to denote the image of a
path $\tilde{p} = \tilde{e}_1\tilde{e}_2\ldots \tilde{e}_n$ in
$\tilde{\Gamma}$, where it is understood that $p$ is the path $p =
e_1e_2\ldots e_n$ in $\Gamma$ and $e_i = f(\tilde{e}_i)$.

 A morphism $f :
\tilde{\Gamma} \rightarrow \Gamma$ between directed graphs induces
maps $f_{\tilde{v}} : s^{-1}(\tilde{v}) \rightarrow
s^{-1}(f(\tilde{v}))$   in the obvious
way. We say that  $f$ is a {\em directed cover} if the induced maps
  $f_{\tilde{v}}$
 are bijections for each $\tilde{v} \in
\tilde{\Gamma}^0$ and  that $f$ is a {\it directed immersion} if the
induced maps  $f_{\tilde{v}}$  are injections for each $\tilde{v}
\in \tilde{\Gamma}^0$.

This is closely related to the classical notion of covers and
immersions of graphs in Stallings' paper \cite{Stall},
the
distinction being that Stallings defines
$f$ to be a {\em cover} if the induced maps
  $f_{\tilde{v}} : s^{-1}(\tilde{v}) \cup r^{-1}(\tilde{v}) \rightarrow
s^{-1}(f(\tilde{v})) \cup r^{-1}(f(\tilde{v}))$
 are bijections for each $\tilde{v} \in
\tilde{\Gamma}^0$ and   $f$ is an {\em immersion} if these induced
maps $f_{\tilde{v}}$  are injections for each $\tilde{v} \in
\tilde{\Gamma}^0$.

There is a significant difference between directed immersions (or
directed covers) of graphs and immersions (or covers) of graphs in
the classical sense. Connected covers of  a connected graph $\Gamma$
are classified via conjugacy classes of subgroups of the fundamental
group $\pi_1(\Gamma)$ of the graph (see  \cite{hatch} or
\cite{Stall}). For example, connected covers of the circle
$B_{\{a\}}$ (the graph with one vertex and one directed edge) are
classified via subgroups of $\mathbb Z$. This yields the circuits
$C_n$ with $n$ edges (the finite covers of $B_{\{a\}}$) and the
universal cover of $B_{\{a\}}$ (the Cayley graph of $\mathbb Z$
relative to the usual presentation $\mathbb Z = Gp \langle a :
\emptyset\rangle$). The only connected immersions into $B_{\{a\}}$
are the connected covers and the connected subgraphs of the
universal cover. However the description of  {\em directed} covers
of $B_{\{a\}}$ and {\em directed} immersions into $B_{\{a\}}$ is
more complicated. Let $L_{\infty}$ be the linear graph with vertices
$v_{-k}, \, k \geq 0$ and an edge $e_k$ from $v_{-k}$ to $v_{-k+1}$
for $k > 0$. For each integer $n \geq 0$ let $L_n$ be the induced
subgraph of $L_{\infty}$ spanned by the vertices $v_{-k}, \, 0 \leq
k \leq n$.

\begin{Theorem}
\label{out1}

A graph $\Gamma$  admits a directed  immersion into $B_{\{a\}}$ if
and only if the out-degree of every vertex of $\Gamma$ is at most
$1$. If $\Gamma$ is a connected graph, all of whose vertices have
out-degree at most $1$, then

$($a$)$ $\Gamma$ has at most one sink. If $\Gamma$ does have a sink
$v_0$ and $v$ is any other vertex in $\Gamma$, then there is a
unique directed path from $v$ to $v_0$ and  $\Gamma$ is a tree.
$\Gamma$ has a  sink if and only if it admits a directed cover onto
the graph $L_n$ where $n$ is the maximum length of a directed path
from some vertex of $\Gamma$ to $v_0$ $($and $n = \infty$ if there are
directed paths of arbitrary length ending at $v_0$$)$.

$($b$)$ If $\Gamma$ is not a tree then  $\Gamma$ has a non-trivial cycle
  and any two non-trivial cycles are  cyclic conjugates of each other.
 Furthermore, if $v'$ is any vertex on one of these cycles $C$ and $v$ is any
 other vertex of $\Gamma$  then  there is a
unique directed path from $v$ to $v'$ that does not include the
cycle $C$ as a subpath. In this case $\Gamma$ is a directed  cover
of $B_{\{a\}}$.

$($c$)$ $\Gamma$ is a directed cover of $B_{\{a\}}$ if and only if it is
either a tree that has no sink or $\pi_1(\Gamma) \cong \mathbb Z$,
in which case $\Gamma$ has a structure as described in part $($b$)$.

\end{Theorem}

\noindent {\bf Proof.} It is clear from the definition of  a
directed immersion that if there is a directed immersion from
$\Gamma$ into $B_{\{a\}}$, then the out-degree of every vertex of
$\Gamma$ is at most $1$ and that if $\Gamma$ covers $B_{\{a\}}$,
then the out-degree of every vertex of $\Gamma$ is $1$. Conversely,
if the out-degree of every vertex is at most $1$ then the obvious
map from $\Gamma$ to $B_{\{a\}}$ is a directed immersion, which is a
directed cover if the out-degree of every vertex is $1$.

(a) Observe first that if $p = e_1e_2\ldots e_n$ is a path in $\Gamma$ such that
$e_1 \in (\Gamma^1)^*$, then we must have $e_i \in (\Gamma^1)^*$ for all $i$ since every vertex has out-degree at most $1$.  Suppose that $v_0$ and $v_1$ are  sinks of $\Gamma$.
Since $\Gamma$ is connected there is a path $p = e_1e_2\ldots e_k$ from
$v_0$ to $v_1$. But since $v_0$ and $v_1$ are both sinks we must have $e_1 \in (\Gamma^1)^*$ and
$e_k \in \Gamma^1$. This violates the observation above unless $v_0 = v_1$, so
$\Gamma$ has a unique sink  if it has one. Suppose that $\Gamma$
does have a sink $v_0$ and that $v$ is any vertex in $\Gamma$ with
$v \neq v_0$. There is a path $p$ containing no circuits from $v$ to $v_0$ that must be
directed by the  argument above. If $p'$ is another directed path
from $v$ to $v_0$, then $p'$ has no circuits since the out-degree of every vertex in $p'$ other than $v_0$ is $1$. We may write $p = p_1p_2$ and $p' = p_1p_2'$
where the first edge of $p_2$ is different from the first edge of
$p_2'$. But this yields a vertex $r(p_1)$ of degree at least $2$, a
contradiction. So there is a unique directed path from $v$ to
$v_0$. If $p = e_1e_2\ldots e_n$ is a reduced circuit such that $s(e_i)
\neq s(e_j)$ for $i \neq j$, then either $p$ or $p^*$ is a cycle
since the out-degree of every vertex in the circuit must be $1$. But the
graph $\Gamma$ cannot contain any non-trivial cycle since $v_0$ is
not a vertex of any such cycle and every vertex $v$ in a cycle must
have out-degree $1$, which is impossible since there is a directed
path from $v$ to $v_0$. It follows  that $\Gamma$ is a tree.

Suppose that $n$ is the maximum length of a directed path in
$\Gamma$ ending at $v_0$. For each vertex $v$ of $\Gamma$ let $d(v)$
be the length of the directed path from $v$ to $v_0$. If $e$ is an
edge of $\Gamma$ with $d(s(e)) = k$, then $d(r(e)) = k-1$. The graph
map that takes such an edge $e$ to the edge $e_k$ of $L_n$ (and
takes $s(e)$ to $v_{-k}$ and $r(e)$ to $v_{-k+1}$) is a covering map,
and every graph that admits a surjective cover onto $L_n$ is of this
form. The argument easily extends to the case when $n = \infty$. A
graph that admits a directed cover of $L_4$ is illustrated in
Diagram 3.1.

\begin{center}
\begin{tikzpicture}
\node at (0, 0) (v1) {$\bullet$}; \node at (-2, 0) (v2) {$\bullet$};
\node at (-2, 0.5) (v3) {$\bullet$}; \node at (-2, -0.5) (v4)
{$\bullet$}; \node at (-4, 1) (v5) {$\bullet$}; \node at (-6, 1)
(v6) {$\bullet$}; \node at (-6, 1.5) (v7) {$\bullet$}; \node at (-6,
2) (v8) {$\bullet$}; \node at (-6, 0.5) (v9) {$\bullet$}; \node at
(-8, 1) (v10) {$\bullet$}; \node at (-4, 0) (v11) {$\bullet$}; \node
at (-4, -1) (v12) {$\bullet$}; \node at (-6, -1.5) (v13)
{$\bullet$}; \node at (-6, -0.5) (v14) {$\bullet$}; \node at (-6,
-1) (v15) {$\bullet$}; \node at (-8, -2) (v16) {$\bullet$};
\draw[->] (v2) to (v1); \draw[->] (v4) to (v1); \draw[->] (v3) to
(v1); \draw[->] (v5) to (v3); \draw[->] (v6) to (v5); \draw[->] (v7)
to (v5); \draw[->] (v8) to (v5); \draw[->] (v9) to (v5); \draw[->]
(v10) to (v6); \draw[->] (v11) to (v2); \draw[->] (v12) to (v4);
\draw[->] (v13) to (v12); \draw[->] (v14) to (v12); \draw[->] (v15)
to (v12); \draw[->] (v16) to (v13);
\end{tikzpicture}

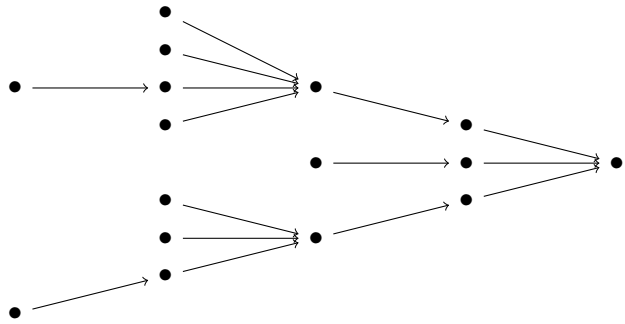
Diagram~3.1 \,\, A directed cover of $L_4$
\end{center}

(b) If $\Gamma$ is not a tree then it must have at least one
non-trivial reduced circuit, and hence $\Gamma$ must have a nontrivial cycle
since every vertex has out-degree $1$ by part (a). Suppose that $\Gamma$
has two distinct cycles $C_1$ and $C_2$ that are not just cyclic
conjugates of each other. These cycles cannot be disjoint. To see
this, note that if $v_1$ is a vertex in $C_1$ and $v_2$ is a vertex
in $C_2$ then there is a path $p = e_1\ldots e_k$ (containing no cycles) from $v_1$ to
$v_2$. Since all vertices
in a cycle have out-degree $1$, there must be indices $i$ and $j$
with $1 \leq i < j \leq k$ such that $e_i \in (\Gamma^1)^*$ and $e_j
\in \Gamma^1$. But this violates the observation in the proof of part (a).
 So the cycles $C_1$
and $C_2$ must have some vertex $v$ in common. Then the cyclic
conjugates of $C_1$ and $C_2$ starting at $v$ must be equal or else
there is some vertex $w$ in $C_1 \cap C_2$ of out-degree at least
$2$, a contradiction. Hence $C_1$ and $C_2$ are cyclic conjugates of
each other.

Suppose that $v'$ is any vertex
 in a non-trivial  cycle $C$ and $v$ is any other vertex. If $v$ is
 on the cycle $C$ then there is a directed path on $C$ from $v$ to $v'$
 that does not include the cycle $C$ as a subpath.
 So suppose that $v$  is not on $C$. Then
 there is a path $p = e_1\ldots e_t$ from $v$ to
$v'$. There is a
largest integer $i$ such that $e_i$ is not an edge in $C$. Since
$r(e_i) \in C$ has out-degree $1$ and there is an edge of $\Gamma$
in $C$ starting at $r(e_i)$ we must have $e_i \in \Gamma^1$. It
follows that all of the edges $e_1,\cdots ,e_i$ are in $\Gamma^1$.
Also there is a unique (possibly empty) directed path $p'$ in $C$
from $r(e_i)$ to $v'$ that does not include $C$ as a subpath. Hence
the path $e_1\ldots e_ip'$ is a directed path from $v$ to $v'$ that
does not include $C$ as a subpath. The uniqueness of such a path
follows easily by an argument similar to that used in part (a). It
is clear that in this case $\Gamma$ is a directed cover of
$B_{\{a\}}$ since every vertex has out-degree $1$. A graph that
admits a directed cover of $B_{\{a\}}$ is illustrated in Diagram
3.2.

\begin{center}
\begin{tikzpicture}
\node at (0, 0) (v1) {$\bullet$}; \node at (1.732, 1) (v2)
{$\bullet$}; \node at (1.732, -1) (v3) {$\bullet$}; \node at (-2, 0)
(v4) {$\bullet$}; \node at (-4, 0.5) (v5) {$\bullet$}; \node at (-4,
-0.5) (v6) {$\bullet$}; \node at (-6, 1) (v7) {$\bullet$}; \node at
(2.7, 2.5) (v8) {$\bullet$};
\draw[->] (v1) to [bend left = 50] (v2); \draw[->] (v2) to [bend
left = 50] (v3); \draw[->] (v3) to [bend left = 50] (v1); \draw[->]
(v4) to (v1); \draw[->] (v5) to (v4); \draw[->] (v6) to (v4);
\draw[->] (v7) to (v5); \draw[->] (v8) to (v2);
\end{tikzpicture}

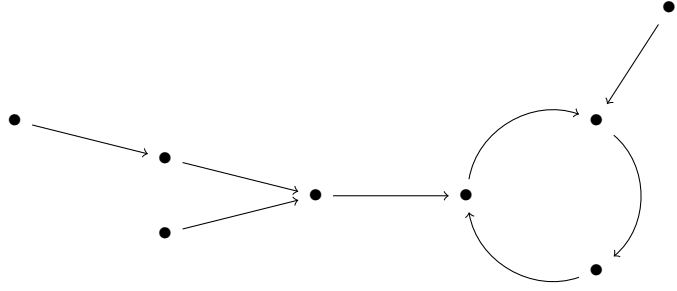
Diagram~3.2 \,\, A directed cover of $B_{\{a\}}$
\end{center}

(c)  If $\Gamma$ is a tree with no sinks, then every vertex of
$\Gamma$ has out-degree $1$, so $\Gamma$ is a directed cover of
$B_{\{a\}}$. If  $\Gamma$ is
not a tree then $\Gamma$ has the structure described in case (b), and hence it is a directed cover of $B_{\{a\}}$. Also, in this case,
since $\Gamma$ has a unique cycle $C$ (up to cyclic conjugates), a
spanning tree for $\Gamma$ contains every edge of $\Gamma$ except
one edge in $C$, so $\pi_1(\Gamma) \cong \mathbb Z$. Conversely, if $\Gamma$ is a directed cover of
$B_{\{a\}}$, then by part (a) it does not have a sink. So if it is a tree, it is a tree with no sinks. If it is not a tree then it has the structure described in part (b), whence $\pi_1(\Gamma) \cong \mathbb Z$ by the argument above. \QED

\medskip

If $f$ is a graph morphism from $\tilde{\Gamma}$ to $\Gamma$ with
$f(\tilde{v}) = s(p)$ for some path $p$ in $\Gamma$ and some vertex
$\tilde{v}$ in $\tilde{\Gamma}$, then we say that $p$ {\em lifts} to
$\tilde{v}$ if there is a path $\tilde{p}$ in $\tilde{\Gamma}$ with
$f(\tilde{p}) = p$ and $s(\tilde{p}) = \tilde{v}$. Note that
directed paths must lift to directed paths if they lift, by the
definition of a graph morphism. It is well-known and easy to prove
that if $f : \tilde{\Gamma} \rightarrow \Gamma$ is a covering map
between graphs, then every path in $\Gamma$ starting at a vertex $v$
lifts to a path at $\tilde{v}$ for every vertex $\tilde{v} \in
f^{-1}(v)$. This is a very special case of the {\em path lifting
theorem} in topology. See Hatcher's book \cite{hatch} for details.
The following easy lemma is the analogous version of this for
directed paths in directed graphs.

\begin{Lemma}
\label{pathlifting}

$($Path lifting lemma for directed covers$)$ A graph morphism $f: \tilde{\Gamma}
\rightarrow \Gamma$ is a directed cover if and only if, for every vertex $v \in \Gamma^0$ and every vertex $\tilde{v} \in f^{-1}(v)$, every directed path $p$ in $\Gamma$ with $s(p) = v$ lifts to a unique path $\tilde{p}$ with $s(\tilde{p}) = \tilde{v}$.

\end{Lemma}

\noindent {\bf Proof.} If $f$ is a directed cover and $e$ is an edge in $\Gamma$ with $s(e) = v$ then by the definition of a directed cover, there is a unique edge $\tilde{e}$ in $\tilde{\Gamma}$ with $f(\tilde{e}) = e$ and $s(\tilde{e}) = \tilde{v}$. This is the basis for an easy inductive proof that directed paths starting at $v$ lift uniquely to directed paths starting at $\tilde{v}$. The proof of the converse statement is equally straightforward. \QED

\medskip

The directed path lifting lemma above does not hold for directed
immersions that are not directed covers in general, but it is easy
to see that maximum initial segments of directed paths in $\Gamma$
lift uniquely to directed paths in $\tilde{\Gamma}$, as described in
the following lemma, the proof of which is a simple adaptation of the proof of Lemma \ref{pathlifting}. The analogous observation for immersions
between graphs may be found in \cite{GM}.

\begin{Lemma}
\label{Liftingimmersions}

$($Path lifting lemma for directed immersions$)$ Let $f: \tilde{\Gamma}
\rightarrow \Gamma$ be a directed immersion between graphs, let $v$
be a vertex of $f(\tilde{\Gamma})$ and let $p$ be a directed path in
$\Gamma$ with $s(p) = v$. Then for every vertex $\tilde{v} \in
f^{-1}(v)$ there is a unique $($possibly empty$)$ maximum initial
segment $p_1$ of $p$  that lifts to a directed path at $\tilde{v}$.
Furthermore, the lift of $p_1$ at $\tilde{v}$ is unique.

\end{Lemma}

For each vertex $v$ of a graph $\Gamma$, let $vI(\Gamma)v$ be the
local submonoid of $I(\Gamma)$ with identity $v$. Since $vpq^*v = 0$
if $pq^*$ is not a circuit at $v$ it follows that the non-zero
elements of $vI(\Gamma)v$ are the circuits of the form $pq^*$ where
$p$ and $q$ are directed (possibly empty) paths with $r(p) = r(q)$
and $s(p) = s(q) = v$. Clearly $vI(\Gamma)v$ is non-trivial (i.e.
does not consist of just $v$ and $0$) if and only if $v$ is not  a
{\em sink} in the graph $\Gamma$ since if $e$ is an edge of $\Gamma$
with $s(e) = v$, then $ee^* \in vI(\Gamma)v$ and $ee^* \neq v$.

Recall our convention that if $f : S \rightarrow T$ is a homomorphism  between inverse semigroups then $f(0_S) = 0_T$. The homomorphism $f$ is
called a {\em $0$-restricted homomorphism} from $S$ to $T$ if in
addition  $f^{-1}(0_T) = \{0_S\}$. We call a function $f : S \rightarrow
T$ a {\em $0$-morphism} if $f(0_S) =0_T$
and $f(st) = f(s)f(t)$ if $st \neq 0$ and we say that it is a {\em
$0$-restricted morphism} if in addition $f^{-1}(0_T) = \{0_S\}$.
Note that a homomorphism from $S$ to $T$ is a $0$-morphism, but
not every $0$-morphism is a homomorphism since we may have
non-zero elements $s,t \in S$ with $st = 0$ but $f(s)f(t) \neq 0$.
For example, let $S$ be the three-element semilattice $S =
\{e_1,e_2,0\}$ where $e_1$ and $e_2$ are idempotents with $e_1e_2 =
0$ and let $T$ be the two-element semilattice $T = \{1,0\}$. The
function $f : S \rightarrow T$ that takes $e_1$ and $e_2$ to $1$ and
$0$ to $0$ is a $0$-morphism that is not a homomorphism. In general,
it is clear that a function $f : S \rightarrow T$ is a
homomorphism if and only if it is a $0$-morphism with the
property that $f(s)f(t) = 0$ if $st = 0$.


A graph morphism $f :\tilde{\Gamma} \rightarrow \Gamma$ induces a
natural  function (which we  denote by $f_*$) from
$I(\tilde{\Gamma})$ to $I(\Gamma)$ in the obvious way. This induced
function $f_*$ maps $0$ to $0$ and maps a nonzero element
$\tilde{p}\tilde{q}^*$ of $I(\tilde{\Gamma})$ to $pq^*$ (where
$f(\tilde{p}) = p$ and $f(\tilde{q}) = q$). By the definition of a
graph morphism it is clear that $pq^*$ is a non-zero element of
$I(\Gamma)$ if  $\tilde{p}\tilde{q}^*$ is a non-zero
element of $I(\tilde{\Gamma})$ since $r(p) = r(q)$ if
$r(\tilde{p}) = r(\tilde{q})$. The induced function $f_*$ is
well-defined by the uniqueness of canonical forms for elements of
$I(\tilde{\Gamma})$ but it is not in general a homomorphism: in fact
by Theorem 20 of \cite{MeMi} it is a homomorphism if and only if the
graph morphism $f :\tilde{\Gamma} \rightarrow \Gamma$ is injective.
However, we have the following fact.

\begin{Theorem}
\label{circuitembed}

Let $f : \tilde{\Gamma} \rightarrow \Gamma$ be a morphism of graphs
with $f(\tilde{v}) = v$ for vertices $\tilde{v} \in \tilde{\Gamma}$
and $v \in \Gamma$.  Then the following hold.

$($a$)$ $f_*$ is a $0$-restricted morphism  from $I(\tilde{\Gamma})$ to
$I(\Gamma)$.

$($b$)$  $f_*$ induces   a $0$-restricted morphism of
$\tilde{v}I(\tilde{\Gamma})\tilde{v}$ into $vI(\Gamma)v$ for all
vertices $\tilde{v}$ of $\tilde{\Gamma}$.

$($c$)$  $f$ is a directed immersion if and only if the  $0$-morphisms
from $\tilde{v}I(\tilde{\Gamma})\tilde{v}$ into $vI(\Gamma)v$
induced by $f_*$ are all injective homomorphisms  $($i.e. embeddings$)$.

$($d$)$  $f$ is a directed cover if and only if the induced
$0$-morphisms from $\tilde{v}I(\tilde{\Gamma})\tilde{v}$ to
$vI(\Gamma)v$ are all full embeddings: that is, the image of
$\tilde{v}I(\tilde{\Gamma})\tilde{v}$ is a full inverse submonoid of
$vI(\Gamma)v$ for all vertices $\tilde{v}$ of $\tilde{\Gamma}$.

\end{Theorem}

\noindent {\bf Proof}. (a) Suppose that $\tilde{p}\tilde{q}^*$ and
$\tilde{p}'\tilde{q}'^*$ are non-zero elements of
$I(\tilde{\Gamma})$ and denote their images under $f_*$ by $pq^*$
and $p'q'^*$ respectively. If
$\tilde{p}\tilde{q}^*\tilde{p}'(\tilde{q}')^*$ is non-zero in
$I(\tilde{\Gamma})$, then  from the multiplication of canonical
forms in graph inverse semigroups  we either have $\tilde{q}$ is a
prefix of $\tilde{p}'$ or $\tilde{p}'$ is a prefix of $\tilde{q}$.
This easily implies that either $q$ is a prefix of $p'$ or $p'$ is a
prefix of $q$. From this and the definition of the multiplication of
canonical forms  it is easy to see that the induced map $f_*$ is a
$0$-morphism. It is in fact a $0$-restricted morphism since
$f_*(\tilde{p}\tilde{q}^*) \neq 0$ if  $\tilde{p}\tilde{q}^*
\neq 0$.

(b) If $\tilde{p}\tilde{q}^*$ is a non-zero element of
$\tilde{v}I(\tilde{\Gamma})\tilde{v}$ then clearly $pq^*$ is a
non-zero element of $vI(\Gamma)v$. It follows from part (a) that the
restriction of $f_*$ to $\tilde{v}I(\tilde{\Gamma})\tilde{v}$ is a
$0$-restricted morphism to $vI(\Gamma)v$.

(c) Now suppose that $f$ is a directed immersion from
$\tilde{\Gamma}$ to $\Gamma$ and let $f(\tilde{v}) = v$. Let
$\tilde{p}\tilde{q}^*$ and $\tilde{p}'\tilde{q}'^*$ be non-zero
elements of $\tilde{v}I(\tilde{\Gamma})\tilde{v}$ and denote their
images under $f_*$ by $pq^*$ and $p'q'^*$ respectively. Suppose that
$\tilde{p}\tilde{q}^*\tilde{p}'(\tilde{q}')^* = 0$ in
$I(\tilde{\Gamma})$. Then $\tilde{q}$ is not a prefix of
$\tilde{p}'$ and $\tilde{p}'$ is not a prefix of $\tilde{q}$. Hence
we may write $\tilde{q} = \tilde{e}_1\ldots
\tilde{e}_k\tilde{e}_{k+1}\ldots \tilde{e}_n$ and $\tilde{p}' =
\tilde{e}_1\ldots \tilde{e}_k\tilde{e}'_{k+1}\ldots \tilde{e}'_m$
for some
    edges $\tilde{e}_i$ and $\tilde{e}'_j$ in $\tilde{\Gamma}$ with
  $s(\tilde{e}_1) = \tilde{v}$ and
$\tilde{e}_{k+1} \neq \tilde{e}'_{k+1}$. (We allow for the
possibility that the common prefix $\tilde{e}_1\ldots \tilde{e}_k$
of $\tilde{q}$ and $\tilde{p}'$ might be empty.) It follows that
$f_*(\tilde{q}) = e_1\ldots e_ke_{k+1}\ldots e_n$ and
$f_*(\tilde{p}') = e_1\ldots e_ke'_{k+1}\ldots e'_m$ (where
$f(\tilde{e}_i) = e_i$ and $f(\tilde{e}'_j) = e'_j$). Then since $f$
is a directed immersion and $\tilde{e}_{k+1} \neq \tilde{e}'_{k+1}$
it follows that $e_{k+1} \neq e'_{k+1}$, whence $q$ is not a prefix
of $p'$ and $p'$ is not a prefix of $q$. Hence $pq^*p'q'^* = 0$.
This implies that the restriction of $f_*$ to
$\tilde{v}I(\tilde{\Gamma})\tilde{v}$ is a $0$-restricted
homomorphism since it is a $0$-restricted morphism by part (b).

Conversely, suppose that the restriction of $f_*$ to
$\tilde{v}I(\tilde{\Gamma})\tilde{v}$ is a homomorphism for all
$\tilde{v}$. Suppose that there are two edges $\tilde{e}_1$ and
$\tilde{e}_2$ in $\tilde{\Gamma}$ with $s(\tilde{e}_1) =
s(\tilde{e}_2) = \tilde{v}$ and $f(\tilde{e}_1) = f(\tilde{e}_2) = e
\in \Gamma^1$. Then $\tilde{e}_1\tilde{e}_1^* \, , \,
\tilde{e}_2\tilde{e}_2^* \in \tilde{v}I(\tilde{\Gamma})\tilde{v}$
and $f_*(\tilde{e}_1\tilde{e}_1^*) = f_*(\tilde{e}_2\tilde{e}_2^*) =
ee^* \in vI(\Gamma)v$. If $\tilde{e}_1 \neq \tilde{e}_2$ then
$\tilde{e}_1^*\tilde{e}_2 = 0$ and so
$\tilde{e}_1\tilde{e}_1^*\tilde{e}_2\tilde{e}_2^* = 0$. But
$f_*(\tilde{e}_1\tilde{e}_1^*\tilde{e}_2\tilde{e}_2^*) =
(ee^*)(ee^*) = ee^* \neq 0$. This violates the assumption that $f_*$
is a homomorphism, and so we must have $\tilde{e}_1 = \tilde{e}_2$.
Hence $f$ is a directed immersion.

Now suppose that $f$ is a directed immersion and
$f_*(\tilde{p}\tilde{q}^*) = f_*(\tilde{r}\tilde{s}^*) = pq^*$ for
some non-zero elements $\tilde{p}\tilde{q}^*$ and
$\tilde{r}\tilde{s}^*$ of $\tilde{v}I(\tilde{\Gamma})\tilde{v}$.
Since $f$ maps directed edges to directed edges, $f(\tilde{p}) = p =
f(\tilde{r})$. That is, $\tilde{p}$ and $\tilde{r}$ are lifts of $p$
at $\tilde{v}$. By the ``uniqueness" part of Lemma
\ref{Liftingimmersions}, this forces $\tilde{p} = \tilde{r}$.
Similarly
$\tilde{q} = \tilde{s}$, so $\tilde{p}\tilde{q}^* =
\tilde{r}\tilde{s}^*$. Hence $f_*$ is an injective map from
$\tilde{v}I(\tilde{\Gamma})\tilde{v}$ to $vI(\Gamma)v$.

(d) Suppose now that $f$ is a directed covering map from
$\tilde{\Gamma}$ to $\Gamma$, let $\tilde{v}$ be a vertex in
$\tilde{\Gamma}$ and $f(\tilde{v}) = v$. By part (c), $f_*$ is an
injective map from $\tilde{v}I(\tilde{\Gamma})\tilde{v}$ to
$vI(\Gamma)v$. From the multiplication in $I(\Gamma)$ it follows
that the non-zero idempotents of $I(\Gamma)$ are of the form $pp^*$
for some directed path $p$ starting at $v$. By Lemma
\ref{pathlifting}, the path $p$ lifts to a unique path $\tilde{p}$
at $\tilde{v}$ and so $pp^*$ lifts to $\tilde{p}\tilde{p}^*$, an
idempotent of $I(\tilde{\Gamma})$, so $f$ induces a full embedding
of $\tilde{v}I(\tilde{\Gamma})\tilde{v}$ into $vI(\Gamma)v$.
Conversely, suppose that  $f_*$ induces a full embedding of
$\tilde{v}I(\tilde{\Gamma})\tilde{v}$ into $vI(\Gamma)v$.
Then if $p$ is a directed path in $\Gamma$ starting at $v$, the
circuit $pp^*$ is an idempotent in $vI(\Gamma)v$, so it is the image
under $f_*$ of some idempotent in
$\tilde{v}I(\tilde{\Gamma})\tilde{v}$, which must be of the form
$\tilde{p}\tilde{p}^*$ for some lift $\tilde{p}$ of $p$ at
$\tilde{v}$. Hence all directed paths in $\Gamma$ starting at all
vertices $v \in \Gamma^1$ lift to all preimages $\tilde{v} \in
f^{-1}(v)$, whence $f$ is a directed covering map by Lemma
\ref{pathlifting}. \QED

\medskip

We remark that the induced maps from
$\tilde{v}I(\tilde{\Gamma})\tilde{v}$ to $ vI(\Gamma)v$ are in
general not surjective since directed circuits in $\Gamma$ do not
necessarily lift to directed circuits in $\tilde{\Gamma}$, even when
$f$ is a cover. However powers of directed circuits do lift to
directed circuits via finite directed covers of $\Gamma$.

\begin{Lemma}
\label{circuitlift}

Let $f: \tilde{\Gamma} \rightarrow \Gamma$ be a directed cover of
finite graphs, let $v$ be a vertex in $f(\tilde{\Gamma})$ and let
$p$ be a directed circuit at $v$. Then there is a vertex $\tilde{v}'
\in f^{-1}(v)$  and a positive integer $n$ such that $p^n$ lifts to
a directed circuit at $\tilde{v}'$.

\end{Lemma}

\noindent {\bf Proof.} Let $\tilde{v}_0$ be any vertex in $f^{-1}(v)$.
By the directed path lifting lemma (Lemma \ref{pathlifting}), $p$
lifts to a  directed  path $\tilde{p}_1$ from $\tilde{v}_0$ to some
vertex $\tilde{v}_1$. Then $f(\tilde{v}_1) = v$ so again $p$ lifts
to a directed path $\tilde{p}_2$ from $\tilde{v}_1$ to some vertex
$\tilde{v}_2$. Continue like this to obtain a sequence of lifted
paths $\tilde{p}_i$ from $\tilde{v}_{i-1}$ to $\tilde{v}_i$. By
finiteness of $\tilde{\Gamma}$ we must have $\tilde{v}_i =
\tilde{v}_{i+n}$ for some $n
> 0$ and $i \geq 0$. Then $p^n$ lifts to the directed circuit
$\tilde{p}_{i+1}\ldots \tilde{p}_{i+n}$ at $\tilde{v}' =
\tilde{v}_i$. \QED

\begin{Theorem}
\label{finitecover}

Let $f : \tilde{\Gamma} \rightarrow B_A$ be a finite directed cover
of the bouquet of $|A|$ circles.  Then for every vertex $\tilde{v}$
in $\tilde{\Gamma}$, $\tilde{v}I(\tilde{\Gamma})\tilde{v}$ contains
a submonoid isomorphic to the polycyclic monoid  $P_A$ if $|A| > 1$,
and $\tilde{v}I(\tilde{\Gamma})\tilde{v}$ contains a submonoid
isomorphic to the bicyclic monoid if $|A| = 1$.

\end{Theorem}

\noindent {\bf Proof}.  Let $e_a$ denote the loop in $B_A$ labeled
by $a \in A$. For each $a \in A$ let $V_a$ be the set of vertices
$\tilde{w}$ in $\tilde{\Gamma}$ that lie on a non-trivial cycle
$\tilde{e}_1\ldots \tilde{e}_s$ such that $s(\tilde{e}_1) =
\tilde{w}$ and $f(\tilde{e}_1) = e_{a}$. Let $A = \{a_1,a_2,\cdots
a_n\}$ and let $V_m = \bigcap_{i = 1,\cdots ,m}V_{a_i}$ for $m \leq
n$. We claim that $V_m \neq \emptyset$ and that if $\tilde{v}$ is
any vertex in $\tilde{\Gamma}$ and $\tilde{e}'$ is any edge with
$s(\tilde{e}') = \tilde{v}$, then there is a directed path
$\tilde{p} = \tilde{e}'\tilde{e}'_2\ldots \tilde{e}'_t$ from
$\tilde{v}$ to some vertex $\tilde{v}'_1 \in V_m$. By the proof of
Lemma \ref{circuitlift}, some power of the loop $e_{a_1}$ lifts to a
directed path $\tilde{p}'$ starting at $r(\tilde{e}')$ and ending at
a vertex in $V_1$, so the directed path $\tilde{e}'\tilde{p}'$ leads
from $\tilde{v}$ to a vertex in $V_1$, and hence the claim is true
if $m = 1$. Assume inductively that it is true if $m = k$. Let
$\tilde{v}$ be any vertex in $\tilde{\Gamma}$, $\tilde{e}'$ an edge starting at $\tilde{v}$, and let $\tilde{e}_1'$
be the (unique) edge in $\tilde{\Gamma}$ with $s(\tilde{e}_1') =
r(\tilde{e}')$ and $f(\tilde{e}_1') = e_{a_{k+1}}$. By the induction
assumption, $\tilde{e}_1'$ can be extended to some directed path
$\tilde{p}_0$ from $\tilde{v}_0 = r(\tilde{e}'_1)$ to a vertex $\tilde{v}_1 \in V_k$.
But then again by the induction hypothesis there is a directed path
$\tilde{p}_1$ from $\tilde{v}_1$ to some vertex $\tilde{v}_2 \in
V_k$ whose first edge projects by $f$ to $e_{a_{k+1}}$. Continue in
this fashion to obtain a sequence of directed paths
$\tilde{p}_1,\tilde{p}_2,\cdots \tilde{p}_i,\cdots $ whose first
edge projects onto $e_{a_{k+1}}$ and with $\tilde{v}_i =
s(\tilde{p}_i)
 \in V_k$ for all $i \geq 1$.   By
finiteness of $\tilde{\Gamma}$ there must be a directed circuit
$\tilde{p}_i\tilde{p}_{i+1}\ldots \tilde{p}_j$ based at
$\tilde{v}_i$ for some $i < j$. If the first edge of $\tilde{p}_i$
is a loop at $\tilde{v}_{i}$ (that projects onto $e_{a_{k+1}}$) then
$\tilde{v}_{i} \in V_{k+1}$. So assume this is not the case.
Choosing $i$ and $j$ minimal, we may
 assume that this circuit
$\tilde{p}_i\tilde{p}_{i+1}\ldots \tilde{p}_j$ is  a cycle. But then
since the first edge of $\tilde{p}_{i}$ projects onto $e_{a_{k+1}}$
we see that in fact $\tilde{v}_i \in V_{k+1}$. The claim then
follows by induction on $k$.

Thus for every vertex $\tilde{v} $ in $\tilde{\Gamma}$ there is a
directed path $\tilde{p}$ from $\tilde{v}$ to some vertex $\tilde{w}
\in V_n$. Denote by $\tilde{q}_a$ a  cycle at $\tilde{w}$ whose
first edge projects onto $e_a$. Then we see that the paths
$\tilde{r}_a = \tilde{p}\tilde{q}_a\tilde{p}^*$ are in
$\tilde{v}I(\tilde{\Gamma})\tilde{v}$ for all $a \in A$. From the
relations in $I(\Gamma)$ it easily follows that
$\tilde{r}_a^*\tilde{r}_a = \tilde{p}\tilde{p}^*$ and
$\tilde{r}_a^*\tilde{r}_b = 0$ if $a \neq b$, so the inverse
subsemigroup of $\tilde{v}I(\tilde{\Gamma})\tilde{v}$ generated by
the elements $\tilde{r}_a$ (for $a \in A)$ is a homomorphic image of
the copy of the polycyclic monoid $P_A$ with identity
$\tilde{p}\tilde{p}^*$ (provided $|A| > 1$). Since the polycyclic
monoid is congruence free (see \cite{Law1}) it follows that this
monoid is isomorphic to $P_A$. A similar argument yields a copy of
 the bicyclic monoid if $|A| = 1$.  \QED

\medskip

\noindent {\bf Remarks} (a) The conclusion of Theorem
\ref{finitecover} is in general false if $\tilde{\Gamma}$ is an
infinite directed cover of $B_A$. For example, if $\tilde{\Gamma}$
is the universal cover of the circle $B_{\{a\}}$ and $\tilde{v}$ is
any vertex of $\tilde{\Gamma}$, then no power of the loop
 in $B_{\{a\}}$ lifts to a circuit at $\tilde{v}$ and $I(\tilde{\Gamma})$ does not contain a copy of
 the bicyclic monoid.

 (b) The
conclusion of Theorem \ref{finitecover} also fails if $f$ is a
directed immersion that is not a directed cover. For example, if
$\tilde{\Gamma}$ is the graph with two vertices $\tilde{v}_1$ and
$\tilde{v}_2$ and one directed edge $\tilde{e}$ from $\tilde{v}_1$
to $\tilde{v}_2$, then there is a directed immersion of
$\tilde{\Gamma}$  into the circle $B_{\{a\}}$, but
$I(\tilde{\Gamma})$ is finite and so does not contain a copy of the
bicyclic monoid.

(c) It is not true in general that if $\tilde{\Gamma}$ is a finite
directed cover of $\Gamma$, then $I(\Gamma)$ embeds in
$I(\tilde{\Gamma)}$. For example, let $\Gamma$ be the graph with two
vertices $v$ and $w$ and two edges $a$ and $b$ from $v$ to $w$, and
let $\tilde{\Gamma}$ be graph with three vertices, $v_1,w_1$ and
$w_2$ and two edges, namely  $a_1$ from $v_1$ to $w_1$ and $b_1$
from $v_1$ to $w_2$. Then the map that sends $v_1$ to $v, \, w_i$ to
$w, \, a_1$ to $a$ and $b_1$ to $b$ is a directed cover but
$I(\Gamma)$ does not embed in $I(\tilde{\Gamma)}$. Thus Theorem
\ref{finitecover} is specific to finite directed covers of a graph
$B_A$.

\section{Universal groups}

Recall that if $S$ and $T$ are inverse semigroups with $0$, then a
function $\theta : S \rightarrow T$ is called a $0$-morphism if
$\theta(0) = 0$ and $\theta(st) = \theta(s)\theta(t)$ if $st \neq
0$. We define the {\em universal group} ${\mathcal U}(S)$ of an
inverse semigroup $S$ with $0$ to be the group generated by the set
$S^\ast = S \setminus \{0\}$  of non-zero elements of $S$ subject to
the relations $s \cdot t = st$ if  $st \neq 0$. Equivalently
(\cite{Law2}), ${\mathcal U}(S)$ may be defined (up to isomorphism)
by the following universal property. Namely, ${\mathcal U}(S)$ is
the group with the property that there is a $0$-morphism $\tau_S : S
\rightarrow {\mathcal U}(S)^0$ such that if $\alpha : S \rightarrow
H^0$ is a $0$-morphism from $S$ to a group $H$ with $0$ adjoined,
then there exists a unique $0$-restricted  homomorphism $\beta : \mathcal{U}(S)^0
\rightarrow H^0$ such that $\beta \circ \tau_S = \alpha$. We say that
$S$ is {\em strongly $E^*$-unitary} if the $0$-morphism $\tau_S$ is
{\em idempotent-pure}, that is $\tau_S^{-1}(1_{{\mathcal U}(S)})$ is
the set of non-zero idempotents of $S$. Lawson shows in \cite{Law2}
that graph inverse semigroups are strongly $E^*$-unitary.

 A homomorphism
$\phi : S \rightarrow T$ between inverse semigroups is called {\em
idempotent-pure} if, for every idempotent $a$ in $T$, $\phi^{-1}(a)$
is a semilattice (i.e. every preimage of an idempotent of $T$ is an
idempotent of $S$). An inverse monoid $S$ is called {\em
factorizable} if for all $a \in S$ there is an element $b$ in the
group of units of $S$ such that $a \leq b$.

\begin{Prop}
\label{univhism}

Let $S$ and $T$ be  inverse semigroups with zero and $\phi$ a
$0$-restricted  homomorphism from $S$ to $T$. Then

$($a$)$ $\phi$ induces a homomorphism $\phi_{\mathcal U}$ from
${\mathcal U}(S)$ to ${\mathcal U}(T)$ such that the following
diagram commutes;
\begin{center}
\begin{tikzpicture} \label{diagram}
\node at (0, 0) (v1) {$T^*$};
\node at (0, 2) (v3) {$S^*$};
\node at (2, 0) (v2) {$\mathcal{U}(T)$};
\node at (2, 2) (v4) {$\mathcal{U}(S)$};
\draw[->] (v1) to node[below] {$\tau_T$} (v2);
\draw[->] (v3) to node[above] {$\tau_S$} (v4);
\draw[->] (v3) to node[left] {$\phi$} (v1);
\draw[->] (v4) to node[right] {$\phi_{\mathcal{U}}$} (v2);
\end{tikzpicture}
\end{center}

$($b$)$  $\phi_{\mathcal U}$ is surjective if $\phi$ is surjective;

$($c$)$  If $\phi_{\mathcal U}$ is injective and $S$ is strongly
$E^*$-unitary, then $\phi$ is idempotent-pure;

$($d$)$ If $S$ is factorizable and $T$ is strongly $E^*$-unitary, then
$\phi_{\mathcal U}$ is injective if $\phi$ is idempotent-pure.

\end{Prop}

\noindent {\bf Proof.}  Since $\phi$ is $0$-restricted it follows
that $\tau_T \circ \phi$ is a $0$-morphism from $S$ to ${\mathcal
U}(T)$. Hence by the universal property of ${\mathcal U}(S)$, there
is a unique homomorphism $\phi_{\mathcal U}$ from ${\mathcal U}(S)$
to ${\mathcal U}(T)$ such that $\phi_{\mathcal U}(\tau_S(a)) =
\tau_T(\phi(a))$ for all $a \in S^*$.  Since $\tau_S$ maps $S^*$
onto the generators of ${\mathcal U}(S)$ and $\tau_T$ maps $T^*$
onto the generators of ${\mathcal U}(T)$, it follows that
$\phi_{\mathcal U}$ is surjective if $\phi$ is surjective.

Suppose that $ \phi_{\mathcal U}$ is injective and $S$ is strongly
$E^*$-unitary. If $\phi(a)$ is an idempotent of $T$ then
$\tau_T(\phi(a)) = \phi_{\mathcal U}(\tau_S(a))$ is the identity of
${\mathcal U}(T)$. Hence  $a$ is an idempotent of $S$, and so $\phi$
is idempotent-pure.

Now suppose that $S$ is factorizable and $T$ is strongly
$E^*$-unitary and that $\phi$ is idempotent-pure. Note that if $S$
is factorizable, then for every element $a \in S^*$, there exists an
element $a'$ in the group of units of $S$ such that $a = e a'$ for
some idempotent $e$. Hence, $\tau_S(a) = \tau_S(a')$. It follows
that, if $a, b \in S^*$, then $\tau_S(a) \tau_S(b) = \tau_S(a')
\tau_S(b') = \tau_S(a'b')$. So every element of ${\mathcal U}(S)$ is
of the form $\tau_S(a)$ for some element $a$ in the group of units
of $S$. If $\phi_{\mathcal U}(\tau_S(a)) = \phi(\tau_T(a))$ is the
identity of ${\mathcal U}(T)$ then since $\phi$ and $\tau_T$ are
both idempotent-pure it follows that $a$ is the identity of $S$, and
so $\tau_S(a)$ is the identity of ${\mathcal U}(S)$, whence
$\phi_{\mathcal U}$ is injective.  \QED

\medskip

\noindent {\bf Remarks.} (1) The converse of part (b) of Proposition
\ref{univhism} is false in general.
For example, let $S$ be the two element semilattice $S = \{e,0\}$
and let $T$ be the three element semilattice $T = \{e,f,0\}$ with
$ef = 0$. Then ${\mathcal U}(S) \cong {\mathcal U}(T)$ is the
trivial group but the obvious embedding of $S$ into $T$ is a
homomorphism that is not surjective.

(2) The converse of part (c) of Proposition \ref{univhism} is also
false in general. For example, let $S = SIM(a,b)$, the symmetric
inverse monoid on two letters, and let $T = SIM(a,b,c)$, the
symmetric inverse monoid on three letters. The identity map on $S$
extends in the obvious way to an idempotent-pure homomorphism $\phi
: S \rightarrow T$. By Example 2.1 in \cite{Law2}, $S$ is strongly
$E^*$-unitary with maximal group image ${\mathcal U}(S) \cong
{\mathbb Z}_2$, the cyclic group of order $2$, while ${\mathcal
U}(T)$ is the trivial group. The homomorphism $\phi_{\mathcal U}$ is
not injective.

\medskip

The following fact is implicit in Lawson's paper \cite{Law2}. We
provide a proof for completeness.

\begin{Theorem}
\label{entire}

For any  graph $\Gamma$, the universal group
$\mathcal{U}(I(\Gamma))$ is isomorphic to $FG(\Gamma^1)$, the free
group on $\Gamma^1$.

\end{Theorem}

\noindent {\bf Proof}. First recall that the non-zero elements of
$\Gamma$ consists of all elements of form $pq^\ast$ where $p, q$ are
directed paths satisfying $r(p) = r(q)$. For each edge $e \in
\Gamma^1$ define $\tau(e) = e$, regarded as a generator for
$FG(\Gamma^1)$ and define $\tau(e^*) = e^{-1} \in FG(\Gamma^1)$. By
the uniqueness of the canonical form for  non-zero  elements of
$I(\Gamma)$, this extends in the obvious way to a well-defined
function $\tau : I(\Gamma) \rightarrow FG(\Gamma^1)^0$ with $\tau(0)
= 0$ and $\tau (pq^\ast) = \mbox{red}(pq^{-1})$, the reduced form of
$pq^{-1}$,  if $r(p) = r(q)$. For any $p_1 q_1^\ast, p_2 q_2^\ast
\in I(\Gamma)^\ast$, $(p_1 q_1^\ast) (p_2 q_2^\ast) \in
I(\Gamma)^\ast$ if and only if either $q_1$ is a prefix of $p_2$ or
$p_2$ is a prefix of $q_1$. In either case, it is routine to see
that $\tau((p_1 q_1^\ast) (p_2 q_2^\ast)) = \tau (p_1 q_1^\ast) \tau
(p_2 q_2^\ast)$. That is to say, $\tau$ is a 0-morphism. Now for any
group $H$ and any 0-morphism $\alpha: I(\Gamma) \rightarrow H^0$, we
easily see that $\alpha(e^\ast) = \alpha(e)^{-1}$ for every $e \in
\Gamma^1$. Since $FG(\Gamma^1)$ is freely generated by the elements
$e \in \Gamma^1$, the map $e \mapsto \alpha(e)$ for $e \in
\Gamma^1$ extends to a unique homomorphism  $\beta : FG(\Gamma^1)
\rightarrow H$ and clearly $\alpha = \beta \circ \tau$. Hence
$FG(\Gamma^1) \cong {\mathcal U}(I(\Gamma))$.  \QED

\begin{Cor}
\label{retract} If $\Delta$ is a subgraph of $\Gamma$ then
 ${\mathcal U}(I(\Delta))$ is a free factor of ${\mathcal
U}(I(\Gamma))$.

\end{Cor}

\noindent {\bf Proof.} Clearly the set of edges of $\Delta$ is a
subset of the set of edges of $\Gamma$, so the result is immediate
from Theorem \ref{entire}. \QED

\medskip

We turn to a description of the universal groups of the local
submonoids in $I(\Gamma)$. The non-zero idempotents of $I(\Gamma)$
are of the form $pp^*$ where $p$ is a directed path in $\Gamma$. We
denote by $\mathcal U(\Gamma,pp^*)$ the universal group of the local
submonoid $pp^*I(\Gamma)pp^*$. In particular, when $p$ is the
trivial path at the vertex $v$, $\mathcal U(\Gamma,v)$ denotes the
universal group of the local submonoid $vI(\Gamma)v$. Recall that
the non-zero elements of the local submonoid $vI(\Gamma)v$ are of
the form $pq^*$ where $p$ and $q$ are directed (or empty) paths with
$s(p) = s(q) = v$ and $r(p) = r(q)$.

Let $V_v = \{w \in \Gamma^0 : $  there is a (possibly empty) directed path
 $p$ in $\Gamma$  from $v$ to $w\}$ and let $\Gamma_v$ denote the
 subgraph of $\Gamma$ induced by the vertices in $V_v$. A subtree
 $T$ of $\Gamma_v$ is called a {\em directed tree at $v$} if $T$
 contains the vertex $v$ and every geodesic path in $T$ from $v$ to some
 other vertex $w$ in $T$ is directed. $T$ is called a {\em directed
 spanning tree of $\Gamma_v$ at $v$} if $T$ is a directed tree at
 $v$ that contains all of the vertices of $\Gamma_v$.

\begin{Lemma}
\label{tree}

Let $v$ be a vertex of a  graph $\Gamma$ and let $T$ be a directed
subtree of $\Gamma_v$ containing the vertex $v$. Then $T$ extends to
a directed spanning tree $T_v$ of $\Gamma_v$.

\end{Lemma}

\noindent {\bf Proof.} The proof of this is a straightforward
application
 of Zorn's Lemma.
 Let $\mathcal T$ be the set of all subtrees $T'$ of $\Gamma_v$ such
 that $T'$
 contains the tree $T$ and $T'$ is directed at $v$. Then
 $\mathcal T$ is a partially ordered set with respect to inclusion
 (i.e. $T_1 \leq T_2$ if and only if $T_1$ is a subtree of $T_2$).
 It is easy to see that the union of  a chain of  trees in $\mathcal T$ is
 another tree in $\mathcal T$, so by Zorn's Lemma $\mathcal T$ has
 a maximal element $T_v$. If $T_v$ is not a spanning tree of
 $\Gamma_v$, then there is some vertex $w$ of $\Gamma_v$ that is not
 in $T_v$. Since there is a directed path in $\Gamma_v$
 from $v$ to $w$, there is some directed path $p$ that
 starts at a vertex $w'$ in $T_v$ and ends at $w$ and has no edge or vertex  other than $w'$ in
 $T_v$. Then $T_v \cup \{p\}$ is a tree  strictly containing $T_v$ as a
 subtree. If $p'$ is the geodesic path in $T_v$ from $v$ to $w'$, then $p'p$ is
 the geodesic path in $T_v \cup \{p\}$ from $v$ to $w$ and $p'p$ is directed,
 so $T_v \cup \{p\} \in \mathcal
 T$. This contradicts the maximality of $T_v$, so $T_v$ is a
 directed
 spanning tree at $v$. \QED

 \begin{Theorem}
 \label{local}
 If $v$ is a vertex of a  graph $\Gamma$  then $\mathcal U(\Gamma,v) \cong \pi_1(\Gamma_v,v)$.

 \end{Theorem}

 \noindent {\bf Proof.} It follows from Lemma \ref{tree}
  (with $T = \{v\}$) that $\Gamma_v$ has a
 directed spanning tree $T_v$ at $v$.  Denote the geodesic path in $T_v$
 from $v$ to a vertex $w$ in $\Gamma_v$ by $p_w$. Thus each path
 $p_w$ is a directed path from $v$ to $w$.
The group $\pi_1(\Gamma_v,v)$ is generated by
 the homotopy classes $[c(e)]$ of circuits of the form
 $c(e) = p_{s(e)}ep_{r(e)}^*$ for each edge $e$ of $\Gamma_v$ that is not in
 $T_v$ (see \cite{Stall} or \cite{hatch} for basic information about homotopy of graphs). We claim that the set $S_v$ consisting of these circuits, together with $\{0\}$ and the
 circuits of the form $pp^*$,
 for $p$ a directed path in $T_v$ starting at $v$,  generates the local submonoid
 $vI(\Gamma)v$ as an inverse submonoid of $I(\Gamma)$.

 To see this, suppose first that $w$
 is a vertex in $\Gamma_v$, $q' = e'_1e'_2\ldots e'_m$ is the geodesic path in $T_v$ from
 $v$ to $w$ and $p = e_1e_2\ldots e_n$ is any other directed
 path from $v$ to $w$. We prove by induction on $n$ that $pq'^*$
 can be expressed as a product of  elements in $S_v$ and their inverses. The result is
 clearly true if $p = q'$ or if $n = 0$ so assume $p \neq q'$ and $n \geq 1$.
If $n =
 1$ then $e_1 \neq q'$ and so $e_1$ is not an edge in $T_v$.  So in this case $pq'^* = c(e_1) \in S_v$.
So assume that $n
 > 1$ and that the result is true for all directed paths $p$ of length less than $n$ from $v$
 to some
 vertex $w$ in $\Gamma_v$.
Since $p \neq q'$ we may
 write $p = e_1e_2\ldots e_je'_ie'_{i+1}\ldots e'_m$ for some $j \leq n$ and
 $e_j \neq e'_{i-1}$. (We allow for the case that $e'_i\ldots e'_m$ is empty).
 Let $p_1$ be the geodesic in $T_v$ from $v$ to $s(e_j)$.
 By the induction
 assumption, the circuit
 $e_1e_2\ldots e_{j-1}p_1^*$ can be written as a product of
 elements of $S_v$ and their inverses. Also $p_1e_j(e'_1\ldots e'_{i-1})^* = c(e_j)$, so
 $e_1\ldots e_{j-1}e_j(e'_1\ldots e'_{i-1})^* =
 e_1\ldots e_{j-1}p_1^*p_1e_j(e'_1\ldots e'_{i-1})^*$ is a product of
 elements in $S_v$ and their inverses.
But then $pq'^* =
 (e_1\ldots e_j(e'_1\ldots e'_{i-1})^*)(q'q'^*)$ is a product of elements of
 $S_v$ and their inverses. Now let $w$ be any vertex in $\Gamma_v$ and $p,q$ any directed paths
 from $v$ to $w$ in $\Gamma_v$. Let $q'$ be the geodesic in $T_v$
 from $v$ to $w$. Then by the argument above, the circuits $pq'^*$
 and $qq'^*$ can be written as products of elements in $S_v$ and their inverses. It
 follows that the circuit $pq^* = (pq'^*)(q'q^*)$ is in the inverse
 submonoid of $I(\Gamma)$ generated by $S_v$. So $vI(\Gamma)v$ is
 generated as an inverse monoid by the elements of $S_v$.

We now claim that every non-zero element of $vI(\Gamma)v$ can be written uniquely  in the form $c(e_1)\ldots c(e_k)pp^*c(e_{k+1})^*\ldots c(e_{k+s})^*$ for some edges $e_i$ in $\Gamma_v$ that are not in $T_v$ and some directed path $p$ in $T_v$ starting at $v$ such that both $p_{r(e_{k})}$ and $p_{r(e_{k+1})}$ are prefixes of $p$. (We allow for the possibility that either $k$ or $s$ (or both) might be zero). To prove this, we first make three  observations.

{\em Observation 1.} If $p $ is a directed path in $T_v$ starting at $v$, $e$ is an edge of $\Gamma_v$ that is not in $T_v$ and $pp^*c(e) \neq 0$, then $pp^*c(e) = c(e)$. To see this, note that if $pp^*c(e) = pp^*p_{s(e)}ep_{r(e)}^* \neq 0$, then  $p$ is a prefix of $p_{s(e)}e$ since $e$ is not an edge in $T_v$.  Hence $p$ is a prefix of $p_{s(e)}$ (again since $e$ is not an edge in $T_v$). Observation 1 follows easily from this.

{\em Observation 2.} If $e$ and $f$ are edges of $\Gamma_v$ that are not in $T_v$ and $c(e)^*c(f) \neq 0$, then $c(e)^*c(f) = p_{r(e)}p_{r(e)}^*$. To see this, note that if $c(e)^*c(f) \neq 0$ then $p_{r(e)}e^*p_{s(e)}^*p_{s(f)}fp_{r(f)}^* \neq 0$.  Since neither $e$ nor $f$ is an edge in $T_v$ it follows  that  $p_{s(e)}e = p_{s(f)}f$ and so $e = f$. This implies that $c(e)^*c(f) = c(e)^*c(e)$, which easily yields Observation 2.

{\em Observation 3.} The set of elements $pp^*$ where $p$ is a directed path in $T_v$ starting at $v$ is a submonoid of $vI(\Gamma)v$. This follows easily since if $p_1p_1^*p_2p_2^* \neq 0$ then either $p_1$ is a prefix of $p_2$ or $p_2$ is a prefix of $p_1$.

It follows from these three observations and the fact that $vI(\Gamma)v$ is generated as an inverse monoid by $S_v$ that every element of $vI(\Gamma)v$ can be written as a product of the form $c(e_1)\ldots c(e_k)pp^*c(e_{k+1})^*\ldots c(e_{k+s})^*$ for some edges $e_i$ in $\Gamma_v$ that are not in $T_v$ and some directed path $p$ in $T_v$ starting at $v$. If $p$ is a prefix of $p_{r(e_{k})}$ then $p^*_{r(e_{k})}pp^* = p^*_{r(e_{k})} = p^*_{r(e_{k})}p_{r(e_{k})}p^*_{r(e_{k})}$ and similarly if $p$ is a prefix of $p_{r(e_{k+1})}$ then $pp^*p_{r(e_{k+1})} = p_{r(e_{k+1})}p_{r(e_{k+1})}^*p_{r(e_{k+1})}$.  So we may assume without loss of generality that both $p_{r(e_{k})}$ and $p_{r(e_{k+1})}$ are prefixes of $p$.

Note that if $c(e_1)c(e_{2}) \neq 0$, then either
$p_{s(e_2)}e_2$ is a prefix of $p_{r(e_1)}$ or $p_{r(e_1)}$ is a
prefix of $p_{s(e_2)}e_2$. Since $e_2$ is not an edge in $T_v$,
 we must have
$p_{r(e_1)}$ is a prefix of $p_{s(e_2)}e_2$. Applying this to all
products $c(e_i)c(e_{i+1})$ we see that if $c(e_1)c(e_2)\ldots
c(e_k) \neq 0$ then we must have $c(e_1)c(e_2)\ldots c(e_k) =
p_{s(e_1)}e_1p_{1,2}e_2p_{2,3}e_3\ldots e_kp_{r(e_k)}^*$ for some
directed paths $p_{i,i+1}$ in $T_v$ from $r(e_i)$ to $s(e_{i+1})$. A similar argument applies to the non-zero product $c(e_{k+1})^*\ldots c(e_{k+s})^*$. Hence we have
\begin{align*}
& c(e_1)\ldots c(e_k)pp^*c(e_{k+1})^*\ldots c(e_{k+s})^* \\
= & p_{s(e_1)}e_1p_{1,2}e_2p_{2,3}e_3\ldots e_kp_{r(e_k)}^*pp^*p_{r(e_{k+1})}e_{k+1}^*p_{k+1,k+2}^*e_{k+2}^*\ldots e_{k+s}^*p_{s(e_{k+s})}^*
\end{align*}
where the $p_{i,i+1}$ are paths in $T_v$. The uniqueness of canonical forms in $I(\Gamma)$ and the fact that each $e_i$ is not in $T_v$ implies that if $$
c(e_1)\ldots c(e_k)pp^*c(e_{k+1})^*\ldots c(e_{k+s})^* = c(e'_1)\ldots c(e'_m)p'p'^*c(e'_{m+1})^*\ldots c(e'_{m+n})^*
$$
then $k = m, s = n, p = p'$ and $e_i = e'_i$ for all $i$.

For $e$ an edge of $\Gamma_v$ not in $T_v$  define $\theta(c(e)) =
[c(e)] \in \pi_1(\Gamma_v,v)$  and also define
  $\theta(pp^*) = 1$ (the identity of
$\pi_1(\Gamma_v,v)$) for
 $p$ a geodesic path in $T_v$ from $v$ to some vertex $w = r(p)$.
By the uniqueness of the expression for non-zero elements of $vI(\Gamma)v$ established above, it follows that  $\theta$ extends to a well-defined function  (again
denoted by $\theta$) from   $vI(\Gamma)v$ to $\pi_1(\Gamma_v,v)^0$. A routine argument, using Observations 1, 2 and 3 above,  shows  that $\theta$ defines a
$0$-morphism from  $vI(\Gamma)v$ to $\pi_1(\Gamma_v,v)^0$. If
$\alpha$ is any other $0$-morphism from $vI(\Gamma)v$ to $G^0$, for
some group $G$ with $0$, then since $\pi_1(\Gamma_v,v)$ is freely
generated by the $\theta(c(e))$'s we see, as in the proof of Theorem
\ref{entire}, that there is a unique homomorphism $\beta :
\pi_1(\Gamma_v,v) \rightarrow G$ that satisfies $\alpha = \beta
\circ \tau$. Hence $\mathcal U(\Gamma,v) \cong \pi_1(\Gamma_v,v)$.
\QED

\begin{Cor}
\label{strong}

If $v$ is a vertex in a graph $\Gamma$ then $\mathcal U(\Gamma,v)$
is a free group. If $u$ and $v$ are vertices in the same strongly
connected component of $\Gamma$, then $\mathcal U(\Gamma,v) \cong
\mathcal U(\Gamma,u)$. In particular, if $\Gamma$ is a strongly
connected graph, then $\mathcal U(\Gamma,v) \cong \pi_1(\Gamma,v)$
is a free group with rank independent of the choice of $v$.

\end{Cor}

\noindent {\bf Proof.} Clearly $\mathcal U(\Gamma,v)$ is a free
group since it is the fundamental group of  a graph by Theorem
\ref{local}. If $u$ and $v$ are in the same strongly connected
component of $\Gamma$ then there is a directed path $p$ from $u$ to
$v$ and a directed path $q$ from $v$ to $u$. If $w$ is any vertex in
$\Gamma_v$, there is a directed path $p'$ from $v$ to $w$ and hence
there is a directed path $pp'$ from $u$ to $w$, whence $w$ is a
vertex in $\Gamma_u$. Hence $\Gamma_v$ is a subgraph of $\Gamma_u$.
Similarly $\Gamma_u$ is a subgraph of $\Gamma_v$, so $\Gamma_u =
\Gamma_v$. Hence by Theorem \ref{local}, $\mathcal U(\Gamma,v) \cong
\pi_1(\Gamma_v,v) \cong \pi_1(\Gamma_u,u) \cong \mathcal
U(\Gamma,u)$. The result about strongly connected graphs follows
immediately since if $\Gamma$ is strongly connected then $\Gamma_v =
\Gamma$. \QED

\begin{Cor}
\label{univsub} If $\Delta$ is a subgraph of a graph $\Gamma$ and
$v$ is a vertex of $\Delta$, then $\mathcal U(\Delta,v)$ is a free
factor of $\mathcal U(\Gamma,v)$.

\end{Cor}

\noindent {\bf Proof.} If there is a directed path in $\Delta$ from
$v$ to some vertex $w$ in $\Delta$, then the same path lies in
$\Gamma$, so $\Delta_v$ is a subgraph of $\Gamma_v$. Let $T_v$ be  a
directed spanning tree for $\Delta_v$ at $v$.  By Lemma \ref{tree},
$T_v$ can be extended to a directed spanning tree $T'_v$ for
$\Gamma_v$ at $v$. Notice that if $e$ is an edge of $\Delta_v$ that
is not in $T_v$ then it is not in $T'_v$ either. Hence the set of
free generators for $\mathcal U(\Delta,v)$ obtained from $T_v$ is
contained in the set of free generators for $\mathcal U(\Gamma,v)$
obtained from $T'_v$. It follows from Theorem \ref{local} that
$\mathcal U(\Delta,v)$ is a free factor of
 $\mathcal U(\Gamma,v)$. \QED

\medskip

\noindent {\bf Remark.} We remark that in general if $\Delta$ is a
 subgraph of the  graph $\Gamma$, then  there may
be vertices of $\Delta$ that are in $\Gamma_v$ but not in $\Delta_v$
since there may be directed paths in $\Gamma$ from $v$ to a vertex
in $\Delta$ but no such directed path in $\Delta$. Also, while the
proof of Corollary \ref{univsub} shows that every directed spanning
tree of $\Delta_v$ at $v$ may be extended to a directed spanning
tree of $\Gamma_v$ at $v$, it is not necessarily true that every
directed spanning tree of $\Gamma_v$ restricts to a directed
spanning tree of $\Delta_v$. This is because in general a geodesic
path from $v$ to some other vertex $w$ in $\Delta$ in  a directed
spanning tree for $\Gamma_v$ may pass through vertices and edges of
$\Gamma \setminus \Delta$.

\begin{Lemma}
\label{Duniv}

If $a$ and $b$ are  $\mathcal D$-related idempotents in an inverse
semigroup $S$ with $0$, then $\mathcal U(aSa) \cong \mathcal U(bSb)$.

\end{Lemma}

\noindent {\bf Proof.} There exists an element $x \in S$ such that
$xx^{-1} = a$ and $x^{-1}x = b$. So  $u$ is a non-zero element of
$aSa$ if and only if $x^{-1}ux$ is a non-zero element of $bSb$. It
follows that the map defined by $u \mapsto x^{-1}ux$ induces an isomorphism
from $\mathcal U(aSa)$ onto $\mathcal U(bSb)$. \QED

\begin{Theorem}
\label{freefact}

Let $p$ be a directed path from  a vertex $v$ to a vertex $w$ in a
graph $\Gamma$. Then

$(a)$ $\mathcal U(\Gamma,pp^*) \cong \mathcal U(\Gamma,w)$.

$(b)$ $\mathcal U(\Gamma,w)$ is isomorphic to a free factor of
$\mathcal U(\Gamma,v)$.

\end{Theorem}

\noindent {\bf Proof.} (a) Note that $pp^* \, \mathcal D \, r(p) =
w$ in $I(\Gamma)$ since $p^*p = r(p)$, so the result of part (a)
follows immediately from Lemma \ref{Duniv}.

(b) If $v$ and $w$ are in the same strongly connected component then
the result follows from Corollary \ref{strong}. So we may assume
that there is a directed path from $v$ to $w$ but no directed path
from $w$ to $v$. Let $p = e_1e_2\ldots e_n$ be a directed path from
$v$ to $w$ and suppose that $k$ is the largest index such that
$s(e_k) \notin \Gamma_w$. That is, there is a directed path from $w$ to
$r(e_k)$ but no directed path from $w$ to $s(e_i)$ for $i = 1,\cdots
, k$.  Clearly every vertex $s(e_i)$ for $i = k+1,\cdots
,n$ is in the same strongly connected component as $w$, so
$\Gamma_w = \Gamma_{s(e_i)}$ for all of these vertices. By Lemma
\ref{tree}, we may choose a directed spanning tree $T_{r(e_k)}$ for
$\Gamma_w$ at $r(e_k) = s(e_{k+1})$.

If in the directed path $e_1e_2\ldots e_{k}$ from $v$ to $r(e_k)$
we have $s(e_i) = s(e_j)$ for some $i \neq j$, then we may omit the
subpath $e_i\ldots e_{j-1}$ to obtain a shorter directed path
$e_1\ldots e_{i-1}e_j\ldots e_{k}$ from $v$ to $r(e_k)$. By
omitting all such circuits in the path $e_1\ldots e_{k}$ we obtain
a directed geodesic path $p'$  from $v$ to $r(e_k)$ consisting of
some of the vertices and edges of the path $e_1\ldots e_{k}$. Let
$T'$ be the subgraph of $\Gamma$ consisting of all of the vertices
and edges of $\Gamma$ contained in the paths $p'q$, for $q$ a
 geodesic path in $T_{r(e_k)}$ starting at $r(e_k)$. Since
no vertex in $p'$ other than $r(e_k)$ lies in $\Gamma_w =
\Gamma_{r(e_k)}$, it follows that $T'$ is a tree with the property
that every geodesic path in $T'$ from $v$ to some vertex in $T'$ is
directed. Clearly, $T'$ contains all of the vertices in $\Gamma_w =
\Gamma_{r(e_k)}$. By Lemma \ref{tree} we may extend $T'$ to a
directed spanning tree $T_v$ for $\Gamma_v$ at $v$.  If $e$ is an
edge in $\Gamma_w$ that is not in $T_{r(e_k)}$, then $e$ is not in
$T_v$ either, so the free generators for $\mathcal U(\Gamma,r(e_k))$
obtained from the spanning tree $T_{r(e_k)}$ are among the free
generators for $\mathcal U(\Gamma,v)$ obtained from the spanning
tree $T_v$. It follows that $\mathcal U(\Gamma,r(e_k))$ is a free
factor of $\mathcal U(\Gamma,v)$. The result then follows since
$\mathcal U(\Gamma,r(e_k)) \cong \mathcal U(\Gamma,w)$ by Corollary
\ref{strong}. \QED

\section{Quotients which are also graph inverse semigroups}

 Recall that if $J$ is an ideal of an inverse semigroup $S$, then $S/J$ denotes the
 {\em Rees quotient} of $S$ by the corresponding {\em Rees congruence}
 $\rho_J$,
 where $a \rho_J b$ if $a = b$ or $a,b \in J$.
  Rees quotients of
graph inverse semigroups  are again graph inverse semigroups as
described in the following theorem \cite[Theorem 7]{MeMi}.

\begin{Theorem} \label{rees}
Let $J$ be an ideal of $I(\Gamma)$. Then $I(\Gamma) / J \cong
I(\Delta)$, where $\Delta^0 = \Gamma^0 \setminus (J \cap \Gamma^0)$,
$\Delta^1 = \{e \in \Gamma^1 : r(e) \not\in J\}$, and the source
mapping and range mapping of $\Delta$ are restrictions of those for
$\Gamma$.
\end{Theorem}

Recall that a congruence $\rho$ on an inverse semigroup $S$ with
$0$ is called {\em $0$-restricted} if $0 \rho = \{0\}$. Notice that
if $\rho$ is an arbitrary congruence on a graph inverse semigroup
$I(\Gamma)$ and $J = 0 \rho$, then $J$ is an ideal of $I(\Gamma)$
and $\rho$ induces in the obvious way a $0$-restricted congruence on
the Rees quotient $I(\Gamma)/J \cong I(\Delta)$ where $\Delta$ is
the graph constructed in Theorem \ref{rees}. Thus  the discussion of
general congruences (other than Rees congruences) on graph inverse
semigroups may be reduced  to that of $0$-restricted congruences on
graph inverse semigroups.

For any $v \in {\Gamma}^0$ with out-degree $1$ we
denote the unique edge in $s^{-1}(v)$ by $e_v$. Let $W$ be a set of vertices with out-degree $1$, let  ${\mathbb Z}^+$ be the set of all positive integers
and let $C(W)$ be the set of all cycles whose vertices lie in $W$. Since
all vertices in $W$ have out-degree one, any two cycles in $C(W)$
are either disjoint or cyclic conjugates of each other. A {\em cycle
function} $f : C(W) \rightarrow {\mathbb Z}^+ \cup \{\infty\}$ is a
function that is invariant under cyclic conjugation. A {\em
congruence pair $(W, f)$ of $\Gamma$} consists of a subset $W$ of
vertices of out-degree $1$ and a cycle function $f$.

Let $\rho$ be a $0$-restricted congruence on $I(\Gamma)$ and $W =
\{v \in \Gamma^0 : e  e^{*} \, \rho \, v = s(e)\}$. Then all vertices of $W$
have out-degree $1$. For $c \in C(W)$, let $f(c)$ be the smallest
positive integer $m$ such that $c^m \, \rho \, s(c)$. If no power of
$c$ is equivalent to $s(c)$, then we define $f(c) = \infty$. Then
$T(\rho) = (W, f)$ is a congruence pair. Conversely, let $(W, f)$ be
a congruence pair of $\Gamma$ and let $\wp(W, f)$ denote the
congruence generated by the relation $\mathbf{R}$ consisting of all
pairs $(e_v e_v^{*}, v)$ for $v \in W$ and $(c^{f(c)}, s(c))$ for $c
\in C(W)$ with $f(c) \in {\mathbb Z}^+$. Then the following theorem
is proved in ~\cite[Theorem~1.3]{Wang}.

\begin{Theorem}\label{cong}
The mapping $T$ from the set of all 0-restricted congruences on
$I(\Gamma)$ to the set of all congruence pairs of $\Gamma$ and the
mapping $\wp$ from the set of all congruence pairs of $\Gamma$ to
the set of all 0-restricted congruences on $I(\Gamma)$ are inverses.
In particular, there exists a one-to-one correspondence between
0-restricted congruences on $I(\Gamma)$ and congruence pairs of
$\Gamma$.
\end{Theorem}

Theorem~\ref{cong} enables us to describe all $0$-restricted
congruences on a graph inverse semigroup for which the quotient is
another graph inverse semigroup.

\begin{Theorem} \label{preserve}
Let $\rho$ be a 0-restricted congruence on $I(\Gamma)$ determined by
the congruence pair $(W, f)$. Then $I(\Gamma) / \rho$ is isomorphic
to a graph inverse semigroup if and only if
\begin{description}
\item[(1)] $W \subseteq \{v \in {\Gamma}^0 : v \mbox{ has out-degree }
1, e_v \mbox{ is a loop at }
v\}$; and
\item[(2)] for any $v \in W$,  $f(e_v) = 1$.
\end{description}
\end{Theorem}

\noindent {\bf Proof.} {\it Sufficiency.} Suppose that conditions
(1) and (2) are satisfied. We proceed to prove that $I(\Gamma) /
\rho$ is isomorphic to the graph inverse semigroup $I(\Delta)$,
where $\Delta$ is the graph with $\Delta^0 = \Gamma^0$, $\Delta^1 =
\Gamma^1 \setminus \{e_v : v \in W \}$, and the source mapping for
$\Delta$ is the restriction of the source mapping for $\Gamma$ and
the range mapping for $\Delta$ is the restriction of the range
mapping for $\Gamma$. By conditions~(1), (2) and Theorem~\ref{cong},
$\rho$ is generated by all pairs $(e_v e_v^{*}, v)$ and $(e_v, v)$
where $v \in W$. However, in the inverse semigroup $I(\Gamma)$, the
relation $(e_v, v) \in \rho$ implies the relation $(e_v e_v^{*}, v)
\in \rho$. Hence $\rho$ is generated by all pairs $(e_v,v)$ where $v
\in W$. Let $\phi$ be the function that maps a loop $e_v$ of
$\Gamma$ at a vertex $v$ in $W$ to the vertex $v$ and fixes all
other vertices and edges of $\Gamma$. Then $\phi$ is a function that
maps the generators of $I(\Gamma)$ to the generators of $I(\Delta)$.
This function $\phi$ extends to a homomorphism which we again denote
by $\phi$ from $I(\Gamma)$ onto $I(\Delta)$. To see this, note that
if $pq^*$ is a non-zero element of $I(\Gamma)$ then $\phi(pq^*) =
pq^*$ if neither $p$ nor $q$ contains an edge $e_v$ that is a loop
at some vertex $v \in W$. If $pq^*$ does contain such a loop $e_v$
then we must  have $p = e_1e_2 \ldots e_ne_v^k$ and  $q =
f_1f_2\ldots f_m e_v^t$  for some $k, t \geq 0$ (with at least one
of $k$ or $t$ greater than $0$). Then we see that $\phi(pq^*)$ is
obtained from $pq^*$ by removing the path $e_v^k(e_v^*)^t$ at the
vertex $r(p) = r(q) = v$. Then from the definition of the
multiplication of canonical forms in $I(\Gamma)$ it is easy to see
that $\phi$ is a homomorphism from $I(\Gamma)$ onto $I(\Delta)$. But
then since $\rho$ is generated by the pairs $(e_v,v)$ where $e_v$ is
a loop at some vertex $v \in W$, it follows that the kernel of
$\phi$ is $\rho$ and so
 the inverse semigroups $I(\Gamma) / \rho$ and $I(\Delta)$ are
 isomorphic.

{\it Necessity}.  We may assume that $W$ is nonempty or else the
congruence $\rho$ determined by the pair $(W,f)$ is the identity
congruence. Suppose that $I(\Gamma)/{\rho} \cong I(\Delta)$ for some
graph $\Delta$.

Recall first  that the idempotents of a graph inverse semigroup
$I(\Gamma)$ are of the form $pp^*$ for some directed  (possibly
empty) path $p$ and that the maximal idempotents in the partial
order correspond to the vertices of $\Gamma$ by
\cite[Lemma~15(3)]{MeMi}. Recall also that in any homomorphism
between inverse semigroups, idempotents lift, and so the idempotents
of the inverse semigroup $I(\Gamma)/{\rho}$ are of the form $(pp^*)
\rho$ for some directed (possibly empty) path $p$ in $\Gamma$.
Suppose that $(pp^*)\rho \geq v\rho$ for some vertex $v$ and
directed path $p$ in $\Gamma$. Then $(pp^*v)\rho = (vpp^*)\rho =
v\rho$. Since $\rho$ is $0$-restricted, $v\rho \neq 0\rho$ and so we
must have $s(p) = v$, in which case it follows that $(pp^*)\rho =
v\rho$. Hence $v\rho$ is maximal in the partial order in the graph
inverse semigroup $I(\Delta) \cong I(\Gamma)/{\rho}$, and so we may
view $v\rho$ as a vertex of $\Delta$ for each vertex $v$ of
$\Gamma$.

Now suppose that condition (1) fails. Then there exists a vertex $v$
in $W$  such that $s(e_v) = v$ and $r(e_v)v = 0$. We have
$e_v^*e_v = r(e_v)$ and $(e_ve_v^*,v) \in \rho$, and so $v\rho \,
\mathcal D \, r(e_v)\rho$ in $I(\Gamma)/{\rho}$. But since $\rho$ is
$0$-restricted, we cannot have $v\rho = r(e_v)\rho$ since $r(e_v) v= 0$
in $I(\Gamma)$. Hence $v\rho$ and $r(e_v)\rho$ are distinct vertices
of $\Delta$ that are ${\mathcal D}$-related in $I(\Delta) \cong
I(\Gamma)/{\rho}$. This, together with condition (1) of the theorem and Theorem \ref{cong}, contradicts \cite[Corollary~2]{MeMi}, and so
condition (1) must hold.  Then for any vertex $v \in W$ the edge
$e_v$ is a loop at $v$. This implies that $(e_ve_v^*)\rho =
(e_v^*e_v)\rho = v\rho$, and so $e_v\rho$ is in the $\mathcal
H$-class of the idempotent $v\rho$ in the graph inverse semigroup
$I(\Delta)$. But it is routine to see that if $pq^*qp^* = qp^*pq^*$
in a graph inverse semigroup, then $p = q$ and so $pq^* = pp^*$, and
so graph inverse semigroups are combinatorial. From this it follows
that $e_v\rho = v\rho$, that is $f(e_v) = 1$. Hence condition (2)
must also hold. \QED

\medskip

\begin{Cor}
\label{conguniv}
 Let $\rho$ be a congruence on $I(\Gamma)$ such that
$I(\Gamma) / \rho$ is isomorphic to a graph inverse semigroup
$I(\Delta')$ and let $J = 0\rho$. Then $\Delta'$ is a subgraph of
$\Gamma$ with set $\Gamma^0\setminus (\Gamma^0 \cap J)$ of vertices.
If $v$ is a vertex of $\Delta'$, then the universal group $\mathcal
U(\Delta',v)$ is a free factor of $\mathcal U(\Gamma,v)$.
\end{Cor}

\noindent {\bf Proof.} $J$ is an ideal of $I(\Gamma)$ and
$I(\Gamma)/J$ is a graph inverse semigroup $I(\Delta)$ as described
in Theorem \ref{rees}. Since $\Delta$ is obtained from $\Gamma$ by
omitting some of the vertices and edges of $\Gamma$, we see that
$\Delta$ is a subgraph of $\Gamma$ with $\Delta^0 = \Gamma^0
\setminus (\Gamma^0 \cap J)$. Furthermore, if $pq^*$ and $p'q'^*$
are non-zero elements of $I(\Gamma)$ that are $\rho$-related, either
$pq^*, p'q'^* \in J$ or $(pq^*,p'q'^*) \in \rho'$ where $\rho'$ is
the $0$-restricted congruence on $I(\Delta)$ that is the restriction
of $\rho$ to $I(\Delta)$. The quotient $I(\Delta)/{\rho'}$ is the
graph inverse semigroup $I(\Delta')$ isomorphic to
$I(\Gamma)/{\rho}$ as described in Theorem \ref{preserve}. By the
proof of Theorem \ref{preserve}, the graph $\Delta'$ is obtained
from $\Delta$ by removing  the loops at some of the vertices of
$\Delta$, so $\Delta'$ is a subgraph of $\Delta$ with the same set
of vertices as $\Delta$.  Hence $\Delta'$ is a subgraph of $\Gamma$
with set $\Gamma^0\setminus (\Gamma^0 \cap J)$ of vertices. The
description of the universal groups then follows from Corollary
\ref{univsub}. \QED

\begin{Cor}
\label{congraphinv}

Let $\rho$ be a congruence on $I(\Gamma)$ such that $I(\Gamma) /
\rho$ is isomorphic to a graph inverse semigroup $I(\Delta')$. Then
$I(\Delta')$ is a retract of $I(\Gamma)$ and ${\mathcal
U}(I(\Delta'))$ is a free factor of ${\mathcal U}(I(\Gamma))$.

\end{Cor}

\noindent {\bf Proof.} From the proof of Corollary \ref{conguniv} we
see that $\Delta'$ is a subgraph of $\Gamma$ so $I(\Delta')$ is an
inverse subsemigroup of $I(\Gamma)$. Again using the notation of
Theorem \ref{preserve} and Corollary \ref{conguniv}, let $\phi$ be
the map from $I(\Gamma)$ to $I(\Delta')$ defined by  $\phi(pq^*) =
0$ if $r(p) \in J$ and $\phi(pe_v^k(e_v^*)^tq^*) = pq^*$ if $r(p)
\not\in J, \, v \in W$, $e_v$ is a loop at $v$ and $f(e_v) = 1$. It is routine to check that $\phi$ is
a semigroup homomorphism from $I(\Gamma)$ onto $I(\Delta')$. Clearly
the restriction of $\phi$ to the inverse subsemigroup $I(\Delta')$
of $I(\Gamma)$ is the identity map, so $\phi$ is a retraction map
and $I(\Delta')$ is a retract of $I(\Gamma)$. The fact that
${\mathcal U}(I(\Delta'))$ is a free factor of ${\mathcal
U}(I(\Gamma))$ follows  from Corollary \ref{retract}. \QED

\section{Leavitt path algebras and Leavitt inverse semigroups}

Let $F$ be a field and let $\Gamma$ be a row finite graph; that is
$|s^{-1}(v)| $ is finite for every vertex $v$ in $\Gamma$. Recall
(see \cite{AAM})  that the Leavitt path algebra $L_F(\Gamma)$ is the
$F$-algebra generated by the set $\Gamma^0 \cup \Gamma^1 \cup
(\Gamma^1)^*$ subject to the relations (1)--(4) defining the graph
inverse semigroup $I(\Gamma)$ and the additional ``Cuntz-Krieger"
relations

\medskip
(5) $ v = \Sigma_{ e \in s^{-1}(v)}ee^*$
for all $v \in
\Gamma^0$ such that $v$ is not a sink.
\medskip

The following fact is immediate from the definition of a Leavitt
path algebra.

\begin{Lemma}
\label{sgpalgebra}

The Leavitt path algebra $L_F(\Gamma)$ corresponding to  a graph
$\Gamma$ and a field $F$ is isomorphic to the algebra
$F_0I(\Gamma)/\langle v - \Sigma_{ e \in s^{-1}(v)}ee^*\rangle$ where the sum is taken over all vertices that are
not sinks and where  $F_0I(\Gamma)$ is the contracted semigroup
algebra of $I(\Gamma)$.

\end{Lemma}

We denote by $LI(\Gamma)$  the multiplicative subsemigroup  of $L_F(\Gamma)$
generated  by  $\Gamma^0 \cup \Gamma^1 \cup
(\Gamma^1)^*$. (Of course $LI(\Gamma)$ is a proper subset of $L_F(\Gamma)$ since the addition and scalar multiplication operations are not used in constructing elements of this subsemigroup.) We will see that $LI(\Gamma)$ is in fact an inverse
semigroup, which we refer to as the {\em Leavitt inverse semigroup}
of the graph $\Gamma$. We will give a  presentation for this
semigroup (as an inverse semigroup) by generators and relations.

To see this, we make use of a  natural basis for $L_F(\Gamma)$ as an
$F$-vector space, as  described in a  paper by Alahmadi, Alsulami,
Jain and Zelmanov \cite{AAJZ}. For each vertex $v$ which is not a
sink, choose an edge $\gamma(v)$ such that $s(\gamma(v)) = v$ and
refer to this as a {\em special} edge. Then the following theorem
was proved in \cite{AAJZ}.

\begin{Theorem}
\label{basis}

The following elements form a basis for the Leavitt path algebra
$L_F(\Gamma)$: $(i)$ $v$, where $v \in \Gamma^0$; $($ii$)$ $p, \, p^*$,
where $p$ is a directed path in $\Gamma$; $($iii$)$ $pq^*$ where $p =
e_1...e_n, q = f_1...f_m, \, e_i, f_j \in \Gamma^1, \, r(e_n) =
r(f_m)$, and either $e_n \neq f_m$ or $e_n = f_m$ but this edge $e_n
= f_m$ is not special.

\end{Theorem}

We refer to the basis constructed in Theorem \ref{basis} as the {\em
natural} basis for $L_F(\Gamma)$.

Let $L(\Gamma)$ be the semigroup generated by the set $\Gamma^0 \cup
\Gamma^1 \cup (\Gamma^1)^*$ subject to the relations (1)--(4) used to
define the graph inverse semigroup $I(\Gamma)$ and the additional
relations:

\medskip

(v)  $e_ve_v^* = v$ for each vertex $v \in \Gamma^0$ of out-degree $1$.

\medskip

 Clearly $L(\Gamma)$ is an inverse semigroup since it is a
 homomorphic image of $I(\Gamma)$.

\begin{Theorem}
\label{leavittinv}

For each graph $\Gamma$, $LI(\Gamma) \cong L(\Gamma)$.
In particular, $LI(\Gamma)$ is an inverse semigroup. Every element of $LI(\Gamma)$ is uniquely expressible in one of the forms

$(a)$ $pq^*$ where $p = e_1...e_n$ and $  q =
f_1...f_m$ are $($possibly empty$)$ directed  paths with $r(e_n) =
r(f_m)$ and $e_n \neq f_m$; or

$(b)$ $pq^* = p'ee^*q'^*$ where $p'$ and
$q'$ are $($possibly empty$)$ directed paths with $r(p') = r(q')$ and
the vertex $s(e) = r(p') = r(q')$ has out-degree at least $2$.

\end{Theorem}

\noindent {\bf Proof.} Since  $LI(\Gamma)$ is generated by
$\Gamma^0 \cup \Gamma^1 \cup (\Gamma^1)^*$,  it
satisfies all of the relations (1)-(4) defining the graph inverse
semigroup $I(\Gamma)$.
Also, from the additional relations (5) used to define a
Leavitt path algebra, it follows that $e_ve_v^* = v$ if $v$ has out-degree $1$, so $LI(\Gamma)$ satisfies the relations (v) also. Hence
$LI(\Gamma)$ is a homomorphic image of $L(\Gamma)$ and so it is an inverse semigroup. Any non-zero
element of $LI(\Gamma)$ is expressible in the form $pq^*$ for some
(possibly empty) directed paths $p, q$ in $\Gamma$. If $p =
e_1e_2...e_ne_v$ and $q = f_1f_2...f_me_v$ are directed paths ending
in the same edge $e_v$ (where $v = s(e_v)$ is a vertex of out-degree $1$), it follows from (v) that $pq^* = p_1q_1^*$ where $p_1 =
e_1...e_n$ and $q_1 = f_1...f_m$.

Thus by induction we see that all elements of $LI(\Gamma)$ are
expressible in the form (a) or (b) in the statement of the theorem.
The elements of the form (a) are in the natural basis for
$L_F(\Gamma)$ and the elements of the form (b) are also in the
natural basis for $L_F(\Gamma)$ provided $e$ is not a special edge. If $e$ is a special edge, then the relations (5) imply that $p'ee^*q'^* =
p'q'^*-\Sigma_{i} p'g_ig_i^*q'^*$ in $L_F(\Gamma)$, where the sum is
taken over all edges $g_i$ such that $s(g_i) = r(p')$ and $g_i \neq
e$. Again by applying the relations (v), we see that $p' q'^*$ is equal to an element in $LI(\Gamma)$ of form (a) or (b). Thus inductively,
$$
p'ee^*q'^* = p_0q_0^* - \Sigma_{i_1} p_1g_{i_1}g_{i_1}^*q_1^* - \ldots - \Sigma_{i_s} p_sg_{i_s}g_{i_s}^*q_s^*
$$
where the sum $\Sigma_{i_j}$ is taken over all edges $g_{i_j}$ such that $s(g_{i_j}) = r(p_j)$ and $g_{i_j} \neq \gamma(r(p_j))$, $p_s = p'$, $q_s = q'$, $p_0q_0^*, p_1 q_1^*, \cdots, p_s q_s^*$ have strictly ascending lengths and $p_0 q_0^*$ is of form (a) or (b) which is in the natural basis. (Note that $p'ee^*q'^*$ is essentially determined by the last sum in the above formula.) Since elements of $L_F(\Gamma)$ can
be expressed uniquely as linear combinations of the elements in the
natural basis, it follows that two elements $pq^*$ and $rs^*$ of
$LI(\Gamma)$ that are written either in form (a) or in form (b) are
equal in $L_F(\Gamma)$ (and hence in $LI(\Gamma)$) if and only if $p
= r$ and $q = s$. But since $L(\Gamma)$ satisfies the relations
(1)-(4) and (v) it follows that every non-zero element of
$L(\Gamma)$ may also be expressed in one of the forms (a) or (b). The same argument that is used to prove uniqueness of canonical forms of non-zero elements in $I(\Gamma)$ shows that
 two such elements $pq^*, rs^*$ of
$L(\Gamma)$ are equal in $L(\Gamma)$ if and only if $p = r$ and $q =
s$. Hence  $LI(\Gamma)$ and $L(\Gamma)$ are isomorphic since they
have the same generators and their elements can be expressed in the
same canonical forms. \QED

\medskip

We remark that an alternative proof of the fact that $LI(\Gamma) \cong L(\Gamma)$ using Gr\"{o}bner-Shirshov bases has been provided by Fan and Wang
\cite{FW}\footnote[1]{Based on an earlier version of this paper, David Milan asked whether the kernel of the natural homomorphism from a graph inverse semigroup $I(\Gamma)$ onto the corresponding Leavitt inverse semigroup $LI(\Gamma)$ coincides with the congruence $\leftrightarrow$ introduced by Lenz \cite{lenz}. This is in fact the case. We present a proof of this  in the addendum to this paper.}

\begin{Cor}
\label{univleavitt}

For each graph $\Gamma$ and each vertex $v$ in $\Gamma$ the
universal groups ${\mathcal U}(I(\Gamma))$ and ${\mathcal
U}(LI(\Gamma))$ are isomorphic and the universal group of the local
submonoid $vI(\Gamma)v$ is isomorphic to the universal group of the
local submonoid $vLI(\Gamma)v$.

\end{Cor}

\noindent {\bf Proof.} Imposing the additional relations $e_ve_v^* =
v$ on the generators for $I(\Gamma)$ does not change  the universal
group since the relation $e_ve_v^{-1} = 1$ holds in any group. \QED

\medskip

We say that the directed path $p = e_1e_2...e_n$ in a graph $\Gamma$
{\em has exits} if at least one of the vertices $s(e_i)$ has
out-degree greater than $1$ (and in this case we say that $p$ has an
exit at $s(e_i)$). In particular, an edge $e \in \Gamma^1$ has exits
if and only if $s(e)$ has out-degree greater than $1$. We say that
the directed path $p = e_1e_2...e_n$   has {\em no exits} (or that
$p$ is an {\em NE path}) if every vertex $s(e_i), \, i = 1,..., n$
has out-degree $1$. We
also define the empty path at any vertex $v$ to be an NE path.

\begin{Cor}
\label{levinvidemp}

For each graph $\Gamma$ the non-zero idempotents of $LI(\Gamma)$ are
the elements of the form $pp^*$ where $p$ is a directed path in
$\Gamma$. Furthermore, $pp^* = qq^*$ in $LI(\Gamma)$ if and only if
either $p = qp_1$ for some NE path $p_1$ or $q = pq_1$ for some NE
path $q_1$. In particular, $pp^* = v$ in $LI(\Gamma)$ for some $v
\in \Gamma^0$ if and only if $v = s(p)$ and $p$ is an NE path.

\end{Cor}

\noindent {\bf Proof.} It is clear from the relations defining a
Leavitt inverse semigroup that every non-zero element of
$LI(\Gamma)$ is of the form $pq^*$ where $p$ and $q$ are directed
paths with $r(p) = r(q)$. It is also routine to see  that $pq^*$ is
a non-zero idempotent of $LI(\Gamma)$ if and only if $p = q$, and
that $pp^* = pp_1p_1^*p^*$ if $p_1$ is an NE path with $s(p_1) =
r(p)$. Suppose conversely that $pp^* = qq^*$ for some directed paths
$p$ and $q$. Then $pp^* = s(p)pp^* = s(p)qq^* \neq 0$ so $s(p) =
s(q)$. If $p$ is not a prefix of $q$ and $q$ is not a prefix of $p$
then there exist edges $e_1,e_2$ and paths $s,p',q'$ with $p =
se_1p', \, q = se_2q'$ and $e_1 \neq e_2$. From $se_1p'p'^*e_1^*s^*
= se_2q'q'^*e_2^*s^*$ we see, on premultiplying by $s^*$ and
postmultiplying by $s$ that $e_1p'p'^*e_1^* = e_2q'q'^*e_2^*$. Hence
$e_2^*e_1p'p'^*e_1^* = e_2^*e_2q'q'^*e_2^* = q'q'^*e_2^* \neq 0$.
But since $e_1 \neq e_2$, we see that $e_2^*e_1 = 0$, a
contradiction. Hence we must have either $q$ is a prefix of $p$ or
$p$ is prefix of $q$. In the first case we have $p = qp_1$ for some
directed path $p_1$. Then from $qq^* = qp_1p_1^*q^*$ we see as above
that $p_1p_1^* = q^*q = r(q) = s(p_1)$. Then by an argument very
similar to the argument above, we see that $p_1$ is an NE path.
Similarly, if $p$ is a prefix of $q$ then $q = pq_1$ for some NE
path $q_1$. \QED

\medskip

Recall  that a graph $\Gamma$ admits a directed immersion into a
circle  $B_{\{a\}}$ if and only if all of its vertices have
out-degree at most $1$: the structure of such graphs is described in
Theorem \ref{out1}.
We next provide a  straightforward classification of the Leavitt
inverse semigroups and Leavitt path algebras of such graphs. We may
assume that such a graph is connected, since the Leavitt inverse
semigroup of a graph is the $0$-direct union of the Leavitt inverse
semigroups of the connected components of the graph.

We recall (see \cite{Law1}) that for each non-empty set $A$ and each
group $G$, the  {\em Brandt semigroup} $B_A(G)$ is  the semigroup
 $B_A(G) = \{(a,g,b) : a, b \in A, \, g \in G\} \cup \{0\}$ with
multiplication

\begin{center}

$(a,g,b)(c,h,d) = (a,gh,d)$ if $b = c$ and $0$ otherwise.

\end{center}

\begin{Theorem}
\label{classout1}

Let $\Gamma$ be a  connected graph that immerses into a circle.

$(a)$ If $\Gamma$ is a tree then $LI(\Gamma) \cong B_{\Gamma^0}(1)$,
the combinatorial $|\Gamma^0| \times |\Gamma^0|$ Brandt semigroup;

$(b)$ If $\Gamma$ is  not a tree then $LI(\Gamma) \cong
B_{\Gamma^0}(\mathbb Z)$, the $|\Gamma^0| \times |\Gamma^0|$ Brandt
semigroup with maximal subgroups isomorphic to $\mathbb Z$.

\end{Theorem}

\noindent {\bf Proof.}  If $e$ is an edge of $\Gamma$ then from the
relations defining $LI(\Gamma)$ and the fact that $s(e)$ has out-degree
$1$ we see that $ee^* = s(e)$ and $e^*e = r(e)$ so $s(e)$ and $r(e)$
are $\mathcal D$-related in $LI(\Gamma)$. Also, if $p$ is a directed
path starting at a vertex $v$, then by induction on the length of
$p$ we easily see that $pp^* = v$ in $LI(\Gamma)$.  These facts, together with the fact that $\Gamma$ is connected, imply that
$LI(\Gamma)$ is a $0$-bisimple inverse semigroup whose idempotents
may be identified with the vertices of $\Gamma$. Since $v_1v_2 = 0$
if $v_1 \neq v_2 \in \Gamma^0$, this implies that $LI(\Gamma)$ is a
homomorphic image of a Brandt semigroup with $|\Gamma^0|$ rows
($\mathcal R$-classes) and $|\Gamma^0|$ columns ($\mathcal
L$-classes). By  Theorem \ref{leavittinv}, we see that
every element of $LI(\Gamma)$ may be expressed uniquely in the form
$pq^*$ where $p$ and $q$ are (possibly empty) directed paths with
$r(p) = r(q)$ and the last edge in the path $p$ is different from
the last edge in $q$. Hence distinct vertices of $\Gamma$ remain
distinct as elements of $LI(\Gamma)$ and so $LI(\Gamma)$ is a Brandt
semigroup with $|\Gamma^0|$ rows and columns. The corresponding
maximal subgroups are trivial if $\Gamma$ is a tree and isomorphic
to a homomorphic image of $\mathbb Z$ otherwise, by Theorem
\ref{out1}. But by the canonical form for elements of $LI(\Gamma)$
described in Theorem \ref{leavittinv}, no two distinct powers of a circuit in $\Gamma$ are
equal in $LI(\Gamma)$, so the maximal subgroups of $LI(\Gamma)$ are
isomorphic to $\mathbb Z$ if $\Gamma$ is not a tree.  \QED

\begin{Cor}
\label{classout2}

Let $\Gamma$ be a connected graph that immerses  into a circle and
let $F$ be a field. Then

$(a)$ If $\Gamma$ is a tree then $L_F(\Gamma) \cong
M_{|\Gamma^0|}(F)$, the algebra of $|\Gamma^0| \times |\Gamma^0|$
matrices with entries in $F$ and only finitely many non-zero entries
in each row and column.

$(b)$ If $\Gamma$ is  not a tree, then $L_F(\Gamma) \cong
M_{|\Gamma^0|}(F[x,x^{-1}])$ where $F[x,x^{-1}]$ is the algebra of
Laurent polynomials over $F$ $($i.e. the semigroup algebra
 $F\mathbb Z)$.

\end{Cor}

\noindent {\bf Proof.} By Lemma \ref{sgpalgebra} and the fact that
all vertices have out-degree at most $1$ we have $L_F(\Gamma) \cong
F_0I(\Gamma)/\langle ee^* - s(e) : e \in \Gamma^1 \rangle$ where
$F_0I(\Gamma)$ is the contracted semigroup algebra of $I(\Gamma)$.
But since the relation $ee^* = s(e)$ holds in $LI(\Gamma)$ for all
$e \in \Gamma^1$, this implies that $L_F(\Gamma) \cong
F_0LI(\Gamma)$, the contracted semigroup algebra of the Leavitt
inverse semigroup $LI(\Gamma)$. The result then follows from Theorem
\ref{classout1}. \QED

\medskip

\noindent {\bf Remark} We remark that the characterization given in
Corollary~\ref{classout2}(b) of Leavitt path algebras of a graph
$\Gamma$ that admits a directed cover of the circle  is  a special
case of the characterization of Leavitt path algebras of a class of
graphs given in Proposition 3.5 of \cite{APPM}. This is because by
Theorem \ref{out1} there is a one-one correspondence between the
vertices of $\Gamma$ and the directed paths that end in a specified
vertex of the unique cycle $C$ in $\Gamma$ and do not include $C$ as
a subpath.

\begin{Theorem}
\label{algiso1}

Let $\Gamma$ and $\Delta$ be connected graphs that immerse into a
circle and let $F$ be a field. Then the following are equivalent.

$(a)$ $LI(\Gamma) \cong LI(\Delta)$;

$(b)$   $L_F(\Gamma) \cong L_F(\Delta)$;

$(c)$ $|\Gamma^0| = |\Delta^0|$ and  either $\Gamma$ and $\Delta$ are
both trees  or
  $\pi_1(\Gamma) \cong \pi_1(\Delta)
\cong \mathbb Z$.

\end{Theorem}

\noindent {\bf Proof.} The equivalence of (a) and (c) follows
immediately from  Theorems~\ref{out1} and \ref{classout1}
since two Brandt
semigroups are isomorphic if and only if they have isomorphic
maximal subgroups and the same number of rows.

Suppose  that $LI(\Gamma) \cong LI(\Delta)$. The non-zero elements
of $LI(\Gamma)$ are precisely the non-zero elements in a natural
basis for $L_F(\Gamma)$ so an isomorphism between $LI(\Gamma)$ and
$LI(\Delta)$ is a bijection between the natural bases of
$L_F(\Gamma)$ and $L_F(\Delta)$ that also preserves multiplication
of basis elements in the algebras, so it induces an isomorphism
between $L_F(\Gamma)$ and $L_F(\Delta)$. Hence (a) implies (b).

Conversely suppose that $L_F(\Gamma) \cong L_F(\Delta)$. If $\Gamma$
is a  tree then in particular $\Gamma$ is acyclic, so from Theorem 1
of \cite{AR} it follows that $L_F(\Gamma)$ is von-Neumann regular.
Hence $L_F(\Delta)$ is von-Neumann regular, from which it follows,
again by Theorem 1 of \cite{AR}, that $\Delta$ is acyclic and hence
since the out-degree of every vertex of $\Gamma$ is at most  $1$,
$\Delta$ is a tree. Thus $\Gamma$ is a tree if and only if $\Delta$
is a tree. If $\Gamma$ and $\Delta$ are both trees, then by Corollary~\ref{classout2}(a) $L_F(\Gamma) \cong M_{|\Gamma^0|}(F)$ and
$L_F(\Delta) \cong M_{|\Delta^0|}(F)$. So if $L_F(\Gamma) \cong
L_F(\Delta)$ it follows that $|\Gamma^0| = |\Delta^0|$ in this
case. If $\Gamma$ and $\Delta$ are not trees, then by Corollary~\ref{classout2}(b), $L_F(\Gamma) \cong M_{|\Gamma^0|}(F[x,x^{-1}])$
and $L_F(\Delta) \cong M_{|\Delta^0|}(F[x,x^{-1}])$. It is
well-known that if $R$ and $S$ are commutative
 rings
 then $M_n(R) \cong M_m(S)$ if and only if $R \cong S$ and $m = n$.
 Hence if $L_F(\Gamma) \cong L_F(\Delta)$ and $\Gamma$ is not a tree
 then it again follows that $|\Gamma^0| = |\Delta^0|$. Hence (b)
 implies (c). \QED

 \medskip

We will prove that the implication (a) implies (b) of Theorem
\ref{algiso1} holds for arbitrary connected graphs; that is, if
$LI(\Gamma) \cong LI(\Delta)$ then $L_F(\Gamma) \cong L_F(\Delta)$
(Theorem \ref{algiso} below). However the converse is false in
general as the following example shows.

\medskip

\noindent {\bf Example} Given the following two graphs,
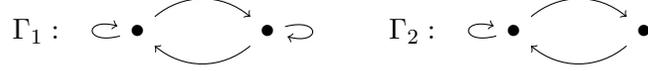
\begin{center}
\begin{tikzpicture}
\node (v1) {$\bullet$}; \node [left of = v1, node distance = 3.5em]
(v0) {$\Gamma_1:$}; \node [right of = v1, node distance = 4.5em]
(v2) {$\bullet$}; \draw[->] (v1) edge [loop left] (v1) to [bend left
= 40] (v2); \draw[->] (v2) edge [loop right] (v2) to [bend left =
40] (v1); \node [right of = v2, node distance = 5em] (v3)
{$\Gamma_2:$}; \node [right of = v3, node distance = 3.5em] (v4)
{$\bullet$}; \node [right of = v4, node distance=4.5em] (v5)
{$\bullet$}; \draw[->] (v4) edge [loop left] (v4) to [bend left =
40] (v5); \draw[->] (v5) to [bend left = 40] (v4);
\end{tikzpicture}
\end{center}
we see from \cite[Example~2.2]{AALP} that $L_F(\Gamma_1) \cong
L_F(\Gamma_2)$. However, $LI(\Gamma_1)$ is not isomorphic to
$LI(\Gamma_2)$.  This is because, according to
Theorem~\ref{leavittinv}, $LI(\Gamma_1) \cong I(\Gamma_1)$ since
every vertex in $\Gamma_1$ has out-degree  2 whereas
$LI(\Gamma_2)$ is not a graph inverse semigroup by
Theorem~\ref{preserve}. Alternatively we can use Theorem
\ref{leavinvclass}  below to see that these Leavitt inverse
semigroups are not isomorphic.

\begin{Lemma}
\label{maxlidemp} For each graph $\Gamma$ we have the following:

$(a)$ $\Gamma^0$ is the set of maximal idempotents in $LI(\Gamma)$.

$(b)$   $ \{p e e^* p^* : p$ is an NE path, $e \in \Gamma^1$ and the
out degree of $s(e)$ is at least $2\}$ is the set of maximal
idempotents of $LI(\Gamma) \setminus \Gamma^0$.

\end{Lemma}

\noindent {\bf Proof.} (a)   By Corollary \ref{levinvidemp},
the non-zero idempotents of $LI(\Gamma)$
are of the form $pp^*$ for some (possibly empty) directed path $p$
in $\Gamma$. Now suppose that $pp^* \geq v$ for some idempotent
$pp^*$ in $LI(\Gamma)$ and some $v \in \Gamma^0$. Then $pp^*v =
vpp^* = v$ in $LI(\Gamma)$. This forces $v = s(p)$, and $vpp^* =
pp^*$, so $v = pp^*$.  Hence $v$ is a maximal idempotent in
$LI(\Gamma)$.

(b) If $q q^* \geq p  e e^* p^*$ where  $q \neq pe$, then we see
from $(q q^*) (p e e^* p^*) = p  e e^* p^*$ that $q$ is a prefix of
$p$. If $p$ is an NE path then $q$ is also an NE path so we get $q
q^\ast = s(q) = s(p) \in \Gamma^0$. Furthermore, it is clear that
any idempotent $p_1 p_1^*$ for which $p_1$ is not an NE path  is
less than or equal to some $p  e e^* p^*$ where $p$ is an NE path,
$e \in \Gamma^1$ and the out-degree of $s(e)$ is at least $2$. \QED

\begin{Lemma}
\label{levisolem} Let $\phi$ be an isomorphism between the Leavitt
inverse semigroups $LI(\Gamma)$ and $LI(\Delta)$ for some graphs
$\Gamma$ and $\Delta$. Then

$(a)$ $\phi$ preserves vertices; that is, $\phi(v) \in \Delta^0$ for
each $v \in \Gamma^0$;

$(b)$ for any nonzero $p q^* \in LI(\Gamma)$, if $\phi(p q^*) = p_1
q_1^*$ and $q$ is an NE path, then $q_1$ is an NE path, $\phi(s(p))
= s(p_1), \phi(p p*) = p_1 p_1^*$ and $\phi(s(q)) = q_1 q_1^* =
s(q_1)$;

$(c)$ for any nonzero $p q^* \in LI(\Gamma)$, if $\phi(p q^*) = p_1
q_1^*$ and $p, q$ are NE paths, then $p_1, q_1$ are NE paths,
$\phi(s(p)) = p_1 p_1^* = s(p_1)$ and $\phi(s(q)) = q_1 q_1^* =
s(q_1)$;

$(d)$ for any $e \in \Gamma^1$, if $s(e)$ has out-degree at least $2$,
then there exist NE paths $p_1, p_2, p_3$ and an edge $\check{e}$
for which  $s(\check{e})$  has out-degree at least $2$ such that $\phi(e) =
p_1 \check{e} p_2 p^*_3$, and there exist NE paths $q_1, q_2, q_3$
such that $\phi^{-1}(\check{e}) = q_1 e q_2 q^*_3$;

$(e)$ for any $v \in \Gamma^0$, if $s^{-1}(v) =\{e_1, \cdots, e_n\}$
with $n \geq 2$, then there exist NE paths $p, p_i, q_i$ and
distinct edges $\check{e}_i, i = 1, \cdots, n$ such that $\phi(e_i)
= p \check{e}_i p_i q^*_i, i = 1, \cdots, n$ and  $s^{-1}(r(p))
=\{\check{e}_1, \cdots, \check{e}_n\}$.

\end{Lemma}

\noindent {\bf Proof.} (a)  This follows from Lemma
\ref{maxlidemp}(a) since $\phi$ must map maximal idempotents in
$LI(\Gamma)$  to maximal idempotents in $LI(\Delta)$.

(b) If $\phi(pq^*) = p_1q_1^*$ and $q$ is NE, then $\phi(pp^*) =
\phi(p r(q) p^*) = \phi(pq^* qp^*) = p_1 q_1^* q_1p_1^* = p_1 r(q_1)
p_1^* = p_1p_1^*$ and $\phi(s(q)) = \phi(qq^*)  = \phi(q r(p) q^*) =
\phi(qp^*pq^*) = q_1 p_1^*p_1 q_1^* = q_1 r(p_1) q_1^* = q_1q_1^*$. Since $\phi(s(q)) \in \Delta^0$ by Lemma
\ref{maxlidemp}(a), this implies that $q_1q_1^* = v$ in $LI(\Delta)$
for some $v \in \Delta^0$. This forces $v = s(q_1)$ and $q_1$ is an
NE path by  Corollary \ref{levinvidemp}. Also, $p_1q_1^* =
\phi(pq^*) = \phi(s(p)pq^*) = \phi(s(p))\phi(pq^*) =
\phi(s(p))p_1q_1^* \neq 0$ so we must have $\phi(s(p)) = s(p_1)$
since $\phi(s(p)) \in \Delta^0$ by part (a) of this lemma.

(c) Note that $\phi(p q^*) = p_1 q_1^*$ implies $\phi(q p^*) = q_1
p_1^*$. This part follows directly from part (b).

(d) Let $e$ be an edge with $s(e)$ of out-degree at least $2$ and
suppose that $\phi(e) = p q^*$. By part (b) we see that $q$ is an NE path.  Also by Lemma \ref{maxlidemp}(b),
$ee^*$ is a maximal idempotent in $LI(\Gamma) \setminus \Gamma^0$ so
$pp^* = pq^*qp^* $ is a maximal idempotent in $LI(\Delta) \setminus
\Delta^0$. Then from  Lemma \ref{maxlidemp}(b) we
see that  there exists an NE path $p_1$
and an edge $\check{e}$ for which $s(\check{e})$ has out-degree at
least $2$ such that $p p^* = p_1 \check{e} \check{e}^* p_1^*$ in
$LI(\Delta)$. By Lemma \ref{maxlidemp}(b) we have $p$ is not a
prefix of $p_1$ since $p_1$ is an NE path and so, again Lemma
\ref{maxlidemp}(b), $p = p_1 \check{e} p_2$ where $p_2$ is an NE
path in $\Delta$. Moreover, we have $e = \phi^{-1}(p_1)
\phi^{-1}(\check{e}) \phi^{-1}(p_2 q^*)$. From part (c) of this
lemma, we observe that $\phi^{-1}(\check{e}) = (\phi^{-1}(p_1))^* e
(\phi^{-1}(p_2 q^*))^*$. This forces the existence of NE paths $q_1,
q_2, q_3$ such that $\phi^{-1}(\check{e}) = q_1 q_2 e q^*_3$

(e) Take $v \in \Gamma^0$ such that $s^{-1}(v) =\{e_1, \cdots,
e_n\}$ with $n \geq 2$. According to part (d) of the lemma, we
have $\phi(e_i) = p_{i,1} \check{e}_i p_{i,2}
p^*_{i,3}$ for some NE paths $p_{i,j}$. Then since  $s(p_{i,1}) = \phi(v)$ for
all $i$  and since each $p_{i,1}$ is an NE path, we see that all
$p_{i,1}$ are the same path.  Since $p_{i,1}\check{e}_i
\check{e}_i^*p_{i,1}^*$ is the image in $LI(\Delta)$ of $e_ie_i^*$ under
$\phi$ and the $e_i$ are distinct, it follows that the $\check{e}_i$
are distinct, for $i = 1,..., n$. Hence the out-degree of
$r(p_{i,1}) = s(\check{e}_i)$ is at least $n$, which is the out-degree of $v = s(e_i)$. Similarly, from the second statement of part (d), we see that the
out-degree of $s(e_i)$ is less than or equal to the out-degree of
$s(\check{e}_i)$. It follows that
$s^{-1}(r(p)) =\{\check{e}_1, \cdots, \check{e}_n\}$. \QED

\medskip

We are now in a position to prove the following theorem.

\begin{Theorem} \label{algiso}
Let $\Gamma$ and $\Delta$ be connected graphs and let $F$ be a
field. If $LI(\Gamma) \cong LI(\Delta)$, then $L_F(\Gamma) \cong
L_F(\Delta)$.
\end{Theorem}

\noindent {\bf Proof.} By the definition of a Leavitt path algebra
we observe that $L_F(\Gamma)$ is isomorphic to the quotient of the
contracted semigroup algebra $F_0LI(\Gamma)$ of $LI(\Gamma)$ by the
ideal $I_1$ generated by elements of the form $\Sigma_{e \in
s^{-1}(v)} ee^* - v$ for $v \in \Gamma^0$ with the out-degree of $v$
at least 2. $L_F(\Delta)$ is isomorphic to the contracted semigroup
algebra $F_0LI(\Delta)$ of $LI(\Delta)$ by the ideal $I_2$ generated
by elements of the form $\Sigma_{d \in s^{-1}(u)} dd^* - u$ for $u
\in \Gamma^0$ with the out-degree of $u$ at least 2. Suppose that
$\phi$ is an isomorphism from $LI(\Gamma)$ onto $LI(\Delta)$. Then
$\phi$ induces an algebra isomorphism, say $\eta$, from
$F_0LI(\Gamma)$ onto $F_0LI(\Delta)$. Now for any $v \in \Gamma^0$
with out-degree greater than $1$ and any $e_i \in s^{-1}(v)$ we see
from Lemma~\ref{levisolem}(d), (e) that there exist NE paths $p,
p_i, q_i$ and edges $\check{e}_i \in s^{-1}(r(p))$  such that
$\phi(e_i) = p \check{e}_i p_i q_i^*$, $\phi(v) = s(p)$ and
$|s^{-1}(v)| = |s^{-1}(r(p))|$.
Distinct $e_i$ correspond to distinct $\check{e}_i$. Thus,
\begin{align*}
\eta(\Sigma_{e_i \in s^{-1}(v)} e_ie_i^* - v) & = \Sigma_{e_i \in s^{-1}(v)}\phi(e_i) (\phi(e_i))^* - s(p) \\
& = \Sigma_{\check{e}_i \in s^{-1}(u)} p \check{e}_i \check{e}_i^* p^* - p p^* \\
& = p(\Sigma_{\check{e}_i \in s^{-1}(u)} \check{e}_i \check{e}_i^* -
u) p^* \in I_2
\end{align*}
which means that $\eta(I_1) \subseteq I_2$. Similarly, one can
obtain that $\eta^{-1}(I_2) \subseteq I_1$. So we have $\eta(I_1) =
I_2$ and $\eta^{-1}(I_2) = I_1$. It follows that $L_F(\Gamma) \cong
L_F(\Delta)$ as required. \QED

\medskip

We remark that  Ruy Exel outlined an
alternative  proof (also suggested by Benjamin Steinberg) of Theorem
\ref{algiso} using his notion \cite{exel} of  tight representations
of  inverse semigroups.

\section{Some structural properties of  Leavitt inverse semigroups}

In this section   we determine some structural properties of Leavitt inverse semigroups culminating in a description of necessary and sufficient conditions for
two graphs to have isomorphic Leavitt inverse semigroups (Theorem
\ref{leavinvclass}) and some applications of that theorem to the structure of Leavitt path algebras for some classes of graphs. We will need some preliminary concepts and
lemmas in order to formulate and prove this result and some other structural properties.

Let $\Gamma$ be an arbitrary (directed) graph. Define a relation $\sim$ on $\Gamma^0$ by $v_1 \sim v_2$ if there
exist
 (possibly empty) NE paths $p$ and $q$
such that $s(p) = v_1, s(q) = v_2$ and $r(p) = r(q)$. Note that this
implies that $v_i \sim r(p) = r(q)$ for $i = 1,2$ even if the out-degree of $r(p)$ is at least 2 since the empty path at $r(p)$ is an NE path.

\begin{Lemma}
\label{equiv}

The relation $\sim$ is an equivalence relation on $\Gamma^0$.

\end{Lemma}

\noindent {\bf Proof.} The relation $\sim$ is reflexive since we
regard the empty path at any vertex $v \in \Gamma^0$ to be an NE
path. Clearly $\sim$ is symmetric. If $v_1 \sim v_2$ and $v_2 \sim
v_3$ then there are NE paths $p_1,q_1,p_2,q_2$ such that $s(p_1) =
v_1, \,  s(q_1) = v_2, \, r(p_1) = r(q_1), \, s(p_2) = v_2, \,
s(q_2) = v_3$ and $r(p_2) = r(q_2)$. If $q_1$ is the empty path then
$r(p_1) = v_2$ and so  $p_1p_2$ is an NE path with $s(p_1p_2) = v_1$
and $r(p_1p_2) = r(q_2)$, so in this case $v_1 \sim v_3$. Similarly
$v_1 \sim v_3$ if $p_2$ is the empty path. If neither $q_1$ nor
$p_2$ is the empty path then since all vertices in an NE path
(except the range vertex) have out-degree $1$ it follows that either
$q_1$ is a prefix of $p_2$ or $p_2$ is a prefix of $q_1$. In the
first case, there is an NE path $t$ such that $s(t) = r(q_1),\, r(t)
= r(q_2)$ and $p_2 = q_1t$, so $p_1t$ is an NE path with $s(p_1t) =
v_1$ and $r(p_1t) = r(q_2)$, and so $v_1 \sim v_3$. Similarly $v_1
\sim v_3$ if $p_2$ is a prefix of $q_1$. Hence $\sim$ is transitive.
\QED

\begin{Cor} \label{simloop}
If $e \in \Gamma^1$ is an edge that has exits then  $s(e) \sim
r(e)$ if and only if $e$ lies in a cycle which has exits only at
$s(e)$.
\end{Cor}

\noindent {\bf Proof.} If $s(e) \sim r(e)$, then there exist NE
paths $p, q$ such that $r(e) = s(p), s(e) = s(q)$ and $r(p) = r(q)$.
This forces that $q$ is trivial since $q$ has no exit and $s(e)$ has
out-degree greater than $1$. So the path $ep$ is a cycle which has
exits only at $s(e)$. The converse part is clear. \QED

\medskip

The equivalence relation $\sim$ enables a description of the Green relations on $LI(\Gamma)$.

\begin{Theorem}
\label{green}

Let $\Gamma$ be a graph and $pq^*, xy^*$ elements of $LI(\Gamma)$ in canonical form as described in Theorem~$\ref{leavittinv}$. Then the Green relations on $LI(\Gamma)$ are described as follows.

\medskip

 $(a)$ $pq^* \, \mathcal R \, xy^*$ iff $pp^* = xx^*$.

 $(b)$ $pq^* \, \mathcal L \, xy^*$ iff $qq^* = yy^*$.

 $(c)$ $pq^* \, \mathcal D \, r(p)$.

 $(d)$ $pq^* \, \mathcal D \, xy^*$ iff $r(p) \sim r(x)$.

 $(e)$ $pq^* \, \mathcal J \, xy^*$ iff there exist vertices $u, v  \in \Gamma^0$ with
$r(p) \sim u$ and $ r(x) \sim v$ such that $u$ and $v$ are in the same strongly connected component of $\Gamma$.

 $(f)$ If $pq^* \, \mathcal H \, xy^*$ and $pq^* \neq xy^*$ in $LI(\Gamma)$ then either there is a $($possibly empty$)$  NE path $p'$ from $r(p)$ to $r(x)$ and a non-trivial  NE cycle $C$ based at $r(x)$ or there is a $($possibly empty$)$ NE path $p'$ from $r(x)$ to $r(p)$ and a non-trivial  NE cycle $C$ based at $r(p)$. In the former case
 $xy^* = pp'C^np'^*q^*$ for some non-zero integer $n $
 $($where $C^{-n} $ is interpreted as $(C^*)^n$ for $n > 0)$: in the latter case $pq^* = xp'C^np'^*y^*$ for some non-zero integer $n $.

 $(g)$ The maximal subgroup of $LI(\Gamma)$ containing the idempotent $pp^*$ is either trivial or is isomorphic to the group $(\mathbb Z, +)$ of integers: it is non-trivial if and only if there is a path of the form $p'C$ where $s(p') = r(p), p'$ is a $($possibly trivial$)$ NE path and $C$ is a non-trivial NE cycle in $\Gamma$ based at $r(p')$.

\end{Theorem}

\noindent {\bf Proof.} Note that $pq^* \, \mathcal R \, xy^*$ iff $pq^*qp^* = xy^*yx^*$. The result of part (a)  follows since $pq^*qp^* = pr(q)p^* = pr(p)p^* = pp^*$ and similarly $xy^*yx^* = xx^*$. The proof of part (b) is similar. For part (c), note that $pq^* \, \mathcal R \, p$ by part (a). But $p^*p = r(p) = r(p)^*r(p)$, so $p \, \mathcal L \, r(p)$. Hence $pq^* \, \mathcal D \, r(p)$. Now suppose that $r(p) \sim r(x)$. Then there exist NE paths $p_1,p_2$ with $s(p_1) = r(p), s(p_2) = r(x)$ and $r(p_1) = r(p_2)$. Since $p_1p_1^* = r(p)$ and $p_1^*p_1 = r(p_1)$ it follows that $r(p) \, \mathcal D \, r(p_1)$. Similarly $r(x) \, \mathcal D \, r(p_2) = r(p_1)$, so $r(p) \, \mathcal D \, r(x)$ and hence by part (c) of this theorem, $pq^* \, \mathcal D \, xy^*$. Conversely, suppose that  $pq^* \, \mathcal D \, xy^*$, so $r(p) \, \mathcal D \, r(x)$, again by part (c). Then there exists $p_1q_1^*$ in canonical form such that $r(p) \, \mathcal R \, p_1q_1^* \, \mathcal L \, r(x)$. This implies that $r(p) = p_1p_1^*$ and $r(x) = q_1q_1^*$ by parts (a) and (b) of this theorem. Then by Corollary \ref{levinvidemp} $p_1$ and $q_1$ are NE paths, so  $r(p) \sim r(x)$. This proves part (d).

Now suppose  that there are vertices $u,v$  satisfying the conditions in part (e).
By Corollary 2 of \cite{MeMi} we know that $u \, \mathcal J \, v$ in $I(\Gamma)$, so $u \, \mathcal J \, v$ in $LI(\Gamma)$. But also by part (d), $r(p) \, \mathcal D \, u$ and $v \, \mathcal D \, r(x)$ so $r(p) \, \mathcal J \, r(x)$ in $LI(\Gamma)$, whence $pq^* \, \mathcal J \, xy^*$ by part (c).  Suppose conversely that $pq^* \, \mathcal J \, xy^*$. Then $r(p) \, \mathcal J \, r(x)$ by part (c). So  there exist $p_1q_1^*, p_2q_2^*$ in canonical form such that $r(p) = p_1q_1^*r(x)p_2q_2^*$. This forces $s(p_1) = s(q_2) = r(p)$ and $s(q_1) = s(p_2) = r(x)$. Also, either $p_2$ is a prefix of $q_1$ or $q_1$ is a prefix of $p_2$. Suppose that $q_1$ is a prefix of $p_2$. So there exists a (possibly empty) directed path $t_1$ with $p_2 = q_1t_1$. Also since $q_2$ is a directed path from $r(p)$ to $r(p_2) = r(q_2) = r(t_1)$ there exist (possibly empty) directed paths $t_2,t_3,t_4$ such that $p_1 = t_2t_3$ and $q_2 = t_2t_4$. Then $r(p) = p_1q_1^*r(x)p_2q_2^* = t_2t_3q_1^*q_1t_1t_4^*t_2^* = p_3q_3^*$ where $p_3 = t_2t_3t_1$ and $q_3 = t_2t_4$. This forces $p_3 = q_3$ to be an NE path by Corollary \ref{levinvidemp}, and hence $t_4 = t_3t_1$ and also $r(p) \sim u' = r(p_3)$ and $p' = p_2$ is a directed path from $r(x)$ to $u'$. A similar argument applies in the case where $p_2$ is a prefix of $q_1$. Similarly, there is some vertex $v'$ with $r(x) \sim v'$ and a directed path $p''$ from $r(p)$ to $v'$. Thus  in all cases we have some vertices $u',v'$ with $u' \sim r(p), v' \sim r(x)$ and directed paths $p'$ from $r(x)$ to $u'$ and $p''$ from $r(p)$ to $v'$.

Since $r(p) \sim u'$, there are NE paths $h_1,h_2$ with $s(h_1) = u', s(h_2) = r(p)$ and $r(h_1) = r(h_2) \sim r(p)$. Since $h_2$ is an NE path starting at $r(p)$ it must be a prefix of $p''$, so there exists a directed path $h_3$ such that $p'' = h_2h_3$. Denote the vertex $r(h_1) = r(h_2) = s(h_3)$ by $u$. Similary, there are directed paths $h_4, h_5, h_6$ such that $h_4,h_5$ are NE, $s(h_4) = v', s(h_5) = r(x), r(h_4) = r(h_5) = s(h_6)$ and $p' = h_5h_6$. Denote the vertex $s(h_6) = r(h_4) = r(h_5)$ by $v$. Then $u \sim r(p), v \sim r(x), h_3h_4$ is a directed path from $u$ to $v$ and $h_6h_1$ is a directed path from $v$ to $u$. Thus $u$ and $v$ are in the same strongly connected component of $\Gamma$. This proves part (e).

Suppose that  $pq^* \, \mathcal H \, xy^*$ and $pq^* \neq xy^*$. Then by parts (a) and (b), $pp^* = xx^*$ and $qq^* = yy^*$ and also either $x \neq p$ or $y \neq q$. Assume that $x \neq p$. (The case $y \neq q$ is similar.) By Corollary \ref{levinvidemp} there is a non-empty NE path $t$ such that either $p = xt$ or $x = pt$. Assume that $x = pt$ since the other case is dual. Since $t$ is an NE path we must have $t = p'C^k$ for some  NE path $p'$ containing no cycles, some NE  cycle $C$, and some integer $k \geq 0$. Since $x \neq p$ we cannot have $p'$ and $C$ both trivial: also, if $p'$ is trivial then $k > 0$.

Case 1: $r(p') = r(p)$. Then $C$ is a cycle based at $r(p) = r(q)$ and $x = pC^k$ for some $k > 0$.  Since $yy^* = qq^*$, Corollary \ref{levinvidemp} implies that there is an NE path $p''$ such that either $y = qp''$ or $q = yp''$. Since $r(q) = r(x) = r(y)$ and $C$ is an NE cycle, this forces $p'' = C^m$ for some $m \geq 0$. If $y = qC^m$ then since the last edge in $C$ is an NE edge, the fact that $xy^*$ is in canonical form forces $m = 0$. So in this case $xy^* = pC^kq^*$. If $q = yC^m$ then $xy^* = pC^kC^{-m}q^* = pC^{k-m}q^*$ since $CC^* = C^*C = s(C)$ in $LI(\Gamma)$. Thus in Case 1, $xy^* = pC^nq^*$ for some non-zero integer $n$.

Case 2: $r(p') \neq r(p)$.  Then $r(x) = r(y) = r(p') \neq r(p) = r(q)$. As in Case 1, there is an NE path $p''$ such that either $y = qp''$ or $q = yp''$. If $q = yp''$ then the path $p'p''$ is a non-trivial  NE circuit based at $r(p)$ so there is some non-trivial cycle $C_1$ based at $r(p)$ and $p'p'' = C^n$ for some $n > 0$. Since $p''$ is an NE path from $r(y) = r(x)$ to $r(p) = r(q) $ in this case we have $p''p''^* = r(y)$ and so $y = yr(y) = yp''p''^* = qp''^*$. Then $xy^* = pp'p''q^* = pC_1^nq^*$. If $y = qp''$ then $p'C^k$ and $p''$ are non-trivial NE paths starting at $r(p) = r(q)$ and ending at $r(p') = r(p'') = r(x) = r(y)$. So $p'' = p' C^l$ for some integer $l \geq 0$. The fact that $xy^*$ is in canonical form forces that one and only one of $k, l$ is nonzero. Thus, we have $xy^* = pp'C^np'^*q^*$ for some non-zero integer $n$. This completes the proof of part (f).

\medskip

Suppose that $xy^* \, \mathcal H \, pp^*$ and $xy^* \neq pp^*$. From part (f) we have an NE path $p'$ and a non-trivial NE cycle $C$ either with $s(p') = r(p), s(C) = r(p') = r(x) = r(y)$ and $ xy^* = pp'C^np'^*p^*$ or with $s(p') = r(x) = r(y), s(C) = r(p') =r(p)$ and $pp^* = xp'C^np'^*y^*$ for some non-zero $n$. This latter condition is impossible by Corollary  \ref{levinvidemp} and the uniqueness of canonical forms, so we must have the former condition. Without loss of generality, we may suppose that $p'$ does not contain an edge in $C$. Otherwise we may assume that $p' = p_1 p_2$ where $p_1$ does not contain an edge in $C$ and all edges in $p_2$ are contained in $C$, and also that $C = p_3 p_2$. Thus, $xy^* = p p_1 p_2 (p_3 p_2)^n p_2^* p_1^* p^* = p p_1 (p_2 p_3)^n p_1^* p^*$ where $p_1$ does not contain an edge in the cycle $p_2 p_3$ which is a conjugate of $C$. Moreover, if $pp'C^mp'^*p^* = pp''(C')^np''^*p^*$, then since $p', p''$ and $C, C'$ are NE, $C'$ must be a conjugate of $C$. Noticing that any edge contained in $C'$ is also contained in $C$, we see that $p' = p''$ and $C^n = C^m$ so that $C^{m-n} = r(C)$. By Corollary \ref{levinvidemp} this implies $m = n$. Since $(pp'C^mp'^*p^*)(pp'C^np'^*p^*) = pp'C^{m+n}p'^*p^*$, we see that the $\mathcal H$-class of $pp^*$ is isomorphic to the group $(\mathbb Z,+)$. Thus part (g) is verified.  \QED

\medskip

\noindent {\bf Remark} We remark that there is a significant difference between the Green relations on $I(\Gamma)$ and the Green relations on $LI(\Gamma)$. The Green relations on $I(\Gamma)$ are given in \cite{MeMi}, Corollary 2. In particular, $I(\Gamma)$ is combinatorial  for every graph $\Gamma$, but by Theorem \ref{green}(f) $LI(\Gamma)$ is  combinatorial if and only if  $\Gamma$   has no non-trivial NE cycles. In particular, if $\Gamma$ is acyclic (which is equivalent to the multiplicative semigroup of $LI(\Gamma)$ being von-Neumann regular by \cite{AR}), then $LI(\Gamma)$ is combinatorial.  The converse is false in general of course since $\Gamma$ may have non-trivial cycles but no non-trivial NE cycles.

\medskip

Denote the $\sim$-equivalence class of a vertex $v \in \Gamma^0$ by
$[v]$. Note that $[v] = \{v\}$ if and only if   $v$ does not have out-degree $1$
and $s(e) $ does not have out-degree $1$  for every edge $e \in r^{-1}(v)$. (In
particular, $[v] = \{v\}$ if $v$ is a source whose out-degree is not
$1$.) We denote by $\Gamma_{[v]}$ the  subgraph of $\Gamma$ induced
by the set of vertices in $[v]$. That is, $\Gamma_{[v]}^0 = [v]$ and
$\Gamma_{[v]}^1 = \{e \in \Gamma^1 : s(e),r(e) \in [v]\}$.

\begin{Lemma} \label{equivdis}
Let $v$ be a vertex in a graph $\Gamma$.

$(a)$ If there are at least two non-conjugate cycles $C_1$ and $C_2$
in $\Gamma_{[v]}$   then $C_1^0 \cap C_2^0$ contains a vertex of
out-degree greater than $1$.

$(b)$ $\Gamma_{[v]}$ contains at most one vertex $w$ of out-degree not
equal to $1$. This vertex $w$ is contained in every cycle $C$ for
which $C \subseteq \Gamma_{[v]}$ $($if there are any such cycles$)$.

$(c)$ If $\Gamma_{[v]}$ contains an NE cycle, then every vertex in
$[v]$ has out-degree $1$. In particular, there is only one conjugacy
class of cycles in $\Gamma_{[v]}$.


\end{Lemma}

\noindent {\bf Proof.} (a) Suppose that $\Gamma_{[v]}$ contains
distinct cycles $C_1$ and $C_2$ that are not cyclic conjugates of
each other. Let $v_1$ be a vertex in $C_1 \setminus C_2$ and $v_2$ a
vertex in $C_2 \setminus C_1$. Then $v_1 \sim v_2$ so there are NE
paths $p$ and $q$ with $s(p) = v_1, \, s(q) = v_2$ and $r(p) = r(q)
\in (C_1 \cap C_2)^0$. If the out-degree of $r(p)$ is greater than
$1$ we are done. If not, then $r(p)$ has out-degree $1$ and the edge $e _1
= e_{r(p)}$ lies in $C_1 \cap C_2$. But then either $r(e_1)$ has
out-degree at least $2$ or $r(e_1)$ has out-degree $1$, and in the latter
case the edge $e_2$ starting at $r(e_1)$ lies on $C_1 \cap C_2$.
Continuing in this fashion we see that there is some vertex $w$ in
$C_1^0 \cap C_2^0$ with out-degree at least $2$.

(b) If $v_1$ and $v_2$ are distinct vertices  in $[v]$ then $v_1
\sim v_2$. By the definition of the equivalence relation $\sim$ this
forces either $v_1 = v_2$ or else at least one of the vertices $v_i$
has out-degree $1$. So there is at most one vertex in $[v]$ of
out-degree not equal to $1$. Suppose that there is such a vertex in
$[v]$ and denote it by $w$. If $C$ is a cycle in $\Gamma_{[v]}$ and
$v_1$ is a vertex in $C$, then $v_1 \sim w$ so there are NE paths
$p$ and $q$ with $s(p) = v_1, \, s(q) = w$ and $r(p) = r(q)$. Since
the out-degree of $w$ is not $1$, $q$ must be the empty path at $w$
and $w = r(p)$. But then since every vertex in $p$ except $w$ has out-degree $1$
 this forces $w$ to be a vertex of the cycle $C$. If $C_1$
and $C_2$ are distinct non-conjugate cycles in $\Gamma_{[v]}$  then
by the proof  above  we see that $w$ is in both cycles.

(c) Suppose that $\Gamma_{[v]}$ contains an NE cycle $C$. If $[v]$
contains a  vertex  of out-degree not equal to $1$ then this vertex
must lie on $C$ by part (b), but this is a contradiction since $C$
is an NE cycle. The fact that there is only one conjugacy class of
cycles in $\Gamma_{[v]}$ follows from part (a). \QED

\begin{Lemma}
\label{equivstructure}

Let $v$ be a vertex of the graph $\Gamma$. Then

$(a)$ If every vertex of $[v]$ has out-degree $1$ then  $\Gamma_{[v]}$
is a directed cover of $B_{\{a\}}$, in which case it is either an
infinite tree or has structure determined by Theorem~$\ref{out1}(b)$;

$(b)$ If $[v]$ has a sink in $\Gamma$,  then $\Gamma_{[v]}$ is an
immersion over $B_{\{a\}}$ whose structure is determined by Theorem~$\ref{out1}(a)$;

$(c)$ If $[v]$ has a vertex $w$ of out-degree greater than $1$ in
$\Gamma$ then either $\Gamma_{[v]}$ contains no cycles, in which
case $w$ is a sink of the graph $\Gamma_{[v]}$ and $\Gamma_{[v]}$ is
an immersion over $B_{\{a\}}$ whose structure is determined by
Theorem~$\ref{out1}(a)$, or else $\Gamma_{[v]}$ has at least one cycle
and $w \in C^0$ for every cycle $C$ in $\Gamma_{[v]}$. In the latter
case, if $v'$ is any vertex of $[v]$ with $v' \neq w$, then there is
a unique directed path from $v'$ to $w$ that does not include a
cycle in $\Gamma_{[v]}$ as a subpath.

\end{Lemma}

\noindent {\bf Proof.} Part (a) is immediate from Theorem
\ref{out1}. If $[v]$ contains a sink of $\Gamma$ then this vertex is
also a  sink of $\Gamma_{[v]}$ so the result of part (b) is also
immediate from Theorem \ref{out1}. Suppose that $[v]$ has a vertex
$w$ of out-degree greater than $1$ in $\Gamma$. There is a unique
such vertex $w$ by Lemma \ref{equivdis}(b). If $\Gamma_{[v]}$
contains no cycles, then  $w$ is a sink in the graph $\Gamma_{[v]}$,
so $\Gamma_{[v]}$ has the structure described in Theorem
\ref{out1}(a). If $\Gamma_{[v]}$ has  at least one cycle then $w
\in C^0$ for every cycle $C$ in $\Gamma_{[v]}$ by Lemma
\ref{equivdis}(b).  If $v' \neq w$ is a vertex of $[v]$  then there
is a directed path $p$ from $v'$ to $w$ since $v' \sim w$ and $w$
has out-degree greater than $1$. Since $w$ is in every cycle in
$\Gamma_{[v]}$, we may assume that the path $p$ does not contain any
cycle in $\Gamma_{[v]}$ as a subpath. The uniqueness of such a
directed path $p$ follows by an argument very similar to the
argument used in the proof of Theorem \ref{out1}. \QED

\begin{Lemma}
\label{equalequiv} If $\phi$ is an isomorphism from $LI(\Gamma)$
onto $LI(\Delta)$, then the following statements hold.

$(a)$ $\phi$ induces a bijection from $\Gamma^0$ onto $\Delta^0$.

$(b)$ For all vertices $v_1,v_2$ of $\Gamma$, $v_1 \sim v_2$ if and
only if $\phi(v_1) \sim \phi(v_2)$.

$(c)$ $\phi$ induces a bijection of the equivalence class $[v]$ in
$\Gamma^0$ onto the equivalence class $[\phi(v)]$ in $\Delta^0$ for
all $v \in \Gamma^0$.

$(d)$ For each integer $n \geq 1$, $[v]$ contains a vertex of
out-degree $n$ if and only if $[\phi(v)]$ contains a vertex of
out-degree $n$.

$(e)$ If $C$ is a cycle in $\Gamma_{[v]}$, then $\phi(C)$ is uniquely
expressible in the form $\phi(C) = pC'p^*$ or $\phi(C) = pC'^*p^*$
in $LI(\Delta)$ for some cycle $C'$ and some NE path $p$ contained
in $\Delta_{[\phi(v)]}$, and moreover $\phi^{-1}(C') = p_1 C_1
p_1^*$ or $\phi^{-1}(C') = p_1 C_1^* p_1^*$ for some cyclic
conjugate $C_1$ of $C$ and some NE path $p_1$ contained in
$\Gamma_{[v]}$.

$(f)$ $\phi$ induces a bijection between the set of distinct conjugacy
classes of cycles in $\Gamma_{[v]}$ and the set of distinct
conjugacy classes of cycles in $\Delta_{[\phi(v)]}$.


\end{Lemma}

\noindent {\bf Proof.} (a) By Lemma \ref{levisolem}(a), $\phi$ maps
vertices of $\Gamma$ to vertices of $\Delta$. The restriction of
$\phi$ to $\Gamma^0$ is clearly injective since $\phi$ is injective.
But by Lemma \ref{levisolem}(a), the map $\phi^{-1}$ maps vertices
of $\Delta$ to vertices of $\Gamma$, so the restriction of $\phi$ to
$\Gamma^0$ is a bijection onto $\Delta^0$.

(b) This follows immediately from Theorem~\ref{green}(d).

(c) By parts (a) and (b) of this lemma, $\phi$ induces an injection
of the equivalence class $[v]$ into the equivalence class
$[\phi(v)]$ for all $v \in \Gamma^0$. This map is surjective since
$\phi^{-1}$ induces an injection of the equivalence class
$[\phi(v)]$ into the equivalence class $[v]$.

(d)  Suppose that $[v]$ contains a
vertex $v$ of out-degree $n > 1$. Then, in the notation of
Lemma \ref{levisolem}(e), the vertex $s(p_{i,1})$ has out-degree
$n$ and since $p_{i,1}$ is an NE path, $s(p_{i,1}) \sim
\phi(v)$. If all vertices of $[v]$ have out-degree $1$ then by what we just proved, applied to $\phi^{-1}$,  all vertices of
$[\phi(v)]$ have out-degree $1$.

(e) Suppose that $\Gamma_{[v]}$ contains a cycle $C = e_1e_2...e_n$.
Let $\phi(e_i) = p_iq_i^*$  where $p_i,q_i$ are directed paths in
$\Delta$ with $r(p_i)=r(q_i)$. If $s(e_i)$ has out-degree $1$ then by
Lemma \ref{levisolem}(c), $p_i,q_i$ are NE paths so all of their
vertices are related via the equivalence relation $\sim$ on
$\Delta^0$. By Lemma \ref{equivdis}(b) $C$ contains at most one
vertex (say $s(e_k)$) whose out-degree is greater than $1$. Then by
Lemma \ref{levisolem}(d), $\phi(e_k) = p'_k\check{e}_kp''_kq_k^*$
where $p'_k,p''_k,q_k$ are NE paths and the out-degree of
$s(\check{e}_k)$ is at least $2$. Then all vertices in $p'_k$ are
$\sim$-related to $r(p'_k) = s(\check{e}_k)$ and all vertices in
$p''_kq_k^*$ are $\sim$-related to $r(\check{e}_k)$. Then since
$s(q_1) = s(p_2), ..., s(q_{k-1}) = s(p_k), s(q_k) =
s(p_{k+1}),...,s(q_n) = s(p_1)$, we see that all vertices in
$p_1q_1^*p_2q_2^*...p'_k$ are $\sim$-related to $s(p_1)$ and all
vertices in $p''_kq_k^*p_{k+1}q_{k+1}^*...q_n^*$ are $\sim$-related
to $s(q_n) = s(p_1)$. Thus all vertices in $\phi(C)$ are in the
$\sim$-class $[\phi(v)]$ of $\Delta^0$, that is $\phi(C)$ is a path in
$\Delta_{[\phi(v)]}$.

Since $C$ represents a non-zero element of $LI(\Gamma)$, we have
$\phi(C) = pq^*$ in $LI(\Delta)$, where $p$ and $q$ are directed
paths in $\Delta$ with $s(p) = s(q) = \phi(s(e_1))$ and $r(p) =
r(q)$. Furthermore, all vertices in $pq^*$ are in $[\phi(v)]$ by the
argument in the previous paragraph since these vertices are among
the vertices in the union of the paths $p_iq_i^*, \, i = 1,...,n$.
Since $C$ is a cycle, $C^2$ is also a non-zero element of
$LI(\Gamma)$ and so $pq^*pq^*$ is a non-zero element in
$LI(\Delta)$. Hence either $p$ is a prefix of $q$ or $q$ is a prefix
of $p$. Note from Lemma~\ref{levisolem}(b) and the multiplication in
$LI(\Gamma)$ that $q$ is NE. If the out-degree of $s(e_k)$ is at
least 2, then by Lemma~\ref{levisolem}(d) $p$ contains the unique
edge $\check{e}_k$ which has exits so that $q$ must be a prefix of
$p$. That is $p = q p_1$ for some directed circuit $p_1$ which
contains $\check{e}_k$. Now $p_1$ must be a cycle since only
$\check{e}_k$ has exits and $\check{e}_k$ appears in $p_1$ only
once. If $C$ is NE and $q$ is a prefix of $p$ which means that $p =
qp_2$ for some directed circuit $p_2$ in $\Delta$, then by
Lemma~\ref{levisolem}(c) $p_2$ is NE. So there must exist some NE
cycle $C'$ such that $p_2 = (C')^k$ for some positive integer $k$.
A similar discussion shows that $\phi^{-1}(q C' q^*) = C^l$ for some
positive integer $l$. These force that $k = l = 1$. If $p$ is a
prefix of $q$, a similar argument shows that $\phi(C) = p C'^* p^*$
for some NE cycle $C'$. The uniqueness of such expression follows
from the canonical forms of elements in $LI(\Delta)$.

Moreover, if $\phi(C) = p C' p^*$ for some NE cycle $C'$ and NE path
$p$, then we get $C = \phi^{-1}(p) \phi^{-1}(C') \phi^{-1}(p^*)$. We
see from Lemma~\ref{levisolem}(c) that there exist NE paths $p_1,
q_1$ in $\Gamma$ such that $C = (q_1 p_1^*) \phi^{-1}(C') (p_1
q_1^*)$. That is, $\phi^{-1}(C') = p_1 q_1^* C q_1 p_1^*$. It
follows that $q_1^* C q_1$ is a cyclic conjugate of $C$ since $q_1$
is NE. A similar argument applies if $\phi(C) = p C'^* p^*$ for some
NE cycle $C'$ and NE path $p$.


(f) By Lemma \ref{equivstructure}(a), if $\Gamma_{[v]}$ is not a
tree and all vertices of $\Gamma_{[v]}$ have out-degree $1$ then
$\Gamma_{[v]}$ contains a unique cycle $C_v$ (up to cyclic
conjugates). By part (c) of this lemma and Lemma \ref{levisolem}(c)
all vertices of $\Delta_{[\phi(v)]}$ have out degree $1$ and
$\Delta_{[\phi(v)]}$ has a unique cycle $C'_{\phi(v)}$ (up to cyclic
conjugates). If $\Gamma_{[v]}$ has $n$ distinct cycles $C_1,...,C_n$
(up to cyclic conjugates) for some $n > 1$, then by Lemma
\ref{equivdis}(b) there exists some unique vertex $w$ (in $[v]$)
whose out-degree is at least 2 in $\Gamma$ such that $w$ lies in all
these cycles. Moreover, these cycles correspond to the edges in
$s^{-1}(w) \cap \Gamma_{[v]}^1$. Then by part (b) of this lemma and
Lemma~\ref{levisolem}(e), $\Delta_{[\phi(v)]}$ has $n$ distinct
cycles (up to cyclic conjugates). \QED

\medskip

For each vertex $v$ of a graph $\Gamma$ let $T_{[v]}$ be a
(directed) spanning tree of the subgraph $\Gamma_{[v]}$. From the
structure of the graphs $\Gamma_{[v]}$ described in Lemma
\ref{equivstructure} it is clear that $T_{[v]} = \Gamma_{[v]}$ if
$[v]$ does not contain any cycle, while on the other hand if
$\Gamma_{[v]}$ does contain a cycle  then $T_{[v]}$ contains all
edges of $\Gamma_{[v]}$ except some particular edge $e_C$ of $C$ for
each cycle $C$ in $\Gamma_{[v]}$. Note that these $e_C$'s can not be
chosen arbitrarily for cycles with exits. For instance, consider the
graph $\Gamma$ represented in Diagram 7.1.
The subgraph $\Gamma_{[v]}$ has two conjugacy classes of cycles,
namely the conjugacy classes of the cycles $C_1$  consisting of
the edges $e_1, e_2, e_3, e_5$, and  $C_2$  consisting of the edges
$e_1, e_2, e_4, e_6$; we can not respectively choose $e_1$ as
$e_{C_1}$ and $e_2$ as $e_{C_2}$ because the remaining subgraph is
not a tree. On the other hand, $\Gamma_{[u]}$ has one conjugacy
class of cycles, namely the conjugacy class of the cycle $C_3$
consisting of the edges $e_7$ and $e_8$, and we can choose either
$e_7$ or $e_8$ as $e_{C_3}$.

\medskip

\begin{center}
\begin{tikzpicture}
\node at (0, 0) (v1) {\color{red}{$\bullet$}}; \node at (2, 0) (v2)
{\color{red}{$\bullet$}}; \node at (4, 0) (v3)
{\color{red}$\bullet$}; \node at (5.3, 1.5) (v4)
{\color{red}$\bullet$}; \node at (5.3, -1.5) (v5)
{\color{red}$\bullet$}; \node at (6, 0) (v6)
{\color{blue}$\bullet$}; \node at (8, 0) (v7)
{\color{blue}$\bullet$}; \node at (7.3, 1.5) (v8) {$\bullet$};
\draw[->, red] (v1) node[above] {$v$} to node[above] {$e_1$} (v2);
\draw[->, red] (v2) to node[above] {$e_2$} (v3); \draw[->, red] (v3)
to node[above,sloped] {$e_3$} (v4); \draw[->, red] (v3) to
node[below,sloped] {$e_4$} (v5); \draw[->, red] (v4) to [bend right
= 20] node[above,sloped] {$e_5$} (v1); \draw[->, red] (v5) to [bend
left = 20] node[below,sloped] {$e_6$} (v1); \draw[->] (v3) to (v6);
\draw[->] (v6) to (v8); \draw[->, blue] (v6) node[above] {$u$} to
[bend left = 30] node[below] {$e_7$} (v7); \draw[->, blue] (v7) to
[bend left = 30] node[below] {$e_8$} (v6);
\end{tikzpicture}

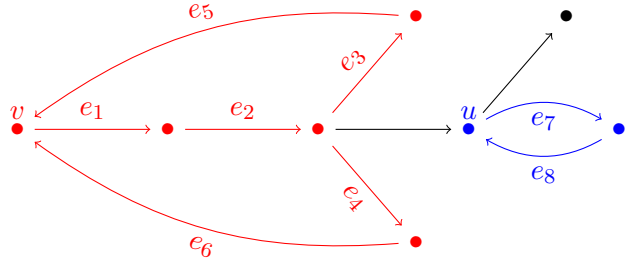
Diagram~7.1 \,\, Two $\sim$-classes containing cycle(s)
\end{center}

\medskip

Since cycles in distinct equivalence classes are clearly disjoint,
the choice of edges $e_C$ for cycles in distinct equivalence classes
are disjoint and the spanning trees $T_{[v]}$ (for $v \in \Gamma^0$)
are uniquely determined by the choice of these edges $e_C$ for each
cycle $C$.

\medskip

We call a set $\{ T_{[v]} : v \in \Gamma^0 \}$ of spanning trees for
the induced graphs $\Gamma_{[v]}$ a {\em set of NE spanning trees}
if every edge in each tree is an NE edge in $\Gamma$.

\begin{Lemma} \label{NEspanning}
For any $v \in \Gamma^0$ and any cycle $C$ in $\Gamma_{[v]}$, one
obtains a set of NE spanning trees by choosing any edge in $C$ as
$e_C$ if $C$ is an NE cycle and choosing the edge with exits as
$e_C$ if $C$ has exits. Every set of NE spanning trees is obtained
this way.
\end{Lemma}

\noindent {\bf Proof.} If every vertex of $C$ has out-degree $1$, then by
Lemma~\ref{equivdis} we know that $\Gamma_{[v]}$ has only one cycle.
It follows from Theorem~\ref{out1}(a) that the subgraph
$\Gamma_{[v]} \setminus \{e_C\}$ is an NE tree (which is a spanning
tree of $\Gamma_{[v]}$) since $s(e_C)$ is a sink and every other
vertex has out-degree 1. If $C$ has exits, then we see from
Lemma~\ref{equivdis}(b) that $C$ contains only one vertex $v_0$ with
out-degree at least 2. In this case, $\Gamma_{[v]}$ may contain more
than one cycle. Again it follows from Theorem~\ref{out1}(a) and
Lemma \ref{equivstructure}(c) that the subgraph $\Gamma_{[v]}
\setminus s^{-1}(v_0)$ is an NE tree (which is a spanning tree of
$\Gamma_{[v]}$) since $v_0$ is a sink and every other vertex has
out-degree 1. In this way, we get a set of NE spanning trees. By the
definition of a set of NE spanning trees, no  set of NE spanning
trees can contain an edge $e \in s^{-1}(v_0)$ since all of these
edges have exits in $\Gamma$. \QED

\medskip

Form a new graph $\bar{\Gamma}$ by contracting each spanning tree
$T_{[v]}, \, (v \in \Gamma^0)$ to a point. More precisely, we may
describe the graph $\bar{\Gamma}$ in the following way:
$\bar{\Gamma}^0 = \{[v] : v \in \Gamma^0\}$ and $\bar{\Gamma}^1$ is
a set in one-one correspondence with $\{ e \in \Gamma^1 \setminus
\bigcup_{v \in \Gamma^0}T_{[v]}^1\}$. We denote the image of $e$ under
this correspondence by $\bar{e}$.  The source and range functions
are defined for $\bar{e} \in \bar{\Gamma}^1$ by $s(\bar{e}) =
[s(e)]$ and $r(\bar{e}) = [r(e)]$. Thus the edge $e_C$ of a cycle
$C$ in $\Gamma_{[v]}$ gives rise to a loop $\bar{e}_C$ at $[v]$ in
$\bar{\Gamma}$. There is a natural function $\chi_{\Gamma} : \Gamma
\rightarrow \bar{\Gamma}$ defined by $\chi_{\Gamma}(v) = [v]$ for
all $v \in \Gamma^0$, $\chi_{\Gamma}(e) = \bar{e}$ for $e \in
\Gamma^1 \setminus \bigcup_{v \in \Gamma^0}T_{[v]}^1$, and
$\chi_{\Gamma}(e) = [v]$ if $e \in T_{[v]}^1$ for some $v \in \Gamma^0$.
The map $\chi_{\Gamma}$ is not a graph morphism since it maps some
edges to vertices.
\medskip

\begin{Lemma} \label{deg1contracted}
In a contracted graph $\bar{\Gamma}$, if the out-degree of $[v]$ is one, then
the only edge in $s^{-1}([v])$ is a loop. Hence, the equivalence relation $\sim$ on
$\bar{\Gamma}^0$ is trivial.
\end{Lemma}

\noindent {\bf Proof.} For any edge $\bar{e}$ in $\bar{\Gamma}$,
either $e$ does not belong to any graph $\Gamma_{[v]}$ or $e$ belongs
to exactly one such graph. In the former case the out-degree of
$s(e)$ is at least 2 and $s(e)$ is not $\sim$ related to $r(e)$. In
the latter case either  $e$ lies in a cycle which has exits only at
$s(e)$ or $e$ lies in an NE cycle. In the first two cases, the
out-degree of $s(\bar{e})$ is also at least 2. In the third case,
$\bar{e}$ is a loop and the out-degree of $s(\bar{e})$ is 1. The
second statement of the lemma thus follows directly from the
definition of the relation $\sim$. \QED
\medskip

We remark that the graph $\bar{\Gamma}$ and the contraction map
$\chi_{\Gamma}$ depends on the choice of the spanning trees
$T_{[v]}, \, v \in \Gamma^0$. However the contracted graphs
corresponding to different choices of NE spanning trees are
isomorphic.

\begin{Lemma}
\label{alternativetree} For any $v \in \Gamma^0$, arbitrarily choose
sets of NE spanning trees $T_{v}$ and $T'_{v}$ for $\Gamma_{[v]}$.
Then the contracted graph $\bar{\Gamma}_1$ corresponding to the
spanning trees $T_{v}, \, v \in \Gamma^0$, is isomorphic to the
contracted graph $\bar{\Gamma}_2$ corresponding to the spanning
trees $T'_{v}, \, v \in \Gamma^0$.
\end{Lemma}

\noindent {\bf Proof.} By Lemma \ref{NEspanning}
an NE spanning tree $T_{[v]}$ for $\Gamma_{[v]}$ is determined by
the choice of an edge $e_C$ in $C$ in an NE cycle $C$ in
$\Gamma_{[v]}$ (if such a  cycle exists) since the choices of the
edges $e_C$ for a cycle $C$ that has exits is fixed. Similarly
another spanning tree $T'_v$ for $\Gamma_{[v]}$ is determined by the
choice of another edge $e'_C$ in each NE cycle $C$
in $\Gamma_{[v]}$. Then the map defined by $[v] \mapsto [v]$, $\bar{e}_C
\mapsto \bar{e}'_C$ for all cycles $C$ in $\Gamma_{[v]}$ (and
all $v \in \Gamma^0$), and $\bar{e} \mapsto \bar{e}$ for all
other edges of $\bar{\Gamma}_1$ induces a graph isomorphism of
$\bar{\Gamma}_1$ onto $\bar{\Gamma}_2$. \QED

\begin{Lemma}
\label{chitilde}
The mapping $\chi_{\Gamma} : \Gamma \rightarrow
\bar{\Gamma}$ naturally induces a $0$-restricted morphism
$\tilde{\chi}_\Gamma$ from $LI(\Gamma)$ onto $LI(\bar{\Gamma})$.
\end{Lemma}

\noindent {\bf Proof.} This is a routine calculation. The map
$\chi_{\Gamma}$ defines a map $\tilde{\chi}_{\Gamma}$ from the
generators of $LI(\Gamma)$ onto the generators of $LI(\bar{\Gamma})$
that is easily seen to extend to a $0$-restricted morphism if we define
$\tilde{\chi}(0) = 0$. \QED

\medskip

We note that in general $\tilde{\chi}_{\Gamma}$ is not a
homomorphism since if $e_1$ and $e_2$ are distinct edges in one of
the spanning trees $T_{[v]}$ in $\Gamma$, then $e_1^*e_2 = 0$ in
$LI(\Gamma)$ but
$\tilde{\chi}_{\Gamma}(e_1^*)\tilde{\chi}_{\Gamma}(e_2) = [v][v] =
[v]$ in $LI(\bar{\Gamma})$. Note also that by the definition of the
edges in $\bar{\Gamma}$, $\tilde{\chi}_{\Gamma}(e) = \bar{e}$ if $e
\in \Gamma^1 \setminus \bigcup_{v \in \Gamma^0}T_{[v]}^1$, so
$\tilde{\chi}_{\Gamma}$  induces a bijection from the edges $e \in
\Gamma^1 \setminus \bigcup_{v \in \Gamma^0}T_{[v]}^1$ onto the edges
$\bar{e}$ in $\bar{\Gamma}$.

\medskip

\begin{Lemma}
\label{baredgeiso}
Let $\phi : LI(\Gamma) \rightarrow LI(\Delta)$ be a semigroup
isomorphism. Choose two sets of NE spanning trees $\{ T_{[v]} : v \in \Gamma^0 \}$
and $\{T'_{[u]} : u \in \Delta^0 \}$. Then
there exists a graph isomorphism $\bar{\phi}$ from $\bar{\Gamma}$ to
$\bar{\Delta}$ such that $\chi_\Delta \psi(v) = \bar{\phi}
\chi_\Gamma(v)$ where $v \in \Gamma^0$ and $\psi$ is the restriction
of $\phi$ to $\Gamma^0$.
\end{Lemma}

\noindent {\bf Proof.} For any $[v] \in \bar{\Gamma}^0$, define
$\bar{\phi}([v]) = \chi_\Delta \phi(v)$. We see from
Lemma~\ref{equalequiv}(b) that
$\bar{\phi}$ is well-defined and is an injection from
$\bar{\Gamma}^0$ to $\bar{\Delta}^0$.  For any $[u] \in
\bar{\Delta}^0$, again by Lemma~\ref{equalequiv}(b), $[\phi^{-1}(u)]
\in \bar{\Gamma}^0$ does not depend on the choice of $u$. Moreover,
$\bar{\phi}([\phi^{-1}(u)]) = \chi_\Delta \phi(\phi^{-1}(u)) =
\chi_\Delta \phi \phi^{-1}(u) = [u]$. So $\bar{\phi}$ is a bijection
from $\bar{\Gamma}^0$ onto $\bar{\Delta}^0$. Clearly, $\chi_\Delta
\psi(v) = \bar{\phi} \chi_\Gamma(v)$ for any $v \in \Gamma^0$.

We now claim that  there exists a bijection $\varphi$ from
$\Gamma^1 \setminus \bigcup_{v \in \Gamma^0}T_{[v]}^1$ onto
$\Delta^1 \setminus \bigcup_{u \in \Delta^0}T_{[u]}^1$ such that for
any $e \in \Gamma^1 \setminus \bigcup_{v \in \Gamma^0}T_{[v]}^1$,
$\phi(e)$ contains $\varphi(e)$ if $e$ has exits and $s(\phi(e))
\sim s(\varphi(e))$ in $\Delta$ if $e$ is NE.
To see this,  note Lemma~\ref{NEspanning} and take a set of
NE spanning trees for $\Gamma$. Then the edges in $\Gamma^1
\setminus \bigcup_{v \in \Gamma^0}T_{[v]}^1$ can be divided into
three types: edges  $e$ where the out-degree of $s(e)$ is at
least 2 and $s(e)$ is not $\sim$ related to $r(e)$;  edges  $e_C$ in a
cycle $C$ which has exits only at $s(e_C)$; and edges $e_C$ in an NE
cycle $C$. Define $\varphi$ as the following: for the first two
types, $\varphi(e)$ is $\check{e}$ as in Lemma~\ref{levisolem}(d);
for the third type, $\varphi(e_C)$ is $e_{C'}$, where $C'$ is the
cycle corresponding to $C$ according to Lemma \ref{equalequiv}(e).
Thus, it follows from Lemmas~\ref{levisolem}(e), \ref{equivdis}(c) and
\ref{equalequiv}(b), (d), (f) that $\varphi$ preserves the types of
elements from $\Gamma^1 \setminus \bigcup_{v \in \Gamma^0}T_{[v]}^1$
to $\Delta^1 \setminus \bigcup_{u \in \Delta^0}T_{[u]}^1$ and is a
bijection. Moreover, we also observe from Lemma~\ref{levisolem}(d) that $\phi(e)$ contains
$\varphi(e)$ if $e$ has exits, and from Lemma~\ref{levisolem}(c) that  $s(\phi(e)) \sim
s(\varphi(e))$ in $\Delta$ if $e$ is NE. This proves the claim.

Take a bijection $\varphi$ from $\Gamma^1
\setminus \bigcup_{v \in \Gamma^0}T_{[v]}^1$ onto $\Delta^1
\setminus \bigcup_{u \in \Delta^0}T_{[u]}^1$ constructed as above.
For any $\bar{e} \in \bar{\Gamma}^1$, define $\bar{\phi}(\bar{e}) =
\chi_\Delta \varphi(e)$. By the
definition of $\chi_\Delta$, $\bar{\phi}$ restricts to a (well-defined) bijection from
$\bar{\Gamma}^1$ onto $\bar{\Delta}^1$.
For any $\bar{e} \in \bar{\Gamma}^1$, we know that $e \in \Gamma^1
\setminus \bigcup_{v \in \Gamma^0}T_{[v]}^1$. It follows from the proof of the claim above that $\phi(s(e)) \sim s(\varphi(e))$ and
$\phi(r(e)) \sim r(\varphi(e))$ in $\Delta$. Thus we observe from the definitions of
$\chi_\Gamma, \chi_\Delta$ and $\bar{\phi}$ that
$s(\bar{\phi}(\bar{e})) = [s(\varphi(e))] = [\phi(s(e))] = \bar{\phi}([s(e)]) = \bar{\phi}(s(\bar{e}))$. Similarly, $r(\bar{\phi}(\bar{e})) = \bar{\phi}(r(\bar{e}))$. Therefore, $\bar{\phi}$ is an graph isomorphism from $\bar{\Gamma}$ to $\bar{\Delta}$. \QED

\medskip

We have proved the direct part of the following theorem, which
classifies graphs with isomorphic Leavitt inverse semigroups.

\begin{Theorem}\label{leavinvclass}
Let $\Gamma$ and $\Delta$ be graphs. Then $LI(\Gamma) \cong
LI(\Delta)$ if and only if there is a bijection $\psi : \Gamma^0
\rightarrow \Delta^0$, sets of NE spanning trees $\{T_{[v]} : v \in
\Gamma^0\}$ and $\{T_{[u]} : u \in \Delta^0\}$, and a graph
isomorphism $\bar{\phi} : \bar{\Gamma} \rightarrow \bar{\Delta}$
such that, for all $v \in \Gamma^0$, $\bar{\phi}(\chi_{\Gamma}(v)) =
\chi_{\Delta}(\psi(v))$; that is, the following diagram is
commutative.
\begin{center}
\begin{tikzpicture}
\node at (0, 0) (v1) {$\bar{\Gamma}^0$};
\node at (0, 2) (v3) {$\Gamma^0$};
\node at (2, 0) (v2) {$\bar{\Delta}^0$};
\node at (2, 2) (v4) {$\Delta^0$};
\draw[->] (v1) to node[below] {$\bar{\phi}$} (v2);
\draw[->] (v3) to node[above] {$\psi$} (v4);
\draw[->] (v3) to node[left] {$\chi_\Gamma$} (v1);
\draw[->] (v4) to node[right] {$\chi_\Delta$} (v2);
\end{tikzpicture}
\end{center}
\end{Theorem}

We need some additional notation and lemmas to prove the converse part of Theorem \ref{leavinvclass}.

\begin{Lemma}
\label{inversethm}
Suppose that there is a bijection $\psi : \Gamma^0
\rightarrow \Delta^0$, sets of NE spanning trees $\{T_{[v]} : v \in
\Gamma^0\}$ and $\{T_{[u]} : u \in \Delta^0\}$, and a graph
isomorphism $\bar{\phi} : \bar{\Gamma} \rightarrow \bar{\Delta}$
such that, for all $v \in \Gamma^0$, $\bar{\phi}(\chi_{\Gamma}(v)) =
\chi_{\Delta}(\psi(v))$. Then

$(a)$ for all $u \in \Delta^0$, $\bar{\phi}^{-1}(\chi_{\Delta}(u)) = \chi_{\Gamma}(\psi^{-1}(u))$; and

$(b)$ $\psi$ restricts to a bijection from $[v]$ to $[\psi(v)]$ for all $v \in \Gamma^0$.

\end{Lemma}

\noindent {\bf Proof.} The proof of part (a) follows in a  routine fashion from the fact that $\bar{\phi}$ and $\psi$ are bijections. If $v_1 \sim v_2$ in $\Gamma^0$, then $\chi_{\Gamma}(v_1) = \chi_{\Gamma}(v_2)$, and so $\chi_{\Delta}(\psi(v_1)) = \bar{\phi}(\chi_{\Gamma}(v_1)) = \bar{\phi}(\chi_{\Gamma}(v_2))  = \chi_{\Delta}(\psi(v_2))$, that is $\psi(v_1) \sim \psi(v_2)$ in $\Delta^0$. Similarly, from part (a) it follows that if $\psi(v_1) \sim \psi(v_2)$ then $v_1 = \psi^{-1}(\psi(v_1)) \sim \psi^{-1}(\psi(v_2)) = v_2$. Part (b) easily follows from this and the fact that $\psi$ is a bijection. \QED

\medskip

Fixing an NE spanning tree of $\Gamma_{[v]}$, for any $v_1, v_2 \in [v]$, we observe that
there exist directed NE
paths $p, q$ in $T_{[v]}$ such that $r(p) = r(q)$, $v_1 = s(p)$ and
$v_2 = s(q)$. These paths $p, q$ are not necessarily unique. However,
there is a unique shortest such directed NE path $p$, and a unique
shortest such path $q$.  We denote this choice of
$p q^*$ by $p[v_1, v_2]$. Clearly $p[v_1, v_2]$ has no non-trivial circuits.

\begin{Lemma}
\label{short}

$(a)$ Let $v_1,v_2,v_3,v$ be vertices of a graph $\Gamma$ such that $v_1, v_2, v_3 \in [v]$. Then as elements of $LI(\Gamma)$, we have $p[v_1, v_2] p[v_2, v_3] = p[v_1, v_3]$, $p[v_1, v_2] p[v_2, v_1] = v_1$ and $p[v_1, v_2]^* = p[v_2, v_1]$.

$(b)$ If $v_1 \sim v_2 \in \Gamma^0$ and $v_2 = s(e)$ for some edge $e$ that is not in $T_{[v_2]}$, then $p[v_1, v_2]$ is a directed NE path from $v_1$ to $v_2$.

\end{Lemma}

\noindent {\bf Proof.} The proof of part (a) follows easily from the definitions. For part (b) there are two cases. Suppose that $p[v_1, v_2] = pq^*$ for NE paths $p,q$ in $T_{[v_2]}$. If $v_2$ has out-degree greater than $1$ then clearly $q$ must be $v_2$ (the empty path  at $v_2$), so $p[v_1, v_2] = p$ is a directed path from $v_1$ to $v_2$. If $v_2$ has out-degree $1$ but $e$ is not in the spanning tree $T_{[v_2]}$, then again $q$ must be the empty path at $v_2$. This is because if $q$ is not empty, then the first edge of $q$ must be $e$, a contradiction since $q$ is in the spanning tree. Hence  again $p[v_1, v_2] = p$ is a directed NE path from $v_1$ to $v_2$. \QED

\medskip

Now we define a mapping $\hat{\chi}_\Gamma$ from $LI(\bar{\Gamma})$
to $LI(\Gamma)$. For any $v \in \Gamma^0$, fix a vertex $v_0$ in the
$\sim$-class $[v]$. Let $\hat{\chi}_\Gamma$ map the vertex $[v]$ in
$\bar{\Gamma}$ to $v_0$ in $\Gamma$. For any directed path $\bar{p}
= \bar{e}_1 \ldots \bar{e}_m$ in $\bar{\Gamma}$, define
\begin{equation*} \label{barp}
\hat{\chi}_{\Gamma}(\bar{p}) = p[\hat{\chi}_\Gamma(s(\bar{e}_1)), s(e_1)] \, e_1 \, p[r(e_1), s(e_2)] \, e_2 \ldots e_{m-1} \, p[r(e_{m-1}), s(e_m)] \, e_m \, p[r(e_m), \hat{\chi}_\Gamma(r(\bar{e}_m))].
\end{equation*}

Notice that since each edge $e_i$ is not in a spanning tree, we have
$p[\hat{\chi}_\Gamma(s(\bar{e}_1)), s(e_1)] = p_1$ for some NE path
$p_1$ and  $p[r(e_i), s(e_{i+1})] = p_{i+1}$ for some NE
path $p_{i+1}$ by Lemma \ref{short}(b). So, as an element of
$LI(\Gamma)$, $\hat{\chi}_{\Gamma}(\bar{p})$ is a non-zero element
of the form $p_1e_1p_2e_2\ldots p_me_mp_{m+1} p_{m+2}^*$, $r(p_{m+1}) = r(p_{m+2})$ and the $p_i$ are NE paths in the spanning trees. In particular, $\hat{\chi}_{\Gamma}(\bar{e}) =
p_1 e p_2 p_3^*$ where $p_1, p_2, p_3$ are in the spanning trees.
For a nonzero element $\bar{p} \bar{q}^*$ in $LI(\bar{\Gamma})$,
define $\hat{\chi}_{\Gamma}(\bar{p} \bar{q}^*) =
\hat{\chi}_{\Gamma}(\bar{p})(\hat{\chi}_{\Gamma}(\bar{q}))^*$. If $\bar{q} = \bar{f}_1...\bar{f}_n$ then
$\hat{\chi}_{\Gamma}(\bar{p} \bar{q}^*) = p_1e_1\ldots p_me_mp'q'^*f_n^*p_n'^*\ldots f_1^*p_1'^*$ for some NE paths $p, p_i,p_j',p'$ in the spanning trees with $r(p') = r(q')$.

\medskip

We call a $0$-morphism $f : S \rightarrow T$ from an inverse semigroup $S$ onto an inverse subsemigroup $T$ of $S$ a {\em $0$-retraction} (and we call $T$ a {\em $0$-retract} of $S$)  if the restriction of $f$ to $T$ is the identity map on $T$.

\begin{Lemma} \label{chihat}
The mapping $\hat{\chi}_{\Gamma}$ is a monomorphism from
$LI(\bar{\Gamma})$ to $LI(\Gamma)$ such that $\tilde{\chi}_\Gamma
\hat{\chi}_\Gamma$ is the identity mapping of $LI(\bar{\Gamma})$. Hence $LI(\bar{\Gamma})$ is isomomorphic to  a $0$-retract of $LI(\Gamma)$.
\end{Lemma}

\noindent {\bf Proof.}  Since the map $\tilde{\chi}_{\Gamma}$
induces a bijection from $\{e \in \Gamma^1 \setminus \bigcup_{v \in
\Gamma^0}T_{[v]}^1\}$ onto $\bar{\Gamma}^1$ the map from $\bar{\Gamma}^1$ into $\Gamma^1$ defined by $\bar{e}
\mapsto e$ is an injection. But also the path $p[v_1,v_2]$ is uniquely
determined by the vertices $v_1,v_2 \in [v]$ and the choice of
spanning tree $T_{[v]}$, so it follows that $\hat{\chi}_{\Gamma}$ is
an injection from $LI(\bar{\Gamma})$ to $LI(\Gamma)$. A routine
argument, using the characterization of $\hat{\chi}_{\Gamma}(\bar{p} \bar{q}^*)$ in the paragraph above, shows that $\hat{\chi}_{\Gamma}$ is a homomorphism, so it is a monomorphism. The fact that
$\tilde{\chi}_\Gamma \hat{\chi}_\Gamma$ is the identity mapping of
$LI(\bar{\Gamma})$ follows immediately from the definitions. Thus $\hat{\chi}\tilde{\chi}$ is a surjective $0$-morphism from $LI(\Gamma)$ onto the inverse subsemigroup $\hat{\chi}(LI(\bar{\Gamma}))$ of $LI(\Gamma)$ and $\hat{\chi}\tilde{\chi}(\hat{\chi}(x)) = \hat{\chi}(\tilde{\chi}\hat{\chi}(x)) = \hat{\chi}(x)$ for all $x \in LI(\bar{\Gamma})$; that is, the restriction of $\hat{\chi}\tilde{\chi}$ to $\hat{\chi}(LI(\bar{\Gamma}))$ is the identity map. Hence $\hat{\chi}_{\Gamma}(LI(\bar{\Gamma}))$ is a $0$-retract of $LI(\Gamma)$. The result follows since $\hat{\chi}$ is an isomomorphism from $LI(\hat{\Gamma})$ onto $\hat{\chi}_{\Gamma}(LI(\bar{\Gamma}))$. \QED

\medskip

 Let $\psi : \Gamma^0 \rightarrow
\Delta^0$ be a bijection, $\{T_{[v]} : v \in \Gamma^0\}$ and
$\{T_{[u]} : u \in \Delta^0\}$ be sets of NE spanning trees and
$\bar{\phi} : \bar{\Gamma} \rightarrow \bar{\Delta}$ be a graph
isomorphism such that, for all $v \in \Gamma^0$,
$\bar{\phi}(\chi_{\Gamma}(v)) = \chi_{\Delta}(\psi(v))$. Then
$\bar{\phi}$ naturally induces an isomorphism $\hat{\phi}$ from
$LI(\bar{\Gamma})$ onto $LI(\bar{\Delta})$ which maps directed paths
to directed paths. From Lemmas~\ref{chitilde} and
\ref{chihat}, we get a 0-restricted morphism $\tilde{\phi} =
\hat{\chi}_{\Delta} \hat{\phi} \tilde{\chi}_\Gamma$ from
$LI(\Gamma)$ into $LI(\Delta)$.

Note that
$\tilde{\phi}(v) = \hat{\chi}_{\Delta} \hat{\phi} \tilde{\chi}_{\Gamma} (v) = \hat{\chi}_{\Delta} \bar{\phi}\chi_{\Gamma}(v) = \hat{\chi}_{\Delta}\chi_{\Delta}\psi(v) \sim \psi(v)$, and also that $s(\tilde{\phi}(p q^*)) = s(\tilde{\phi}(p))$ since $\tilde{\phi}$ is a 0-morphism and $s(\tilde{\phi}(p)) = \tilde{\phi}(s(p))$ since $\hat{\phi} \tilde{\chi}$ is a 0-morphism. Similarly, $r(\tilde{\phi}(p q^*)) = s(\tilde{\phi}(q)) = \tilde{\phi}(s(q))$. In view of these facts we may define, for any nonzero element $pq^* \in LI(\Gamma)$,
$$
\phi(p q^*) = p[\psi(s(p)), \tilde{\phi}(s(p))] \, \tilde{\phi}(p q^*) \, p[\tilde{\phi}(s(q)), \psi(s(q))]
$$
and $\phi(0) = 0$. Then $\phi$ is a well-defined function from $LI(\Gamma)$ to $LI(\Delta)$ and
 $\phi(p q^*)$ is nonzero for a
nonzero element $p q^* \in LI(\Gamma)$. In particular, for any
directed path $p$ in $\Gamma$ (which we may think of as $p r(p)^*$, where $r(p)$ is the empty path at the vertex  $r(p)$), we have
\begin{equation} \label{phiofp}
\begin{split}
\phi(p) & = p[\psi(s(p)), \tilde{\phi}(s(p))] \, \tilde{\phi}(p) \, p[\tilde{\phi}(s(r(p))), \psi(s(r(p)))] \\
& = p[\psi(s(p)), \tilde{\phi}(s(p))] \, \tilde{\phi}(p) \, p[\tilde{\phi}(r(p)), \psi(r(p))]
\end{split}
\end{equation}
and for any vertex $v$ in $\Gamma$, we have $\phi(v) = \psi(v)$.

\medskip

Lemma \ref{phi} below provides a proof of the converse part of
Theorem~\ref{leavinvclass}.

\begin{Lemma} \label{phi}
The map $\phi$ is an  isomorphism from $LI(\Gamma)$ onto $LI(\Delta)$.
\end{Lemma}

\noindent {\bf Proof}.
Let $p_1 q_1^*$ and $ p_2 q_2^*$ be arbitrary non-zero elements in
$LI(\Gamma)$. Then we obtain from the definition of $\phi$ that
\begin{equation} \label{phiofpqpq}
\begin{split}
\phi(p_1 q_1^*) \phi(p_2 q_2^*) = &
\, p[\psi(s(p_1)), \tilde{\phi}(s(p_1))] \, \tilde{\phi}(p_1 q_1^*) \, p[\tilde{\phi}(s(q_1)), \psi(s(q_1))] \, \bullet \\
& \, p[\psi(s(p_2)), \tilde{\phi}(s(p_2))] \, \tilde{\phi}(p_2 q_2^*) \, p[\tilde{\phi}(s(q_2)), \psi(s(q_2))].
\end{split}
\end{equation}
If $(p_1 q_1^*)(p_2 q_2^*) \neq 0$, then $s(q_1) = s(p_2)$. By Lemma \ref{short}(a) and the fact that $\tilde{\phi}$ is a 0-morphism, we obtain
\begin{equation} \label{phiofpqpq2}
\begin{split}
\phi(p_1 q_1^*) \phi(p_2 q_2^*) = & \, p[\psi(s(p_1)), \tilde{\phi}(s(p_1))] \, \tilde{\phi}((p_1 q_1^*) (p_2 q_2^*)) \, p[\tilde{\phi}(s(q_2)), \psi(s(q_2))] \\
= & \, \phi((p_1 q_1^*) (p_2 q_2^*)).
\end{split}
\end{equation}
Suppose that $(p_1q_1^*)(p_2q_2^*) = 0$.
If $s(q_1) \neq s(p_2)$,
then since $\psi$ is injective, we see that $\psi(s(q_1)) \neq \psi(s(p_2))$ which means by (\ref{phiofpqpq})
that $\phi(p_1 q_1^*) \phi(p_2 q_2^*) = 0$. If $s(q_1) =
s(p_2)$, then we know from (\ref{phiofpqpq}) that
\begin{align} \label{phiofpqpq1}
\phi(p_1 q_1^*) \phi(p_2 q_2^*) = p[\psi(s(p_1)), \tilde{\phi}(s(p_1))] \, \tilde{\phi}(p_1 q_1^*) \tilde{\phi}(p_2 q_2^*) \, p[\tilde{\phi}(s(q_2)), \psi(s(q_2))].
\end{align}
Since $(p_1 q_1^*)(p_2 q_2^*) = 0$, then neither of $q_1, p_2$ is a prefix
of the other which means that in both $q_1$ and $p_2$, there exists
a vertex whose out-degree is at least 2. Since our chosen
spanning trees in both $\Gamma$ and $\Delta$ are NE and
$\bar{\phi}$ is an isomorphism, we see that in both
$\tilde{\phi}(q_1)$ and $\tilde{\phi}(p_2)$, there exists a vertex
whose out-degree is at least 2. So we obtain from (\ref{phiofpqpq1})
that $\phi(p_1 q_1^*) \phi(p_2 q_2^*) = 0$. Thus
  $\phi$ is a semigroup
homomorphism.

To see that $\phi$ is surjective, we only need prove that every edge
in $\Delta$ has a  preimage under $\phi$ since $LI(\Delta)$ is
generated by $\Delta^0$ and $\Delta^1$ and we already established that $\phi(v) = \psi(v)$ for each vertex $v$ in $\Gamma$.  Given an edge $d$ in
$\Delta$, we take
$$
x = p[\psi^{-1}(s(d)), \hat{\chi}_\Gamma \bar{\phi}^{-1} \chi_\Delta(s(d))] \, \hat{\chi}_\Gamma \bar{\phi}^{-1} \chi_\Delta(d) \, p[\hat{\chi}_\Gamma \bar{\phi}^{-1} \chi_\Delta(r(d)), \psi^{-1}(r(d))]
$$
which is nonzero by a similar discussion as in the first paragraph of the proof.
Note that each $p[v_1, v_2]$ involves only NE paths in spanning trees and that $\bar{\phi}^{-1}(\chi_{\Delta}(u)) = \chi_{\Gamma}(\psi^{-1}(u))$ for any $u \in \Delta^0$ by Lemma \ref{inversethm}(a).
Then we observe from (\ref{phiofpqpq2}), Lemma~\ref{chihat} and the related definitions that
\begin{align}
\phi(x) = & \, p[\psi(\psi^{-1}(s(d))), \tilde{\phi}(\psi^{-1}(s(d)))] \, \tilde{\phi}(p[\psi^{-1}(s(d)), \hat{\chi}_\Gamma \bar{\phi}^{-1} \chi_\Delta(s(d))]) \tilde{\phi}(\hat{\chi}_\Gamma \bar{\phi}^{-1} \chi_\Delta(d)) \, \bullet \notag \\
& \, \tilde{\phi}(p[\hat{\chi}_\Gamma \bar{\phi}^{-1} \chi_\Delta(r(d)), \psi^{-1}(r(d))]) \, p[\tilde{\phi}(\psi^{-1}(r(d))), \psi(\psi^{-1}(r(d)))] \notag \\
= & \, p[s(d), \tilde{\phi}(\psi^{-1}(s(d)))] \, \tilde{\phi}(\hat{\chi}_\Gamma \bar{\phi}^{-1} \chi_\Delta(d)) \, p[\tilde{\phi}(\psi^{-1}(r(d))), r(d)] \notag \\
= & \, p[s(d), \hat{\chi}_\Delta \chi_\Delta(s(d))] \, \hat{\chi}_\Delta \chi_\Delta(d) \, p[\hat{\chi}_\Delta \chi_\Delta(r(d)), r(d)]. \label{phiofx}
\end{align}
If $d$ is in a spanning tree, then we see from (\ref{phiofx}) that $\hat{\chi}_\Delta \chi_\Delta(s(d)) = \hat{\chi}_\Delta \chi_\Delta(d) = \hat{\chi}_\Delta \chi_\Delta(r(d))$ so that $\phi(x) = p[s(d), r(d)] = d$. If $d$ is not in any spanning tree, then again we observe from (\ref{phiofx}) that
$$
\phi(x) = p[s(d), \hat{\chi}_\Delta \chi_\Delta(s(d))] p[\hat{\chi}_\Delta \chi_\Delta(s(d)),s(d)] \, d \, p[r(d), \hat{\chi}_\Delta \chi_\Delta(r(d))] p[\hat{\chi}_\Delta \chi_\Delta(r(d)), r(d)] = d.
$$
Therefore, $d$ has a preimage under $\phi$.

We have seen that nonzero elements map to nonzero ones by $\phi$. Note the canonical forms of
Leavitt inverse semigroups in Section~5. For any nonzero (reduced) elements
$p_1 q_1^*, p_2 q_2^*$ in $LI(\Gamma)$, if $p_1 q_1^* \neq p_2 q_2^*$, then we may suppose
that $p_1 \neq p_2$. If $s(p_1) \neq s(p_2)$, then $s(\phi(p_1)) = \psi(s(p_1)) \neq \psi(s(p_2)) = s(\phi(p_2))$ which means that $\phi(p_1 q_1^*) \neq \phi(p_2 q_2^*)$. If $s(p_1) = s(p_2)$, then we assume that $p_1 = e_1 \ldots e_m$, $p_2 = e_1' \ldots e_n'$, $e_1 = e_1', \cdots, e_{i-1} = e_{i-1}'$ but $e_i \neq e_i'$ for some $i \in \{1, \cdots, m\}$. Thus, $s(e_i)$ has out-degree at least 2. Since the chosen spanning trees for $\Gamma$ are NE, we obtain that $\tilde{\phi}(e_i) \neq \tilde{\phi}(e_i')$ which leads to
\begin{align*}
\phi(p_1) & = p[\psi(s(e_1)), \tilde{\phi}(s(e_1))] \tilde{\phi}(e_1) \ldots \tilde{\phi}(e_m) p[\psi(r(e_m)), \tilde{\phi}(r(e_m))] \\
& \neq p[\psi(s(e_1')), \tilde{\phi}(s(e_1'))] \tilde{\phi}(e_1') \ldots \tilde{\phi}(e_n') p[\psi(r(e_n')), \tilde{\phi}(r(e_n'))] = \phi(p_2)
\end{align*}
because $\tilde{\phi}(e_i)$ and $\tilde{\phi}(e_i')$ contain
distinct edges which have the same source by
Lemma~\ref{levisolem}(d). So in this case, we also have $\phi(p_1
q_1^*) \neq \phi(p_2 q_2^*)$. We proved that $\phi$ is injective, so
 it is a semigroup isomorphism. This completes the proof of the lemma, and hence of Theorem \ref{leavinvclass}. \QED

\medskip

\noindent {\bf Example.} We  illustrate the use of Theorem \ref{leavinvclass} to construct an isomorphism $\phi$ from the Leavitt inverse semigroup $LI(\Gamma)$ onto the Leavitt inverse semigroup $LI(\Delta)$ for the two graphs $\Gamma$ and $\Delta$ in Diagram~7.2. By Theorem \ref{algiso}, the Leavitt path algebras $L_F(\Gamma)$ and $L_F(\Delta)$ of these graphs are also isomorphic.

$\Gamma^0$ has two $\sim$-classes $[v_1] = \{v_1, v_2, v_3\}$ and $[v_4] = \{v_4, v_5, v_6\}$; $\Delta^0$ also has two $\sim$-classes $[u_1] = \{u_1, u_2, u_3\}$ and $[u_4] = \{u_4, u_5, u_6\}$. As we discussed, $e_{C_1}$ and $e_{C_2}$ respectively determine a spanning tree of $\Gamma_{[v_1]}$ and $\Gamma_{[v_4]}$; $e_{C'_1}$ and $e_{C'_2}$ respectively determine a spanning tree of $\Delta_{[u_1]}$ and $\Delta_{[u_4]}$. Since $\Gamma$ and $\Delta$ have isomorphic contracted graphs, for convenience, we denote them as $\bar{\Gamma} = \bar{\Delta}$, where $w_1 = [v_1] = [u_1]$, $w_2 = [v_4] = [u_4]$, $f_1 = \bar{e}_{C_1} = \bar{e}_{C'_1}$, $f_2 = \bar{e}_5 = \bar{d}_5$ and $f_3 = \bar{e}_{C_2} = \bar{e}_{C'_2}$.

\begin{center}
\begin{tikzpicture}
\node at (-3, 0) {$\Gamma:$}; \node at (0, 0) (v3) {\color{red}{$\bullet$}}; \node at (-1.731, 1) (v2)
{\color{red}{$\bullet$}}; \node at (-1.732, -1) (v1)
{\color{red}$\bullet$}; \node at (2, 0) (v4) {$\color{blue}\bullet$};
\node at (3.732, 1) (v5) {\color{blue}$\bullet$};
\node at (3.732, -1) (v6) {\color{blue}$\bullet$};
\draw[->, red] (v3) node[above] {$v_3$} to node[below, sloped] {$e_{C_1}$} (v1) node[below] {$v_1$}; \draw[->, red] (v1) to node[left] {$e_1$} (v2) node[above] {$v_2$}; \draw[->, red] (v2)
to node[above,sloped] {$e_2$} (v3); \draw[->] (v3) to
node[below] {$e_5$} (v4); \draw[->, blue] (v4) node[above] {$v_4$} to node[above, sloped] {$e_3$} (v5) node[above] {$v_5$}; \draw[->, blue] (v5) to node[right] {$e_{C_2}$} (v6) node[below] {$v_6$}; \draw[->, blue] (v6) to node[below, sloped] {$e_4$} (v4);
\end{tikzpicture}
\bigskip
\begin{tikzpicture}
\node at (-5, 0) {$\Delta:$}; \node at (-1.732, 0) (u3) {\color{red}{$\bullet$}}; \node at (-3.732, 0) (u2)
{\color{red}{$\bullet$}}; \node at (0, 0) (u1)
{\color{red}$\bullet$}; \node at (2, 0) (u5) {$\color{blue}\bullet$};
\node at (2, 2) (u6) {\color{blue}$\bullet$};
\node at (4, 0) (u4) {\color{blue}$\bullet$};
\draw[->, red] (u2) node[below] {$u_2$} to node[above] {$d_1$} (u3) node[below] {$u_3$}; \draw[->, red] (u3) to [bend left = 30] node[above] {$d_2$} (u1); \draw[red] (u1) node[below] {$u_1$}; \draw[->, red] (u1)
to [bend left = 30] node[below] {$e_{C'_1}$} (u3); \draw[->] (u1) to
node[above] {$d_5$} (u5); \draw[->, blue] (u5) node[below] {$u_5$} to node[above] {$d_4$} (u4) node[below] {$u_4$}; \draw[->, blue] (u6) node[left] {$u_6$} to node[right] {$d_3$} (u5); \draw[->, blue] (u4) edge [loop right] node[right] {$e_{C'_2}$} (u4);
\end{tikzpicture}
\bigskip
\begin{tikzpicture}
\node at (-2.5, 0) {$\bar{\Gamma} = \bar{\Delta}:$}; \node at (0, 0) (w1) {\color{red}{$\bullet$}}; \node at (2, 0) (w2) {\color{blue}{$\bullet$}};
\draw[->, red] (w1) edge [loop left] node[left] {$f_1$} (w1) node[below] {$w_1$}; \draw[->] (w1) to node[above] {$f_2$} (w2);
\draw[->, blue] (w2) edge [loop right] node[right] {$f_3$} (w2) node[below] {$w_2$};
\end{tikzpicture}

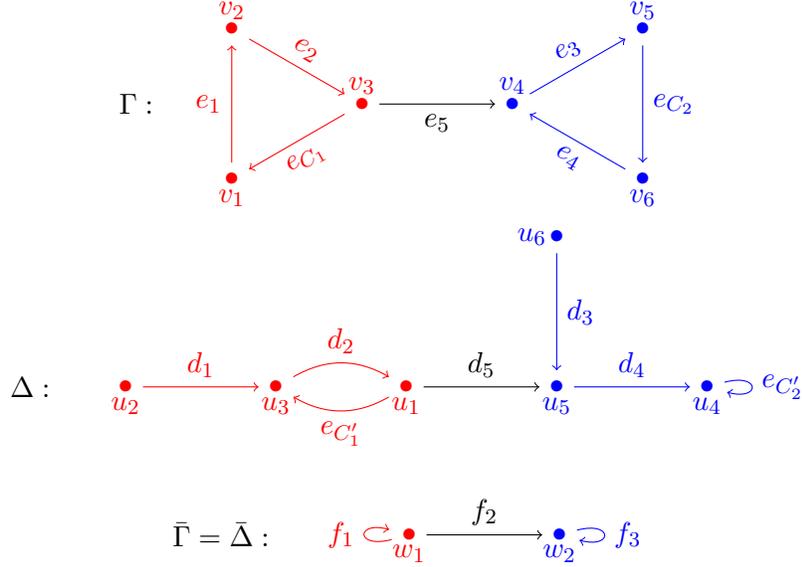
Diagram~7.2 \,\, Two graphs with isomorphic Leavitt inverse semigroups\\
and their contracted graph(s)
\end{center}

 Let $\bar{\phi}$ be the identity mapping of $\bar{\Gamma}$ and $\psi(v_i) = u_i$ for $i = 1, 2, \cdots, 6$. To construct a semigroup isomorphism $\phi$ from $LI(\Gamma)$ to $LI(\Delta)$ according to Theorem~\ref{leavinvclass}, we only need to list the image of generators, in fact only edges, under $\phi$. First, by (\ref{phiofp}) and Lemma~\ref{short}(a), for any $e \in \Gamma^1$, if $e$ is in a spanning tree, then we see from $\chi_\Gamma(s(e)) = \chi_\Gamma(e) = \chi_\Gamma(r(e))$ that $\phi(e) = p[\psi(s(e)), \hat{\chi}_\Delta \chi_\Gamma(s(e))] \, p[\hat{\chi}_\Delta \chi_\Gamma(r(e)), \psi(r(e)] = p[\psi(s(e)), \psi(r(e))]$, and if $e$ is not in a spanning tree, then $\phi(e) = p[\psi(s(e)), s(\hat{\chi}_\Delta(\bar{e}))] \, \hat{\chi}_\Delta(\bar{e}) \, p[r(\hat{\chi}_\Delta(\bar{e})), \psi(r(e)]$. Thus,
\begin{align*}
& \phi(e_1) = p[u_1, u_2] = d_2^* d_1^*; \\
& \phi(e_2) = p[u_2, u_3] = d_1; \\
& \phi(e_3) = p[u_4, u_5] = d_4^*; \\
& \phi(e_4) = p[u_6, u_4] = d_3 d_4; \\
& \phi(e_5) = p[u_3, u_1] \, d_5 \, p[u_5, u_4] = d_2 d_5 d_4; \\
& \phi(e_{C_1}) = p[u_3, u_1] \, e_{C'_1} \, p[u_3, u_1] = d_2 e_{C'_1} d_2; \\
& \phi(e_{C_2}) = p[u_5, u_4] \, e_{C'_2} \, p[u_4, u_6] = d_4 e_{C'_2} d_4^* d_3^*.
\end{align*}

We close the paper with several results that follow from Theorem \ref{leavinvclass}.

\begin{Cor}
Let $\bar{\Gamma}$ and $\bar{\Delta}$ be contracted graphs. Then
$LI(\bar{\Gamma})$ is isomorphic to $LI(\bar{\Delta})$ if and only
if $\bar{\Gamma}$ is isomorphic to $\bar{\Delta}$.
\end{Cor}

\noindent {\bf Proof.} The direct part follows from
Lemma~\ref{deg1contracted} and the direct part of
Theorem~\ref{leavinvclass}. The converse part is trivial. \QED

\begin{Cor}
\label{leavittbouquet}

If $\Gamma$ and $\Delta$ are graphs such that $\bar{\Gamma} \cong \bar{\Delta} \cong B_X$ $($the bouquet of $|X|$ circles$)$, then $LI(\Gamma) \cong LI(\Delta)$ if and only if $|\Gamma^0| = |\Delta^0|$.

\end{Cor}

\noindent {\bf Proof.}  This is immediate from Theorem~\ref{leavinvclass} since there is only one equivalence class of vertices of $\Gamma$ (or $\Delta$) under the equivalence relation $\sim$. \QED

\medskip

We may view Corollary \ref{leavittbouquet} as a generalization of Theorem \ref{algiso1} since a connected graph that immerses into a circle has only one $\sim$-class of vertices.

\medskip

\noindent {\bf Remark} We remark that the hypotheses of  Corollary \ref{leavittbouquet} do not classify graphs whose contracted graphs are isomorphic to a bouquet of circles and which have  isomorphic Leavitt path algebras. It clearly follows from Corollary \ref{leavittbouquet} and Theorem \ref{algiso} that if $\Gamma$ and $\Delta$ are two graphs with $\bar{\Gamma} \cong \bar{\Delta} \cong B_X$ and $|\Gamma^0| = |\Delta^0|$ then $L_F(\Gamma) \cong L_F(\Delta)$. But the conditions $L_F(\Gamma) \cong L_F(\Delta)$ and $\bar{\Gamma} \cong \bar{\Delta} \cong B_X$ do not necessarily imply that $|\Gamma^0| = |\Delta^0|$. In fact  the following result follows easily from Theorem~\ref{leavinvclass} and some results in the paper by  Abrams, \'Anh and Pardo \cite{AAP}.

\begin{Cor}
\label{levalgb2}

Let $\Gamma$  be a finite graph with $\bar{\Gamma} \cong B_X$ where $|X| = n \geq 2$. Then

$(a)$ $L_F(\Gamma) \cong M_{|\Gamma^0|}(L_F(B_X))$ $($where $L_F(B_X)$ is the Leavitt algebra $L_F(1,n))$;

$(b)$ If $n = 2$ then $L_F(\Gamma) \cong L_F(B_X) \cong L_F(1,2)$;

$(c)$ If $\Gamma$ and $\Delta$ are two graphs with $\bar{\Gamma} \cong \bar{\Delta} \cong B_X$ where $|X| = n > 2$, then $L_F(\Gamma) \cong L_F(\Delta)$ iff $g.c.d(|\Gamma^0|,n-1) = g.c.d(|\Delta^0|,n-1)$.

\end{Cor}

\noindent {\bf Proof.}
By Theorem~\ref{leavinvclass}, $LI(\Gamma) \cong LI(B^k(B_X))$ where $B^k(B_X))$ denotes the graph obtained from $B_X$ by attaching a directed NE path of length $k = |\Gamma^0|-1$ ending in  the (unique)  vertex  in the graph $B_X$. Hence $L_F(\Gamma) \cong L_F(B^k(B_X))$ by Theorem \ref{algiso}.   By Lemma 5.1 of \cite{AALP}, $L_F(B^k(B_X)) \cong M_k(L_F(1,n))$, the algebra of $k \times k$ matrices over the Leavitt algebra $L_F(1,n)$. The results of parts (b) and (c) follow immediately from Theorem 12 of \cite{A} (Theorems 4.14 and 5.12 of \cite{AAP}). \QED


\medskip

The ideas employed in the proof of Corollary \ref{levalgb2} may be extended somewhat to obtain a result relating the structure of the Leavitt path algebra $L_F(\Gamma)$ to the structure of the algebra $L_F(\bar{\Gamma})$ of the contracted graph $\bar{\Gamma}$ for any  graph $\Gamma$ with finite $\sim$-equivalence classes.

We define a function $f : L_F(\Delta) \rightarrow L_F(\Gamma)$ between Leavitt path algebras to be a {\em $0$-morphism} if $f$ is a  linear transformation between the underlying vector spaces that restricts to a $0$-morphism $LI(\Delta) \rightarrow LI(\Gamma)$ between the corresponding Leavitt inverse semigroups. We call  a function $f : L_F(\Delta) \rightarrow L_F(\Gamma)$  a {\em $0$-retraction} if $L_F(\Gamma)$ is a subalgebra  of $L_F(\Delta)$ and $f$ is a  $0$-morphism from $L_F(\Delta)$ onto $L_F(\Gamma)$ that restricts to the identity function on $L_F(\Gamma)$: equivalently, we say that $L_F(\Gamma)$ is a $0$-retract of $L_F(\Delta)$. We extend the notation slightly by saying that an $F$-algebra $A_1$ is a $0$-retract of an $F$-algebra $A_2$ if there are graphs $\Gamma$ and $\Delta$ such that $A_1 \cong L_F(\Gamma), A_2 \cong L_F(\Delta)$ and $L_F(\Gamma)$ is a $0$-retract of $L_F(\Delta)$. Recall from \cite{AAM} that a subgraph $\Gamma$ of a graph $\Delta$ is called a {\em complete} subgraph of $\Delta$ if $e \in \Gamma^1$ for every edge $e \in \Delta^1$ such that $s(e) \in \Gamma^0$.

\begin{Theorem}
\label{algstructure}

Let $\Gamma$ be a  graph with finite equivalence classes $[v]$ for each $v \in \Gamma^0$ and let $n$ be the maximum order of any equivalence class $[v]$ for $v \in \Gamma^0$. Then $L_F(\Gamma)$ is a $0$-retract of $M_n(L_F(\bar{\Gamma}))$.

\end{Theorem}

\noindent {\bf Proof.} By Theorem \ref{leavinvclass}, $L_F(\Gamma) \cong L_F(\Gamma_1)$ where $\Gamma_1$ is obtained from $\bar{\Gamma}$ by attaching an NE path $p_{[v]}$ of length $|[v]| -1$ ending at the vertex $[v]$ of $\bar{\Gamma}$ for each $[v] \in \bar{\Gamma}^0$. Let $\Delta$ be the graph obtained from $\bar{\Gamma}$ by attaching an NE path $q_{[v]}$ of length $n$ ending at  $[v] \in \bar{\Gamma}^0$ for each vertex $[v] \in \bar{\Gamma}^0$. By Proposition 9.3 of a paper by Abrams and Tomforde \cite{AT}, $L_F(\Delta) \cong M_n(L_F(\bar{\Gamma}))$. By identifying $p_{[v]}$ with the suffix of $q_{[v]}$ of length $|[v]|-1$, we may view $\Gamma_1$ as a subgraph of $\Delta$, in fact a {\em complete} subgraph in the sense of \cite{AAM}. Hence by Lemma 1.6.6 of \cite{AAM}, $L_F(\Gamma_1)$ is a subalgebra of $L_F(\Delta)$. We show by induction on the number of edges in $\Delta^1 \setminus \Gamma_1^1$ that in fact there is a $0$-retraction of $L_F(\Delta)$ onto $L_F(\Gamma_1)$.

Suppose that $\Delta_1$ is obtained by attaching one NE edge $e$ to $\Gamma_1$ with $r(e) = s(p_{[v]})$ and $s(e) \notin \Gamma_1$ and such that $\Delta_1$ is a subgraph of $\Delta$. We claim that we may construct a well-defined map  $f : L_F(\Delta_1) \rightarrow L_F(\Gamma_1)$ by contracting the edge $e$ and the vertex $s(e)$ to the vertex $r(e) = s(p_{[v]})$ in $\Gamma_1^0$. More precisely, we proceed as follows. Choose the special edges of $\Delta_1$ (in the sense of \cite{AAJZ}) in such a way that all special edges of $\Delta_1$ other than $e$ are special edges of $\Gamma_1$. Note that since $s(e)$ is a source in $\Delta_1$ of out-degree $1$, a directed path $p$ in $\Delta_1$ contains the edge $e$ if and only if $p = et$ for some (possibly empty) directed path $t$ in $\Gamma_1$ with $s(t) = r(e)$.  From this and the choice of special edges in $\Delta_1$ it is easy to see that the corresponding natural basis of the algebra $L_F(\Delta_1)$ consists of the elements in the natural basis of $L_F(\Gamma_1)$ together with the elements $s(e),e,e^*$ and the non-zero elements $epq^*, pq^*e^*$ and $epq^*e^*$ where $pq^*$ is a non-empty natural basis element for the algebra $L_F(\Gamma_1)$. Define $e' = s(e)' = r(e)$ and for each other directed path $p$ in $\Delta_1$, define  $p'$  to be the directed path in $\Gamma_1$ obtained by deleting the first edge of $p$ if this edge is $e$ and $p'=p$ if the first edge of $p$ is not $e$. Then define $f$ on natural basis elements of $L_F(\Delta_1)$ by $f(pq^*) = p'q'^*$. This extends by linearity to a linear transformation (that we again denote by $f$) from $L_F(\Delta_1)$ to $L_F(\Gamma_1)$, and this linear transformation is surjective since every natural basis element of $L_F(\Gamma_1)$ arises as the image of a natural basis element of $L_F(\Delta_1)$. Furthermore, $f$ fixes every natural basis element of $L_F(\Gamma_1)$ so the restriction of $f$ to $L_F(\Gamma_1)$ is the identity map.

 We prove by induction on $|p|+|q|$ that $f(pq^*) = p'q'^*$ for {\em all} non-zero elements of $LI(\Delta_1)$. We may assume that $pq^*$ is not a natural basis element of $LI(\Delta_1)$. It follows from Theorem~\ref{leavittinv} that $|p|+|q|$ is at least 2 since any element $pq^*$ with $|p|+|q| < 2$ is a natural basis element. If $|p| + |q| = 2$, then again by Theorem~\ref{leavittinv}, $pq^*$ is $\gamma \gamma^*$ for some special edge $\gamma$ of $\Delta_1$. If $\gamma = e$, then $\gamma \gamma^* = s(e)$; if $\gamma \neq e$, then $\gamma \gamma^*$ is in $LI(\Gamma_1)$. In both cases, $f(\gamma \gamma^*) = \gamma' \gamma'^*$. Hence the result holds if $|p|+|q| \leq 2$. This is a basis for the induction.

 Now suppose that $|p|+|q| > 2$ and that $f(wz^*) = w'z'^*$ for all non-zero elements $wz^*$ of $LI(\Delta_1)$ such that $|w|+|z| < |p|+|q|$. We aim to show that $f(pq^*) = p'q'^*$. We may assume that the edge $e$ is the first edge of either $p$ or $q$ or both and that $pq^*$ is {\em not} a natural basis element of $L_F(\Delta_1)$, so $p'q'^*$ is not a natural basis element of $L_F(\Gamma_1)$. Then $p'q'^* = xe_1e_1^*y^*$ for some directed paths $x$ and $y$ in $\Gamma_1$, where $e_1$ is the special edge $e_1 = \gamma(s(e_1))$. If $s^{-1}(s(e_1)) = \{e_1,e_2,...e_n\}$ then $e_1e_1^* = s(e_1) - e_2e_2^* - \ldots - e_ne_n^*$, so $p'q'^* = xy^* - xe_2e_2^*y^* - \ldots \ - xe_ne_n^*y^*$. Consider the case where $e$ is the first letter of $p$ but not the first letter of $q$, so $p = ep'$ and $q = q'$.  We have $pq^* = ep'q'^* = exy^* - exe_2e_2^*y^* - \ldots - exe_ne_n^*y^*$. So $f(pq^*) = f(exy^*) - f(exe_2e_2^*y^*) - \ldots - f(exe_ne_n^*y^*)$. Since the terms $exe_ie_i^*$ are basis elements of $L_F(\Delta_1)$, we have $f(exe_ie_i^*y^*) = xe_ie_i^*y^*$ for $i = 2,...,n$. By the induction hypothesis, $f(exy^*) = xy^*$ since $|ex|+|y| < |exe_1|+|e_1y|$. Hence $f(pq^*) = xy^* - xe_2e_2^*y^* - \ldots - xe_ne_n^*y^* = p'q'^*$. A similar argument shows that $f(pq^*) = p'q'^*$ if $e$ is the first letter of $q$ but not the first letter of $p$ or if $e$ is the first letter of both $p$ and $q$, as required.

Now suppose that $p_1q_1^*$ and $p_2q_2^*$ are non-zero elements of $L_F(\Delta_1)$ such that $p_1q_1^*p_2q_2^* \neq 0$. Assume without loss of generality that $p_2$ is a prefix of $q_1$ (the other case is dual), so $q_1 = p_2t$ for some directed path $t$ in $\Delta_1$ which leads to $f(p_1q_1^*p_2q_2^*) = f(p_1(q_2 t)^*) = p'_1 (q_2 t)'^*$. If $p_2$ starts with $e$, then so does $q_1$ and hence $f(p_1q_1^*p_2q_2^*) = p'_1 (q_2 t)'^* = p'_1 t^* q'^*_2 = p'_1 t^* p'^*_2 p'_2 q'^*_2 = p'_1 q'^*_1 p'_2 q'^*_2 = f(p_1 q^*_1) f(p_2 q^*_2)$. If $p_2$ does not start with $e$ but $q_1$ does, then $p_2, q_2$ must be trivial and $p_2 = q_2 = s(e)$. So we have  $f(p_1q_1^*p_2q_2^*) = p'_1 (q_2 t)'^* = p'_1 t'^* = p'_1 q'^*_1 p'_2 q'^*_2 = f(p_1 q^*_1) f(p_2 q^*_2)$. If $q_1$ does not start with $e$, then neither does $p_2$. So we have $f(p_1q_1^*p_2q_2^*) = p'_1 (q_2 t)'^* = p'_1 t^* q'^* = p'_1 t^* p^*_2 p_2 q'^*_2 = p'_1 q'^*_1 p'_2 q'^*_2 = f(p_1 q^*_1) f(p_2 q^*_2)$. In summary, $f$ restricts to a 0-morphism from $LI(\Delta_1)$ to $LI(\Gamma_1)$ so that $f$ is a 0-retraction of $L_F(\Delta_1)$ onto $L_F(\Gamma_1)$.

If $\Delta_1 = \Delta$ we stop.  If $\Delta_1 \neq \Delta$, we may attach another NE-edge to $\Delta_1$ to obtain another subgraph $\Delta_2$ of $\Delta$ and proceed as before, obtaining a $0$-retraction $h$  of $L_F(\Delta_2)$ onto $L_F(\Delta_1)$. Then $hf$ is a $0$-retraction of $L_F(\Delta_2)$ onto $L_F(\Gamma_1)$. Continue adding NE paths to intermediate graphs until we eventually reach $\Delta$ and obtain a $0$-retraction of $L_F(\Delta)$ onto $L_F(\Gamma_1)$. This completes the proof since  $M_n(L_F(\bar{\Gamma})) \cong L_F(\Delta)$ and $L_F(\Gamma) \cong L_F(\Gamma_1)$. \QED

\section{Addendum: The kernel of the map from $I(\Gamma)$ onto $LI(\Gamma)$}

In this addendum, we show that the kernel of the map from $I(\Gamma)$ onto $LI(\Gamma)$ is the congruence $\leftrightarrow$ introduced by Lenz \cite{lenz}, answering a question raised by Milan based on an earlier version of this paper. We first recall the definition of Lenz's congruence $\leftrightarrow$.

Let $S$ be an inverse semigroup (with $0$) and for each $a \in S$ let $a^{\downarrow} = \{x \in S : x \leq a$ in the natural partial order on $S\}$. Given $a,b \in S$ we define $a \rightarrow b$ if, whenever $0 < x \leq a$, $a^{\downarrow} \cap b^{\downarrow} \neq \{0\}$, and define $a \leftrightarrow b$ if $a \rightarrow b$ and $b \rightarrow a$. Then $\leftrightarrow $ is a $0$-restricted congruence on $S$ (\cite{LMS}, Proposition 3.4). The congruence $\leftrightarrow$ was introduced by Lenz \cite{lenz} in connection with his construction of various topological groupoids associated with inverse semigroups: it has been studied by several authors, in particular by Lalonde, Milan and Scott \cite{LMS} who used it to study the ideal structure of the tight $C^*$-algebra of an inverse semigroup.

\begin{Theorem}
\label{lenz}

For any graph $\Gamma$, the congruence $\leftrightarrow$ is the kernel of the natural homomorphism from $I(\Gamma)$ onto $LI(\Gamma)$: that is $LI(\Gamma) \cong I(\Gamma)/{\leftrightarrow}$.

\end{Theorem}

\noindent {\bf Proof}. From the definition of the Leavitt inverse semigroup $LI(\Gamma)$ it is clear that the kernel of the natural homomorphism from $I(\Gamma)$ onto $LI(\Gamma)$ is the congruence $\rho$ generated by $\{(ee^*,s(e)): s(e)$ has out-degree $1$ in $\Gamma^0\}$. We aim to show that $\rho \,  = \,\, \leftrightarrow$.

Suppose first that $e$ is an edge of $\Gamma$ such that $s(e)$ has out-degree $1$. Then $s(e)\, \rho \, ee^*$. Let $0 < x \leq s(e)$ in $I(\Gamma)$. Then $x = pp^*$ for some directed path $p$ with $s(p) = s(e)$. Either $p = s(e)$ or $p = eq$ for some directed path $q$ with $s(q) = r(e)$. Clearly $s(e)^{\downarrow} \cap (ee^*)^{\downarrow} \neq \{0\}$. Also, if $p = eq$ then $(pp^{*})^{\downarrow} = \{eqtt^*q^*e^* : s(t) = r(q)\} \cap (ee^{*})^{\downarrow} \neq \{0\}$, so $s(e) \rightarrow ee^*$. By essentially the same argument, $ee^* \rightarrow s(e)$. Hence $s(e) \leftrightarrow ee^*$ and so
$\rho \, \subseteq \, \leftrightarrow$.

Now suppose that $p_1q_1^*, p_2q_2^*$ are non-zero elements of $I(\Gamma)$ such that $p_1q_1^* \leftrightarrow p_2q_2^*$. Then $(p_1q_1^{*})^{\downarrow} \, \cap \, (p_2q_2^{*})^{\downarrow} \neq \{0\}$, so there exist paths $t_1,t_2$ such that $p_1t_1t_1^*q_1^* = p_2t_2t_2^*q_2^*$. So either $p_1$ is a prefix of $p_2$ or $p_2$ is a prefix of $p_1$. Suppose without loss of generality that $p_2 = p_1z$ for some path $z = e_1e_2...e_n$. We claim that $z$ is an NE path. If not, there is some index $i$ with $1 \leq i \leq n$ such that $s(e_i)$ has out-degree at least $2$, and so there is some edge $f$ with $f \neq e_i$ and $s(f) = s(e_i)$. Let $z_1 = e_1...e_{i-1}$ and $z_2 = e_i...e_n$. Then $(p_1z_1ff^*z_1^*q_1^{*})^{\downarrow} \, \cap \, (p_2q_2^{*})^{\downarrow} \neq \{0\}$. Hence there exist paths $t_3,t_4$ such that $p_1z_1ft_3t_3^*f^*z_1^*q_1^* = p_2t_4t_4^*q_2^* = p_1z_1z_2t_4t_4^*q_2^*$. This implies that $z_2t_4 = ft_3$. But this is impossible since the first edge of $z_2$ is $e_i \neq f$. Hence $z$ is an NE path. Now we have $p_1t_1t_1^*q_1^* = p_2t_2t_2^*q_2^* = p_1zt_2t_2^*q_2^*$, so $t_1 = zt_2$ and hence $q_2t_2 = q_1t_1 = q_1zt_2$. This implies that $q_2 = q_1z$. So we have $p_2 = p_1z$ and $q_2 = q_1z$ for some NE path $z$. But this means that $(p_1q_1^*) \, \rho \, (p_2q_2^*)$. Hence $\leftrightarrow \, \subseteq \, \rho$. It follows that the congruences $\leftrightarrow $ and $\rho$ coincide. \QED

\medskip

\noindent {\bf Acknowledgements} The authors thank Ruy Exel, David Milan,
Benjamin Steinberg and Efim Zelmanov for discussions concerning
various aspects of this paper.   The second author is grateful to the
Mathematics Department of UNL (the University of Nebraska - Lincoln)
for the excellent working conditions and to the China Scholarship
Council for the financial support to enable him to visit UNL.

\bigskip

\noindent John Meakin \,\,  jmeakin@unl.edu

\noindent Department of Mathematics, University of Nebraska, Lincoln, Nebraska 68588, USA.

\medskip

\noindent Zhengpan Wang \,\, zpwang@swu.edu.cn

\noindent School of Mathematics and Statistics, Southwest University, Chongqing 400715, China.

\end{document}